\documentclass[11pt,a4paper,reqno]{amsart}
\usepackage{amsmath,amssymb,amsfonts,epsfig,mathrsfs,cite, hyperref}
\usepackage[T1]{fontenc}
\usepackage{color}
\usepackage{array}
\usepackage{amsthm}
\usepackage{amstext}
\usepackage{graphicx}
\usepackage{setspace}

\usepackage[title]{appendix}

\usepackage{booktabs}

\usepackage{tcolorbox}

\usepackage{float}

\makeatletter
\@namedef{subjclassname@2020}{%
  \textup{2020} Mathematics Subject Classification}
\makeatother

\usepackage[margin=2.5cm]{geometry}
\usepackage{color}
\usepackage{enumitem}
\usepackage{amscd,psfrag}
\usepackage{yhmath}
\usepackage[mathscr]{eucal}

\usepackage{comment}

\allowdisplaybreaks[4]

\usepackage{slashed}

\usepackage{stackengine}

\makeatletter
\pdfpageheight\paperheight
\pdfpagewidth\paperwidth

\usepackage{epstopdf}

\usepackage{indentfirst}	

\usepackage[normalem]{ulem}
\theoremstyle{plain}

\newtheorem{definition}{Definition}
\newtheorem{theorem}[definition]{Theorem}
\newtheorem*{theorem*}{Theorem}

\newtheorem{remark}[definition]{Remark}

\newtheorem*{remark*}{Remark}
\newtheorem*{sideremark*}{Side Remark}

\newtheorem*{claim*}{Claim}
\newtheorem*{lemma*}{Lemma}
\newtheorem*{q*}{Question}
\newtheorem{lemma}[definition]{Lemma}
\newtheorem{corollary}[definition]{Corollary}
\newtheorem*{corollary*}{Corollary}

\newtheorem{proposition}[definition]{Proposition}

\newcommand{\R}{\mathbb{R}}

\newcommand{\na}{\nabla}

\newcommand{\p}{\partial}
\newcommand{\loc}{{\rm loc}}
\newcommand{\weak}{\rightharpoonup}
\newcommand{\e}{\varepsilon}
\newcommand{\C}{\mathbb{C}}
\newcommand{\dd}{{\rm d}}

\newcommand{\G}{\Gamma}

\newcommand{\M}{{\mathcal{M}}}
\newcommand{\bra}{\left\langle}
\newcommand{\ket}{\right\rangle}
\newcommand{\CC}{{\mathscr{C}}}

\newcommand{\1}{{\mathbf{1}}}

\newcommand{\mres}{\mathbin{\vrule height 1.6ex depth 0pt width
0.13ex\vrule height 0.13ex depth 0pt width 1.3ex}}

\newcommand{\weakstar}{\mathrel{\ensurestackMath{\stackon[1pt]{\rightharpoonup}{\scriptstyle\ast}}}}

\newcommand{\kbs}{{\mathcal{K}}_{\rm BS}}

\newcommand{\glue}{u^{\heartsuit}}

\def\XXint#1#2#3{{\setbox0=\hbox{$#1{#2#3}{\int}$ }
\vcenter{\hbox{$#2#3$ }}\kern-.6\wd0}}

\newcommand{\zz}[1]{{u^{\mathrm{good}}_{\varepsilon_{#1}}}}

\newcommand{\leb}{{\mathcal{L}}}

\newcommand{\B}{{\mathbf{B}}}

\newcommand{\mcomp}{\M_{\rm comp}}

\newcommand{\curl}{{\rm curl}}

\newcommand{\srad}{\sigma_{\rm rad}}
\newcommand{\sym}{{\bf Sym}_{2 \times 2}(\R)}

\newcommand{\gdata}{{\omega^{\rm good}_{0,\e}}}

\newcommand{\gdataa}{{\omega^{\rm good}_{0,\eta}}}

\newcommand{\gsol}{{\omega^{\rm good}_{\e}}}

\newcommand{\ugood}{u^{\rm good}_\e}
\newcommand{\ugoodd}{u^{\rm good}_\eta}

\newcommand{\dee}{{\mathscr{DR}_{\e,\eta}}}

\newcommand{\gsoll}{{\omega^{\rm good}_{\eta}}}

\newcommand{\rinit}{{\mathscr{R}^{\rm initial}_{\e,\eta}}}

\newcommand{\gdatan}{\omega^{\mathrm{good}}_{0,\e_n}}

\newcommand{\rey}{{\mathscr{R}^\heartsuit}}

\newcommand{\crey}{{c_{\rm Rey}}}

\newcommand{\leray}{{\mathbf{P}}}

\newcommand{\lerayy}{\leray\circ{\operatorname{div}}}

\newcommand{\cen}{\mathtt{c}}

\newcommand{\flf}{\underline{\mathtt{f}}}

\newcommand{\corner}{\mathsf{v}}

\newcommand{\ci}{{\corner_i}}

\newcommand{\stripeij}{\mathsf{x}_{ij}}

\newcommand{\bigtime}{{\mathsf{T}_{\mathrm{cell}}}}

\newcommand{\band}{{\mathsf{Band}}}

\newcommand{\timeint}{{\left[t^*_{m,i,j},\,t^*_{m,i,j}+d_{ij}^{(m)}\right]}}

\newcommand{\vcorr}{V^{\mathrm{corr}}}

\newcommand{\biga}{{\mathbf{A}^{(m)}}}

\title[Delort's general vortex sheet problem with signed vorticity measure]{Global weak solutions to the two-dimensional Euler
equations with initial vorticity of arbitrary sign}
\author{Siran Li}

\address{Siran Li: School of Mathematical Sciences $\&$ CMA-Shanghai, Shanghai Jiao Tong University, No.~6 Science Buildings,
800 Dongchuan Road, Minhang District, Shanghai, China (200240)}

\email{\texttt{siran.li@sjtu.edu.cn}} 

\keywords{Euler equations; vorticity; vortex sheet; generalised vortex sheet problem; global existence}

\subjclass[2020]{Primary: 76B47; Secondary: 35Q31, 35D30, 35A01}
\date{\today}

\begin{document}

\begin{abstract}
We address the general vortex sheet problem for the two-dimensional incompressible Euler equations: Given a compactly supported signed Radon measure $\omega_0$ in $H^{-1}(\mathbb{R}^2)$, does there exist a global classical weak solution to the Euler equations whose initial vorticity is $\omega_0$? It has been shown by J.-M. Delort (1991) that the answer is affirmative when $\omega_0$ has a distinguished sign. In this paper, without imposing any sign condition on $\omega_0$, we construct a global Euler flow with initial vorticity $\omega_0$ that is eventually $\CC^\infty$. Key ingredients of the proof are a temporal glueing argument for approximate solutions and a convex integration scheme that effectively reduces the $L^1$-norm of the Reynolds stress tensor in each iteration stage. This convex integration scheme, based on Lamb--Chaplygin dipoles traversing straight line segments, is an adaptation of the construction by Bru\`e--Colombo--Kumar (2026). New features are introduced to control the $L^1$-norm of the vorticity and to go beyond the spatially periodic setting. It is also proved that the vorticity is not renormalised by the constructed Euler flow.

\end{abstract}

\maketitle

\tableofcontents

\section{Introduction}\label{sec: intro}

\subsection{General Vortex Sheet Problem}

This paper is concerned with the Euler equations for incompressible flows in two spatial dimensions:
\begin{equation}\label{Euler}
\begin{cases}
\p_t u + {\rm div}(u \otimes u) + \na p = 0,\\
{\rm div}\, u = 0 \qquad\qquad  \text{ in } [0,\infty[ \times \R^2,
\end{cases}
\end{equation}
where $u=u(t,x): [0,\infty[ \times \R^2 \to \R^2$ is the velocity field and $p=p(t,x): [0,\infty[ \times \R^2 \to \R$ is the pressure field of the fluid. These equations, first proposed by Euler in 1757,  have been central to fluid mechanics and PDE (Partial Differential Equations) ever since. We refer the reader to Chemin~\cite{chemin}, Lions~\cite{lions}, and Majda--Bertozzi~\cite{MajdaBertozzi}, among many other classical texts, for a comprehensive survey on the developments of the mathematical theory for the Euler equations.

Cauchy's memoir on water waves (written in 1815 and published in 1827) gave an early account of the vorticity evolution for the Euler equations; \emph{cf}.~\cite{Cauchy1827,FrischVillone2014}.
This theory was later recast into the geometric framework by Helmholtz~\cite{Helmholtz1858}: his vortex theorem ascertains that the vortex lines in an ideal fluid are carried by the flow.  Hankel in 1861~\cite{Hankel1861} derived the conservation of circulation around material loops based on the aforementioned works by Cauchy, and Lord Kelvin independently established the same law  in 1869 in his memoir on vortex motion~\cite[\S~60(a)]{Thomson1869}. For 2D Euler flows, vorticity is a scalar field transported by the
velocity as long as the flow remains smooth. The singular counterpart goes back to Helmholtz's study of discontinuous
motion~\cite{Helmholtz1868}.

In modern terms, the singular vortex structure analysed in~\cite{Helmholtz1868} is known as a \emph{vortex sheet}: across an interface, the tangential components of the velocity fields become discontinuous, giving rise to the concentration of the vorticity, \emph{i.e.}, the distributional curl of the velocity $$\omega := \curl\, u \equiv \p_1u^2 - \p_2 u^1,$$ on the interface. This is a ubiquitous phenomenon in nature. For a vortex sheet carried by a smooth curve $\G_0 \subset \R^2$, its vorticity takes the form
\begin{align*}
\omega = \gamma\,\mathcal H^1\mres{\Gamma_0},
\end{align*}
where $\mathcal{H}^1$ is the 1D Hausdorff measure, and the vortex sheet strength function $\gamma$ may change sign. Classical works by Prandtl, Kaden, and Pullin~\cite{PrandtlSpiral,Kaden,Pullin1978}, among others, made the first attempts to provide rigorous  mathematical characterisations beyond the phenomenological descriptions of the roll-up of vortex sheets, and they motivated the studies in mathematical hydrodynamics of vorticity that is singular with respect to the area
measure. When a vortex sheet persists as a suitably regular curve, its evolution is described by the Birkhoff--Rott equation which, loosely speaking, is equivalent to the weak formulation of the 2D Euler system; see, \emph{e.g.}, Lopes Filho--Nussenzveig Lopes--Schochet~\cite{LNS2007}.

If, however, the vortex sheet becomes irregular and loses the curve structure in later time, one cannot study the geometrical evolution of the sheet itself. Instead, the starting point for the  analysis of vorticity is the weak formulation for the Euler system~\eqref{Euler}. The following definition has been proposed  by DiPerna and Majda in 1987--88 in the seminal papers~\cite{dm2,dm3,dm1}. 

\begin{definition}\label{def: weak sol}
A function $u: [0,\infty[\times \R^2 \to \R^2$ is a classical global weak solution to the Euler system~\eqref{Euler} if the following holds:
\begin{itemize}
\item
$u \in L^\infty_{\mathrm{loc}}\left([0,\infty[;L^2_{\mathrm{loc}}(\R^2)\right)$;
\item
$u$ satisfies~\eqref{Euler} in the sense of distributions: for any test vector field \(\phi\in\CC^\infty_c([0,\infty)\times\R^2;\R^2)\) with
\(\operatorname{div}_x\phi=0\), it holds that 
\begin{align*}
 \int_0^\infty\!\!\int_{\R^2}
     \bigl(u\cdot\p_t\phi+(u\otimes u):\na_x\phi\bigr)\,\dd x\,\dd t
       +\int_{\R^2}u_0(x)\cdot\phi(0,x)\,\dd x=0;
\end{align*}

\item
$u \in \operatorname{Lip}_\loc\left([0,\infty[;H^{-L}_{\mathrm{loc}}(\R^2)\right)$ for some $L>0$.
\end{itemize}
\end{definition}

A significant problem in the above framework, which is also the central theme of this paper, is known as the \emph{general vortex sheet problem}. It asks about the existence of a classical global weak solution to~\eqref{Euler} whose vorticity  remains bounded as a signed measure. More precisely, consider the initial vorticity as a compactly supported signed Radon measure:
\begin{equation}\label{initial vorticity}
 \omega_0\in\mcomp(\R^2)\cap H^{-1}(\R^2).
\end{equation}
The $H^{-1}$-regularity ensures that $\omega_0$ is atomless (\emph{a.k.a.} diffuse; \emph{i.e.}, $\omega_0 (\{x\}) = 0$ for any $x \in \R^2$); see \cite[Lemme~1.2.5]{Delort}. Denote the Biot--Savart kernel as
\begin{equation}\label{eq: KBS}
\kbs(x):=\frac{x^\perp}{2\pi|x|^2}.
\end{equation}
We then seek classical global weak solutions as in Definition~\ref{def: weak sol} with the initial datum
\begin{equation}\label{u0, id}
u\big|_{\{t=0\}} = \kbs\star\omega_0,
\end{equation}
 where $\star$ denotes the convolution.

To put things into perspective, the general vortex sheet problem concerns the existence of solutions to the 2D Euler equations at a substantially lower regularity than the classical well-posedness theory of
Yudovich~\cite{Yudovich1963}, which requires $\omega \in L^\infty([0,\infty[\times\R^2)$. DiPerna--Majda~\cite{dm2,dm3,dm1} investigated the effects of the concentration
of kinetic energy on the passage to the limit in regularisations $\{u_\e\}$ of the 2D Euler equations; \emph{e.g.}, vanishing viscosity limit or large-eddy simulation. The key issue here is about \emph{weak continuity}: the 2D Euler system~\eqref{Euler} admits cancellations of certain  concentrations in the weak momentum equation, even in the absence of strong velocity convergence. A uniform bound on the total variation of the vorticity and a local kinetic energy bound would imply weak compactness for the velocities $u_\e \weak \overline{u}$, but cannot identify the weak limit of the quadratic term; \emph{i.e.}, one may have ${\rm div}(u_\e \otimes u_\e) \not\rightarrow  {\rm div}\left(\overline{u} \otimes \overline{u}\right).$

\subsection{Vortex sheet with vorticity of distinguished sign \emph{\`{a} la} Delort}

In 1991,  Delort~\cite{Delort}  established the global existence of classical weak solutions to~\eqref{Euler}, under the assumption that
\begin{align*}
\omega_0 = \omega_0' + \omega_0''
\end{align*}
where
\begin{equation*}
\begin{cases}
\omega_0' \in \mcomp(\R^2)\cap H^{-1}(\R^2), \qquad \omega_0' \geq 0 \text{ as measures}; \\
\omega_0'' \in L^p_{\rm comp}(\R^2)\qquad\text{for any } p \geq 1.
\end{cases}
\end{equation*}
The key assumption here is that $\omega_0'$ has a distinguished sign. Vecchi--Wu~\cite{VW} extended this to the endpoint case $p=1$ via the Dunford--Pettis theorem. Evans--M\"uller~\cite{EvansMuller} (see also Semmes~\cite{Semmes1994}) gave a harmonic analytic
proof of Delort's theorem by controlling specific
quadratic expressions in the  local Hardy space. Lopes Filho--Nussenzveig Lopes--Xin~\cite{Reflection} further extended the global existence theorem to the case of ``dipoles'', namely a distinguished-sign vorticity measure in a half-plane together with its odd reflection across the boundary.

To understand the role played by the distinguished sign of the initial vorticity component $\omega_0'$, let us first follow the arguments of Delort~\cite{Delort}; see also  Schochet~\cite{Schochet1995}. It is well known that for smooth 2D Euler flows, the vorticity is transported by the velocity. Nevertheless, for vortex sheets the product $u\omega$ may fail to be well defined \emph{per se}, but the symmetries of the transport term allow $u \omega$ to be defined in the sense of distributions. Indeed,  since the kernel $\kbs$ is odd, for $\omega_\e =\omega_\e (t,x)\in \CC^\infty_c([0,\infty[ \times \R^2)$ and $\varphi=\varphi(x) \in \CC^\infty_c(\R^2)$, it holds that  
\begin{align}\label{schochet, kernel}
 \int_{\R^2}\omega_\e(t,x)u_\e(t,x)\cdot\na\varphi(x)\,\dd x
&=
 \frac12\iint_{\R^2\times\R^2}
       \kbs(x-y)\cdot\bigl(\na\varphi(x)-\na\varphi(y)\bigr) \omega_\e(t,x)\omega_\e(t,y)\,\dd x\,\dd y \nonumber\\
       &=: \iint_{\R^2\times\R^2}
       H_\varphi(x, y) \omega_\e(t,x)\omega_\e(t,y)\,\dd x\,\dd y.
\end{align}
Here $\omega_\e := \curl u_\e$, and $\{u_\e\}$ is a smooth approximation of $u$ as above. The identity~\eqref{schochet, kernel} is referred to as the \emph{weak vorticity formulation} of the 2D Euler system~\eqref{Euler}. The term $\na\varphi(x)-\na\varphi(y)$ on the right-hand side  cancels the
$\mathcal{O}\left(|x-y|^{-1}\right)$-singularity of \(\kbs\), so the kernel $H_\varphi(x,y)$ is symmetric and bounded near the diagonal $\{x=y\}$. Nonetheless, $H_\varphi(x,y)$ has no continuous extension across the diagonal in general. To pass to the limit in~\eqref{schochet, kernel}, one
must account for the part of $\omega_\e \otimes \omega_\e$ that lies near the diagonal, and a key observation by Delort~\cite{Delort} is that the strong $L^1$-convergence of $u_\e^1 u_\e^2$ would suffice. This is achieved for $\omega_0'$ possessing a distinguished sign by Delort \cite[Lemme~1.2.7]{Delort}, and has gradually acquired further understandings in subsequent developments, \emph{e.g.}, the maximal vorticity function by Majda~\cite{majda}, the local Hardy space estimates by Evans--M\"{u}ller~\cite{EvansMuller} and Semmes~\cite{Semmes1994}, and the weak vorticity approach to  the point vortex by Schochet~\cite{Schochet1996}, among many other works.

Another approach seeks function spaces that
quantify vorticity concentration. Tadmor~\cite{Tadmor} introduced a scale of spaces for compactness in \(H^{-1}\), now known as the Riesz--Morrey--Tadmor spaces. This approach has been refined by Dom\'{i}nguez--Milman~\cite{DM} through quantitative sparseness criteria. Dom\'{i}nguez--Spector~\cite{DS} constructed steady mixed-sign examples of vorticity uniformly bounded
in the critical logarithmic Morrey space \(M^{\log,1}\). Their velocities converge weakly to zero but fail to converge strongly
in \(L^2_\loc\). This illustrates the distinction
between concentration cancellation and absence of concentration.

Let us refer the reader to, \emph{e.g.}, Marchioro--Pulvirenti~\cite{mp} for a modern, comprehensive account of the vortex dynamics. We also mention the celebrated work~\cite{BourgainLi2015} by Bourgain and Li, in which the strong ill-posedness of 2D Euler equations with velocity in the critical space $H^2$. The arguments in~\cite{BourgainLi2015} feature a delicate analysis of the vorticity from the Lagrangian perspective. 

\subsection{Convex integration and temporal glueing}

In this paper, we resolve the general vortex sheet problem for the 2D Euler flows without any sign condition on the initial vorticity (Corollary~\ref{cor: general vortex sheet}). The key ingredients of the proof are convex integration, especially the Lamb--Chaplygin dipole-based scheme developed by Bru\`{e}--Colombo--Kumar~\cite{BCK}, and in-time glueing.

We first give a brief survey on the convex integration techniques in fluid PDEs, which is by no means exhaustive. For more comprehensive treatment, we refer the reader to  \cite{DLSHPrinciple,BuckmasterVicolSurvey}.

 In 1954, Nash  developed this technique in his seminal work~\cite{nash} solving the $\CC^1$-isometric embedding problem, which has also motivated Gromov~\cite{gromov} to develop the ground-breaking theory of partial differential relations and $h$-principle. In 2009, De Lellis--Sz\'{e}kelyhidi~\cite{DLS2009} first 
found applications of this geometrically originated theory in the analysis of PDE and hydrodynamics. By recasting the incompressible Euler equations in $\R^d$ for any dimension $d$ as a differential inclusion, De Lellis--Sz\'{e}kelyhidi constructed nonunique, ``wild'' weak solutions to the Euler equations. Subsequently, Wiedemann~\cite{Wiedemann2011} further obtained global weak
Euler solutions with arbitrary solenoidal \(L^2\)-initial velocity on the torus (but without any characterisations for the vorticity, in contrast to the general vortex sheet problem). The construction of continuous dissipative Euler flows by De Lellis
and Sz\'{e}kelyhidi~\cite{DLSContinuous} led to successive improvements
in the regularity of such solutions, culminating in the final resolution of Onsager's conjecture (1949)~\cite{o}: For the 3D incompressible Euler equations, if a weak solution has regularity higher than the Hölder exponent~$1/3$, then the kinetic energy is conserved; if its regularity is lower than~$1/3$, then nonconservative weak solutions exist. The first half the conjecture was proved by Constantin--E--Titi~\cite{CET}, and the second half was established by Isett~\cite{IsettOnsager} using convex integration; \emph{cf.} also Buckmaster--De Lellis--Sz\'{e}kelyhidi--Vicol~\cite{BDLSV}. Later, Buckmaster--Vicol~\cite{BuckmasterVicol2019} 
further extended the convex integration technique to prove the nonuniqueness of
finite-energy weak solutions to the three-dimensional
Navier--Stokes equations.

Let us also discuss the glueing constructions in fluid PDEs. Indeed, the glueing of Euler flows over short time intervals has been used in Isett's proof of the
Onsager conjecture~\cite{IsettOnsager} and in the work of Buckmaster--Colombo--Vicol~\cite{BCV}. As for the Navier--Stokes equations, Buckmaster--Vicol~\cite{BuckmasterVicol2019} introduced intermittent building blocks to prove nonuniqueness of
finite-energy weak solutions in 3D, and Cheskidov--Luo~\cite{CheskidovLuoSharp}
used temporal intermittency and time correctors to obtain sharp nonuniqueness results in the periodic setting, and developed accelerating jets to approach the 2D
\(L^2\)-velocity threshold from below~\cite{CheskidovLuo2D}. An alternative approach through instability, rather than through
convex integration, gives nonuniqueness of Leray--Hopf weak solutions with
a prescribed body force in~\cite{ABC} by Albritton--Bru\`{e}--Colombo. See Colombo's ICM talk~\cite{ColomboICM2026}
for a discussion of the relation between these approaches.

In regard to the vortex sheet evolutions, Sz\'ekelyhidi~\cite{SzekelyhidiVortexSheet} used convex integration
to construct nonunique admissible solutions emanating from a flat vortex sheet. Note, however, that the evolving vorticity in~\cite{SzekelyhidiVortexSheet}  is not a bounded Radon measure. Mengual--Sz\'ekelyhidi~\cite{MS} treated vortex sheet initial data concentrated on suitably regular simple closed curves without a sign restriction on the sheet strength, obtaining solutions that are smooth outside a region that expands from the initial vortex sheet. Mengual~\cite{MengualLp} proved nonuniqueness for admissible weak solutions for
particular compact vortex data with \(L^p\)-initial vorticity for each \(2<p<\infty\). Recently, Cheskidov--Zeng--Zhang~\cite{CZZEuler} used a  backward convex integration scheme to obtain decreasing energy profiles and prescribed
initial velocities for selected classes of periodic data. In all these works above, a condition on the \emph{initial} vorticity does not lead to a bound on the total variation of vorticity
at later times.

On the other hand, for regularity of the \emph{evolving} vorticity, Vishik proved nonuniqueness of Euler solutions with \(L^p\)-vorticity for any $p<\infty$, in the presence of an external forcing in the Euler equations; see the exposition
in~\cite{AlbrittonEtAlVishik}. For the unforced equations,
Bru\`e--Colombo~\cite{BC} constructed periodic nonunique
solutions with weak \(L^1\)-vorticity, and Buck--Modena~\cite{BM} constructed spatially compactly supported flows with prescribed energy profiles and vorticity in the Hardy space $H^p(\R^2)$ for \(0<p<1\). The recent important paper~\cite{BCK} by Bru\`e--Colombo--Kumar proved the flexibility and nonuniqueness of $L^\infty_t L^2_x$-weak solutions to the Euler equations on the 2D torus  with vorticity in $L^\infty_t L^p_x$ for some $p>1$. The convex integration scheme utilised in our current work is an adaptation of~\cite{BCK}; in particular, the principal building block of the scheme is a rescaled Lamb--Chaplygin dipole, whose amplitude, radius, and location are time-dependent functions chosen to cancel the spatial $L^1$-norm of the Reynolds stress tensor. See~\S\ref{subsec: dipole} and \S\ref{sec: perturbation lemma} for details.

We also bring to the reader's attention two recent preprints, in which 
distinct solutions to the Euler equations on $\R^2$ from the same integrable initial vorticity are constructed. Choi~\cite{ChoiUnbounded2026} treated the case of initial vorticity in \(\bigcap_{1\le p<\infty}L^p(\R^2)\), and Hou--Qin--Song--Zhang~\cite{HQSZ2026} constructed for each fixed \(1<p<\infty\) two vorticities continuous in \(L^1\cap L^p\), transported by their
respective regular Lagrangian flows. These results concern
nonuniqueness of 2D Euler flows with measure-valued initial vorticity $\omega_0$ \emph{absolutely continuous} with respect to the 2D Lebesgue measure. In contrast, for the general vortex sheet problem, $\omega_0$ may admit singular parts.

In passing, we remark that the existence problem for general vortex sheets, especially the weak vorticity formulation in~\eqref{schochet, kernel}, is distinct from a selection principle for solutions by the  vanishing viscosity limit. In fact, Cheskidov--Lopes Filho--Nussenzveig Lopes--Shvydkoy~\cite{CLNS} established energy conservation for periodic
Euler flows with vorticity uniformly in $L^{3/2}$, and
for physically realisable solutions with initial vorticity in $L^p$; \(p>1\). Elgindi--Lopes Filho--Nussenzveig Lopes~\cite{ELNL} proved the absence of anomalous dissipation on the torus for $u_0 \in L^2(\mathbf{T}^2)$ and $\omega_0 = \omega_0' + \omega_0''$ with $\omega_0' \in \M$, $\omega_0' \geq 0$ in $\M$, and $\omega_0'' \in L^1$. De Rosa--Park~\cite{RP} related the defect measure associated with the inviscid limit to defects of velocity compactness and
vorticity concentration. Moreover, Crippa--Nobili--Seis--Spirito~\cite{CNSS} established the renormalisation property for $L^1$-vorticity obtained via vanishing viscosity under certain hypotheses on the initial data. We do not claim any results about vanishing viscosity limits or renormalisation in this paper.

\subsection{Main results}

The main theorem of this paper is the following:

\begin{theorem}\label{mainthm}
Given any $\omega_0\in\mcomp(\R^2)\cap H^{-1}(\R^2)$ and $\e_\star>0$. There exists a classical weak solution $u \in L^\infty_{\mathrm{loc}}\left([0,\infty[;L^2_{\mathrm{loc}}(\R^2)\right) \cap \operatorname{Lip}_\loc\left([0,\infty[;H^{-L}_{\mathrm{loc}}(\R^2)\right)$ for any $L > 2$, such that the vorticity $\omega = \curl u$ lies in $\CC\left(]0,\infty[; L^1(\R^2)\right)$ and satisfies
\begin{align*}
&\sup_{t>0}\|\omega(t,\cdot)\|_{L^1(\R^2)}\le |\omega_0|(\R^2)+\e_\star,\\
& \omega(t,\cdot) \weakstar \omega_0(\cdot) \quad \text{and} \quad \|\omega(t,\cdot)\|_{L^1(\R^2)} \longrightarrow |\omega_0|(\R^2)\qquad \text{as } t \searrow 0. 
\end{align*}
In addition, the following characterisations hold:
\begin{itemize}
    \item  After some positive time, $u$ agrees with a smooth global Euler flow.

    \item 
    There exists $\srad$, a smooth radial vector field of circulation $\omega_0(\R^2)$, such that $$u-\srad \in \CC\left([0,\infty[; L^2(\R^2)\right) \cap \CC^1_\loc\left([0,\infty[; H^{-L}(\R^2)\right).$$
\end{itemize} 
\end{theorem}

Theorem~\ref{mainthm} addresses the generalised vortex sheet problem. The vorticity is continuous at any  positive time as a mapping into \(L^1_x\), with an
arbitrarily small excess over the initial total variation in the uniform \(L^1_x\)-bound. The initial condition is attained in the sense of weak-star convergence, together with the convergence of the total variation norms, of the vorticity measures. Moreover, the initial velocity fails to lie in $L^2_x$ in general.

Once we establish Theorem~\ref{mainthm}, the general vortex sheet problem with initial vorticity equal an arbitrary signed measure plus an $L^p$-function can be resolved, \emph{e.g.}, by directly adapting the arguments in Delort~\cite{Delort}. We summarise it as Corollary~\ref{cor: general vortex sheet} below, with the proof safely omitted.

\begin{corollary}\label{cor: general vortex sheet}

Given any initial vorticity $\omega_0 \in H^{-1}(\R^2)$ such that $\omega_0 = \omega_0' + \omega_0''$ with $
\omega_0' \in \mcomp(\R^2)\cap H^{-1}(\R^2)$ and $\omega_0'' \in L^p_{\rm comp}(\R^2)$ for some $p \geq 1$. There exists a global classical weak solution to the Euler system~\eqref{Euler} such that the vorticity $\omega(t,\cdot) \weakstar \omega_0(\cdot)$ in $\M(\R^2)$ as $t \searrow 0$.
\end{corollary}

Our proof of Theorem~\ref{mainthm} also provides a negative answer to the following question quoted verbatim from De Lellis' S\'{e}minaire Bourbaki report~\cite[972-22, Section~6.2]{bourbaki}:

\begin{quote}
A very interesting open question, naturally linked to Euler's equations, is whether
the renormalization property holds for divergence–free fields $b \in L^\infty\big(\R; L^2(\R^2;\R^2)\big)$ when
the vorticity of $b$ is a measure. 
\end{quote}

\begin{corollary}\label{cor: bourbaki}
Let $u = \kbs \star \omega$ be the global classical weak solution to the 2D Euler equations constructed in Theorem~\ref{mainthm}, where the Radon--Nikod\'{y}m derivative $\dd\omega_0/\dd\leb^2 \notin L^\infty(\R^2)$. Then the vorticity $\omega$ is not renormalised on the time interval $[a,b]$ for some $0<a<b$.
\end{corollary}

\subsection{Strategy of the  proof} 
Our proof of Theorem~\ref{mainthm} relies on a convex integration scheme based on the moving dipole construction and the time averaging technique pioneered by Bru\'{e}--Colombo--Kumar~\cite{BCK}. The core of the arguments is to obtain the following:

\begin{lemma}\label{lem: key}
Let $I \subset \R$ be a closed finite interval. Let $(u,p,R)$ satisfy
\begin{align*}
    &\partial_tu+\operatorname{div}(u\otimes u)+\nabla p
       =\operatorname{div}R,\\
    &\operatorname{div}u=0,
\end{align*}
where $u \in \left(\CC^\infty_{t,x} \cap \CC^0_t L^3_x\right) \left(I \times \R^2;\R^2\right)$, $\omega \in \CC^0_t L^1_x\left(I \times \R^2\right)$, and $R \in\left(\CC^\infty_{t,x} \cap\CC^0_t L^1_x\right)\left(I \times \R^2; \sym\right)$. Assume that $u-\srad \in \CC^0_t L^2_x(I \times \R^2;\R^2)$ for some fixed smooth radial vector field $\srad$ and $p(t, 0) \equiv 0$. Also assume that $R(t,\cdot) \equiv 0$ for $t$ in an open neighbourhood of $\p I$, and     that 
\begin{align*}
    \sup_{t \in I} \|R(t,\cdot)\|_{L^1_x} \leq \delta\qquad\text{for some $\delta>0$}.
\end{align*}
Then, for any given $\e_\star>0$, we can find a smooth triple $\left(\widetilde{u}, \widetilde{p}, \widetilde{R}\right)$ satisfying
\begin{align*}
    &\partial_t\widetilde{u}+\operatorname{div}\left(\widetilde{u}\otimes \widetilde{u}\right)+\nabla \widetilde{p}
       =\operatorname{div}\widetilde{R},\\
    &\operatorname{div}\widetilde{u}=0,
\end{align*}
such that for some  universal constant $C_0>0$,
\begin{align}
&\left(\widetilde{u}, \widetilde{p}, \widetilde{R}\right)(t,\cdot) \equiv (u,p,0)(t,\cdot)\quad \text{ for each $t$ in an open neighbourhood of $\p I$}\label{6},\\
&\sup_{t\in I}\left\|\widetilde{R}(t,\cdot)\right\|_{L^1_x} <\e_\star,\label{eq:scan-target-stress}\\
& \sup_{t\in I}\left\|\widetilde{u}(t,\cdot)-u(t,\cdot)\right\|_{L^2_x} +\sup_{t\in I}\left\|\operatorname{curl}\left(\widetilde{u}-u\right)(t,\cdot)\right\|_{L^1_x}
 \leq C_0\sqrt\delta+\e_\star.\label{eq:scan-target-increment}
\end{align} 
Moreover, for each $t>0$, the correction term $\left(\widetilde{u}-u\right)(t,\cdot)$ has zero circulation, and $\curl\left(\widetilde{u}-u\right)(t,\cdot)$ has the vanishing vector first moment and scalar second moment, and is compactly supported in space. The updated triple lies in the same regularity classes as for $(u,p,R)$; that is, $\widetilde{u} \in \left(\CC^\infty_{t,x} \cap \CC^0_t L^3_x\right) \left(I \times \R^2;\R^2\right)$, $\curl\widetilde{u} \in \CC^0_t L^1_x\left(I \times \R^2\right)$, and $\widetilde{R} \in\left(\CC^\infty_{t,x} \cap\CC^0_t L^1_x\right)\left(I \times \R^2; \sym\right)$.  
\end{lemma}

The above perturbation Lemma~\ref{lem: key}, which allows us to decrease the $L^1_x$-norm of the Reynolds stress tensor $R$ on a given time interval $I=\left[t^-_I, t^+_I\right]$ at a controllable cost of the $L^2_x$-norm of velocity plus the $L^1_x$-norm of vorticity. This is achieved by a convex integration scheme, in which  perturbations are added iteratively on some time interval strictly in the interior of $I$ in countably many ``stages'', whose basic building blocks consist of the \emph{Lamb--Chaplygin dipole} as in~\cite{BCK}. The Lamb--Chaplygin dipole~\cite{Chaplygin,Lamb}  is a classical special solution to the Euler equations; see, \emph{e.g.}, 
\cite{MeleshkoVanHeijst} for its history, \cite{AbeChoi} for
stability, and \cite{HuangTong,ChoiJeongSim} for related
travelling vortex-patch configurations.

In the convex integration scheme, the centre of a rescaled version of the Lamb--Chaplygin dipole traverses finitely many stripes in four designated directions in an appropriately chosen ``\emph{time cell}''. The amplitude, radius, and centre position of the dipole are carefully chosen. In Bru\`{e}--Colombo--Kumar~\cite{BCK}, constructions of this type have been exploited to prove the flexibility and nonuniqueness of $L^\infty_t L^2_x$-weak solutions to the Euler equations on $\mathbf{T}^2$ with vorticity in $L^\infty_t L^p_x$; \emph{e.g.}, with $p = 1+(6500)^{-1}$. We inherit from Bru\`{e}--Colombo--Kumar~\cite{BCK} the moving dipole construction with variable radius, the compact representative and auxiliary force, as well as the cancellation mechanism by averages over time cells and a further temporal correction.

In contrast to~\cite{BCK}, we are working over the whole plane $\R^2$, so the constructions in~\cite{BCK} must be essentially modified to be applicable for our purpose. In particular, we forgo the arguments in~\cite[\S 2.3]{BCK} for cutting off the dipole at spatial scale $r^{1/5}$. Periodic recurrence of  trajectories of the moving dipole in~\cite{BCK} is replaced by traversing along finitely many equal-length, equi-spaced stripes in four good directions. The allocation of traversal times, the length of the gaps between adjacent stripes, and the design of ``resets'' for the dipole centres near the endpoints of each stripe all require extra fine-tuning in our construction. In addition, we need to control the $L^1_x$-norm of the vorticity, rather than $L^{p}_x$ for $p>1$ as in~\cite{BCK}. Endpoint estimates for the Reynolds stress tensor, in particular, Lemma~\ref{lemma: IO}, need to be invoked.

By utilising Lemma~\ref{lem: key} on its own, we may only effectively reduce the $L^1_x$ of the vorticity up to any given precision, while we cannot prove the attainability of the initial vorticity $\omega_0$ as $t\searrow 0$. Indeed, in view of~\eqref{6}, this construction only alters the Euler--Reynolds triple $(u,p,R)$ in the strict interior of the given time interval $I$. To justify the initial trace \(\omega_0\), we  shall join smooth Euler flows whose initial vorticities approximate \(\omega_0\) while matching its circulation, first moment, and scalar second moment. The matching allows the difference
of initial velocities to be written as the divergence of a
symmetric tensor, and the glueing places the Reynolds stress in the interiors of carefully chosen time windows accumulating at \(t=0\).

\subsection{Organisation}
The remaining parts of the paper are organised as follows. 

In \S\ref{sec: prelim} we collect the basic notations, constructions, and analytic lemmata.

In \S\ref{sec: perturbation lemma} we prove our key perturbation Lemma~\ref{lem: key} by the moving dipole construction.

In \S\ref{sec: glue} we introduce the temporal glueing  construction, whose output will serve as the initial stage of the convex integration scheme in the following section.

In \S\ref{sec: proof}, we conclude the proof of the Main Theorem~\ref{mainthm} using a convex integration scheme and Lemma~\ref{lem: key}.

Finally, some concluding remarks are presented in \S\ref{sec: concluding remarks}.

\section{Preliminaries}\label{sec: prelim}

\subsection{Notation and nomenclature}

For a vector field $a$ on $\R^2$, we write $a(x) = \left[a^1(x), a^2(x)\right]^\top$ and write $\leray a$ as the Leray projection of $a$; \emph{i.e.}, the $L^2$-orthogonal projection onto the solenoidal part of $a$. Denote a  Euclidean ball as $\B_r(x):=\{y \in \R^2:\, |x-y|<r\}$. Our convention for the divergence of a square matrix field $S$ is $(\operatorname{div}S)_j=\partial_kS_{jk}$; that is, we always take the row-wise divergence. Here and throughout, Einstein's summation convention is adopted. For vector fields $a$ and $b$ on $\R^2$, we take $a\wedge b = a_1b_2-a_2b_1$ and $\curl a = \p_1a^2 - \p_2a^1 \equiv \na^\perp \cdot a$, where $\nabla^\perp=(-\partial_2,\partial_1)$. In particular, the curl of a 2D vector field is a scalar field. We also write $\sym$ for the space of symmetric $2\times 2$ matrices of real entries.

We shall denote by $\kbs(z):=z^\perp/(2\pi|z|^2)$  the Biot--Savart kernel on $\R^2$. The star $\star$ always denotes the convolution operator. Also, ${\rm div}$, $\curl$, and $\star$ are always taken with respect to the spatial variable only, unless otherwise specified.

 For a 2D vector field $a$, its mean is $$\bra a\ket := \int_{\R^2}a(x)\,\dd x,$$ and its torque is $$\operatorname{Torque}(a)
    := \int_{\R^2}
    x \wedge a(x)\,\dd x \equiv \int_{\R^2}
    \bigl(x^1a^2(x)-x^2 a^1(x)\bigr)\,\dd x.$$ For a scalar field $\psi$ on $\R^2$, the (vector) first moment is $\int_{\R^2}x\psi(x)\,\dd x,$ and the (scalar) second moment is the impulse: $\int_{\R^2} |x|^2\psi(x)\,\dd x.$

Let $\mu \in \M(\R^2)$ be a signed Radon measure on the plane. Then $|\mu| \in \M_{+}(\R^2)$ is the total variation measure of $\mu$, and $|\mu|(\R^2)$  total variation norm of $\mu$. When $\mu \ll \leb^2$, namely when $\mu$ is absolutely continuous with respect to the 2D Lebesgue measure, by an abuse of notations we tacitly identify $\dd \mu$ with $\mu\, \dd \leb^2$; that is, we identify $\mu$ with $\dd\mu/\dd\leb^2$.

We write $\CC^0_t L^p_x(I \times \Omega)$, where $I \subset \R$  and $\Omega \subset \R^2$, for the Bochner spaces of spacetime variables; similarly for $L^r_t L^p_x(I \times \Omega)$, $L^r_t {W^{k,p}_x}(I \times \Omega)$ and the like. In such notations, by omitting $\Omega$ we always mean that $\Omega \equiv \R^2$. 

By an \emph{Euler--Reynolds triple} $(u,p,R)$ we mean a $\CC^\infty$-solution to the system
\begin{align*}
    &\partial_tu+\operatorname{div}(u\otimes u)+\nabla p
       =\operatorname{div}R,\\
    &\operatorname{div}u=0,
\end{align*}
where $u$ is a solenoidal vector field, $p$ is a locally integrable scalar field, and $R$ is a field of symmetric $2 \times 2$ matrices.

\subsection{Endpoint $L^1$-estimate for the inverse of matrix divergence}

We now collect a lemma on the inversion of the matrix divergence: One can solve for a \emph{symmetric} matrix $S$ from the equation ${\rm div}\,S=\mathfrak{f}$, provided that $\mathfrak{f}$ satisfies suitable compatibility conditions. The arguments essentially follow from a representation formula for ${\rm div}^{-1}$ in Isett--Oh~\cite[Section~10.3, Proposition~10.1 and Equations~(229)--(231)]{IO}, while we focus on the endpoint $L^1$-estimate for $S$ here. The non-endpoint estimates, namely $L^p$ with $p>1$, are proved in \cite[Appendix]{BCK}. 


\begin{lemma}
\label{lemma: IO}
Let $\mathfrak{f}\in\CC_c^\infty(\B_r(x_0);\R^2)$ be a vector field with zero mean and zero torque:
\begin{align*}
&\int_{\R^2}\mathfrak{f}(x)\,\dd x = \mathbf{0},\\
&\operatorname{Torque}(\mathfrak{f})
    :=\int_{\R^2}
    \bigl(x^1\mathfrak{f}^2(x)-x^2\mathfrak{f}^1(x)\bigr)\,\dd x
    =0.
\end{align*}
There exists a symmetric tensor field 
$S\in\CC_c^\infty(\B_r(x_0);\sym)$ such that
\[
    \operatorname{div}S=\mathfrak{f},
    \qquad
    \|S\|_{L^1_x}\le Cr\|\mathfrak{f}\|_{L^1_x},
\]
where $C$ is a universal constant independent of $r$, $x_0$, and $\mathfrak{f}$.
\end{lemma}

Observe that the zero mean and zero torque conditions are  necessary for the invertibility of the matrix divergence: Given any compactly supported, symmetric tensor field 
$S\in\CC_c^\infty(\B_r(x_0);\sym)$, and integration by parts yields that
\[
 \int_{\R^2}\operatorname{Div}S(x)\,\dd x=0\quad\text{ and } \quad
 \int_{\R^2} x\wedge\operatorname{Div}S(x)\,\dd x=-\int_{\R^2}(S_{21}-S_{12})\,\dd x=0.
\]

\begin{proof}[Proof of Lemma~\ref{lemma: IO}]

Our arguments are adapted from Isett--Oh
\cite[\S~10.3, formulas~(229)--(231) and Proposition~10.1]{IO}, in which an explicit representation formula for the inverse of the divergence of symmetric matrices has been derived. Utilising the time-independent, 2D version of the above reference, we verify below the divergence identity and derive the endpoint estimates.  

\smallskip
\noindent
{\bf Step~1.}  Fix a nonnegative test function $\zeta\in\CC_c^\infty(\B_1(0))$
with $\int_{\R^2}\zeta(x)\,\dd x=1$, and set
\[
    \zeta_r(y)=r^{-2}\zeta\!\left(\frac{y-x_0}{r}\right).
\]
It holds that $\|\zeta_r\|_{L^1}=1$ and $\|\nabla\zeta_r\|_{L^1}=r^{-1}\|\nabla\zeta\|_{L^1}$. For any $y,z\in\R^2$ and $\lambda \in ]0,1[$, write
\[
    v=z-y\qquad\text{and}\qquad X_\lambda=y+\lambda v,
\]
and define the symmetric matrix
\[
    \mathsf{M}(y,z)
    =
    \bigl[-2\zeta_r(y)+\nabla\zeta_r(y)\cdot v\bigr]\left(\frac{v\otimes \mathfrak{f}(z) +  \mathfrak{f}(z)\otimes v}{2}\right)
    -\bigl[\nabla\zeta_r(y)\cdot \mathfrak{f}(z)\bigr]v\otimes v.
\]
The representation formula for $S$ in \cite[\S~10.3, formulas~(229)--(231)]{IO} by Isett--Oh can be expressed as follows:
\begin{align*}
     S(x)=\int_0^1 S_\lambda(x)\,\dd\lambda
\end{align*}
where 
\begin{align}\label{Slambda}
      S_\lambda(x)
    =
    \int_{\R^2}\lambda^{-2}
    \mathsf{M}\!\left(y,y+\frac{x-y}{\lambda}\right)\,\dd y.
\end{align}
Notice here that for fixed $y$ and $\lambda$, the map $z\mapsto X_\lambda$ has Jacobian determinant $\lambda^2$. 

Where $\mathsf{M}(y,z)$ is nonzero, both $y$ and $z$ belong to
$\B_r(x_0)$, and hence $|v|\le 2r$. As
$|v\otimes g|\le |v||g|$ and $|v\otimes v|=|v|^2$, we thus have
\[
    |\mathsf{M}(y,z)|
    \le
    \bigl(4r\zeta_r(y)+8r^2|\nabla\zeta_r(y)|\bigr)
    |\mathfrak{f}(z)|.
\]
Hence, we deduce that
\begin{align*}
    \int_0^1 \|S_\lambda\|_{L^1}\,\dd\lambda
    &\le
    \int_0^1\iint_{\R^2\times\R^2}
       |\mathsf{M}(y,z)|\,\dd y\,\dd z\,\dd\lambda\\
    &\le
    \bigl(4r\|\zeta_r\|_{L^1}
          +8r^2\|\nabla\zeta_r\|_{L^1}\bigr)
       \|\mathfrak{f}\|_{L^1}\\
    &\le Cr\|\mathfrak{f}\|_{L^1}
\end{align*}
by a change of variables and Fubini's theorem. Here $C$ depends only on the choice of $\zeta$.

\smallskip
\noindent
{\bf Step~2.}  As $X_\lambda$ parametrises the segment joining $y$ to $z$, one has that
\[
    \operatorname{supp}S
    \subset
    \operatorname{conv}\bigl(
       \operatorname{supp}\zeta_r\cup
       \operatorname{supp}\mathfrak{f}
    \bigr)
    \Subset \B_r(x_0).
\]

To see that $S \in \CC^\infty$, notice that we can differentiate under the integral signs in the formulae for $S$ and $S_\lambda$. This follows from the dominated convergence theorem, as the integrands are compactly supported in the $x$ variable.

Indeed, when $\lambda\ge 1/2$, in the above expression~\eqref{Slambda} for $S_\lambda$,  all powers of $\lambda^{-1}$ generated by $x$-derivatives are bounded. When $\lambda\le 1/2$, by shifting the integration variable from $y$ to 
$z=y+(x-y)/\lambda$ yields that
\[
    S_\lambda(x)
    =
    \int_{\R^2}(1-\lambda)^{-2}
    \mathsf{M}\!\left(\frac{x-\lambda z}{1-\lambda},z\right)
    \,\dd z.
\]
Here, all powers of $(1-\lambda)^{-1}$
generated by $x$-derivatives are bounded. Therefore, $S \in \CC^\infty$.

\smallskip
\noindent
{\bf Step~3.} Take any test vector field $\varphi\in\CC_c^\infty(\R^2;\R^2)$ and put $E_{j\ell}:=\frac12
    \bigl(\partial_j\varphi^\ell+\partial_\ell\varphi^j\bigr)$. As $S$ is symmetric, we have  $S:\nabla\varphi=S:E$, so by a change of variables
\begin{align}
    \int_{\R^2}S:\nabla\varphi\,\dd x
    =\int_0^1\iint_{\R^2\times\R^2}
    \Bigl[
       &-2\zeta_r v^j \mathfrak{f}^\ell E_{j\ell}
       +(\partial_k\zeta_r)v^jv^k \mathfrak{f}^\ell E_{j\ell}
       \notag\\
       &-(\partial_k\zeta_r)v^jv^\ell \mathfrak{f}^k E_{j\ell}
    \Bigr]\,\dd y\,\dd z\,\dd\lambda.
    \label{eq:IO-pairing}
\end{align}
Here and hereafter in this proof, $E$ is evaluated at $X_\lambda$, $\zeta_r$ at $y$, and $\mathfrak{f}$ at $z$.

Let us integrate the terms containing $\partial_k\zeta_r$ by parts in
$y$ while holding $z$ fixed. To this end, observe that $   \partial_{y^k}v^j=-\delta_k^j$, $    \partial_{y^k}(v^jv^k)=-3v^j$, as well as $\partial_{y^k}E_{j\ell}(X_\lambda)
    =(1-\lambda)\partial_k E_{j\ell}(X_\lambda)$. We thus infer from~\eqref{eq:IO-pairing} that
\[
   \int_{\R^2}S:\nabla\varphi\,\dd x
    =\int_0^1\iint_{\R^2\times\R^2} \zeta_r\Bigl[
        -v^j \mathfrak{f}^\ell E_{j\ell}
        +(1-\lambda)
        \bigl(v^jv^\ell \mathfrak{f}^k\partial_k E_{j\ell}
              -v^jv^k \mathfrak{f}^\ell\partial_k E_{j\ell}\bigr)
    \Bigr]\,\dd y\,\dd z\,\dd\lambda.
\]

To proceed, for fixed $y$ and $z$, denote
\begin{align*}
\begin{cases}
 A(\lambda):=v^j \mathfrak{f}(z)^\ell\partial_j\varphi^\ell(X_\lambda),\\
 B(\lambda):=v^j \mathfrak{f}(z)^\ell\partial_\ell\varphi^j(X_\lambda).
\end{cases}
\end{align*} 
Since $\frac{\dd}{\dd\lambda}X_\lambda=v$, we have $    v^j \mathfrak{f}(z)^\ell E_{j\ell}=\frac{A+B}{2}$, $   v^jv^\ell \mathfrak{f}(z)^k\partial_k E_{j\ell}=B'$,  and $v^jv^k \mathfrak{f}(z)^\ell\partial_k E_{j\ell}=\frac{A'+B'}{2}$. Hence, the integrand in the above formula equals
\begin{align*}
&-v^j \mathfrak{f}^\ell E_{j\ell}+(1-\lambda)
        \bigl(v^jv^\ell \mathfrak{f}^k\partial_k E_{j\ell}  -v^jv^k \mathfrak{f}^\ell\partial_k E_{j\ell}\bigr)\\
        &\qquad = -\frac{A+B}{2}
    +\frac{1-\lambda}{2}(B'-A')
    =
    -A+\frac{\dd}{\dd\lambda}
          \left[\frac{1-\lambda}{2}(B-A)\right].
\end{align*}
Since $A=\frac{\dd}{\dd\lambda}
(\mathfrak{f}(z)\cdot\varphi(X_\lambda))$, we then have
\begin{align*}
  \int_0^1\left\{ -A+\frac{\dd}{\dd\lambda}
          \left[\frac{1-\lambda}{2}(B-A)\right]\right\}\,\dd\lambda &=  -\mathfrak{f}(z)\cdot\bigl(\varphi(z)-\varphi(y)\bigr)
    +\frac12\bigl(A(0)-B(0)\bigr)\\
    &=  -\mathfrak{f}(z)\cdot\bigl(\varphi(z)-\varphi(y)\bigr)
    +\frac12 [v\wedge \mathfrak{f}(z)] \curl \varphi(y),
\end{align*} 
where $v\wedge \mathfrak{f}=v^1\mathfrak{f}^2-v^2\mathfrak{f}^1$ and $\operatorname{curl}\varphi=\partial_1\varphi^2-\partial_2\varphi^1$. As $\mathfrak{f}$ has zero mean, the term containing $\varphi(y)$ vanishes after integration in $z$. Also, the skew term vanishes because
\[
    \int_{\R^2}(z-y)\wedge\mathfrak{f}(z)\,\dd z
    =
    \operatorname{Torque}(\mathfrak{f})
    -y\wedge\int_{\R^2}\mathfrak{f}(z)\,\dd z
    =0.
\]
Finally, applying $\int\zeta_r=1$ in the remaining term, we obtain that
\[
    \int_{\R^2}S:\nabla\varphi\,\dd x
    =-\int_{\R^2}\mathfrak{f}\cdot\varphi\,\dd x,
\]
which proves $\operatorname{div}S=\mathfrak{f}$. All the arguments above involve no derivatives
of $\mathfrak{f}$, so the same formula also gives a compactly
supported symmetric $L^1$-tensor field for compactly supported $\mathfrak{f} \in L^1$. This completes the proof. 
\end{proof}

\subsection{Lamb--Chaplygin dipoles}\label{subsec: dipole}

We now present an exact solution to the Euler equations on $\R^2$, known as the Lamb--Chaplygin dipole (\textit{cf}. \cite{Chaplygin, Lamb}), which will serve as the principal building block for the convex integration scheme. Our scheme is an adaptation of the construction in Bru\`e--Colombo--Kumar~\cite{BCK}; see~\cite[\S 2]{BCK} for detailed expositions on the Lamb--Chaplygin dipole.

The Lamb--Chaplygin dipole arises as a natural desingularisation of the doublet flow. The latter is a divergence-free potential flow on $\R^2 \setminus \{0\}$; see~\cite{mvh, Lamb}. It glues the doublet flow in $\R^2 \setminus \B_1(0)$ with a desingularised flow with nonvanishing vorticity in $\B_1(0)$.

It is most convenient to describe the Lamb--Chaplygin dipole in polar coordinates $(R,\vartheta)$: Let $k$
be the first positive zero of the Bessel function $J_1$. Set
\[
 \psi(R,\vartheta):=
 \begin{cases}
 \left(\frac{2J_1(kR)}{kJ_0(k)}-R\right)\sin\vartheta&\qquad
 \text{when }R\le1,\\
 -\frac{1}{R}\sin\vartheta &\qquad  \text{when } R>1.
 \end{cases}
\]
Then, $\psi$ and its radial derivative match at $R=1$. It is the stream function of the Lamb--Chaplygin dipole: the velocity field $U$, the vorticity $\omega$, and the pressure $P$ are given by
\begin{align}\label{U, omega, P, unrotated}
&U = \na^\perp \psi \in \CC^{1,1}(\R^2;\R^2),\nonumber\\
&\omega = \curl\, U = -k^2\left(\psi + x^2\right) \1_{\overline{\B}_1(0)} \in \CC^{0,1}(\R^2),\\
&P = U^1- \frac{(U^1)^2 + (U^2)^2}{2} - \frac{k^2}{2}\left(\psi+x^2\right)^2 \1_{\overline{\B}_1(0)}.\nonumber
\end{align}

Define the scaled velocity and pressure with the scaling parameter $r \in ]0,1[$ to be specified:
\begin{equation}\label{Ur, dipole, def}
\hat{U}_r(x):=\frac{1}{r}U\left(\frac{x}{r}\right),\qquad \hat{P}_r(x) := \frac{1}{r^2}P\left(\frac{x}{r}\right).
\end{equation}
They satisfy
\begin{align}\label{S2'}
-\frac{1}{r}\p_1 \hat{U}_r+
 \operatorname{div}(\hat{U}_r\otimes \hat{U}_r)+\nabla \hat{P}_r=0,
\end{align}
as well as
\begin{equation}\label{S3'}
\hat{U}_r = r\na \Phi, \quad \Phi(x) = -\frac{x^1}{|x|^2},\quad\text{and } \hat{P}_r= \p_{x^1} \Phi - \frac{r^2 \left|\na\Phi\right|^2}{2}\qquad \text{when } |x|>r.
\end{equation}
See \cite[Eqs.~(2.8), (2.9)]{BCK}.

In this paper, we also need to  apply an $SO(2)$-change of variables and express $\psi$, $U_r$, and $P_r$, etc. in the new variables. For any unit vector $\xi=(\xi_1,\xi_2)^{\mathsf T}$, write
$\xi^\perp=(-\xi_2,\xi_1)$ and set
\[
 Q_\xi=
 \begin{pmatrix}\xi_1&-\xi_2\\ \xi_2&\xi_1\end{pmatrix}
       =[\xi\;\;\xi^\perp]\in SO(2).
\]
Note that $Q_\xi \mathcal{J}=\mathcal{J}Q_\xi$, where $\mathcal{J} = [x \mapsto x^\perp]$ is the standard almost complex structure on $\C$. Then $Q_\xi \mathbf{e}_1=\xi$ and
$Q_\xi^{\mathsf T}x=(x\cdot\xi,x\cdot\xi^\perp)$. Define the rotated streamfunction and velocity, as well as their
rescaled versions, by
\begin{align}\label{def, Ur, psi-r}
 \psi_\xi(x)&=\psi(Q_\xi^{\mathsf T}x),&
 U_\xi(x)&=\nabla^\perp\psi_\xi(x)
          =Q_\xi U(Q_\xi^{\mathsf T}x),\nonumber\\
 \psi_{\xi,r}(x)&=\psi_\xi(x/r),&
 U_{\xi,r}(x)&=\nabla^\perp\psi_{\xi,r}(x)
          =r^{-1}Q_\xi U(Q_\xi^{\mathsf T}x/r).
\end{align}
Here $\psi$ is evaluated as a function on $\R^2$.
The velocity identities follow from the chain rule and
$Q_\xi \mathcal{J}=\mathcal{J}Q_\xi$.

Below we
abbreviate $U_{\xi,r}$ as $U_r$ when $\xi$ is fixed. Note that $U_{{\mathbf{e}_1}, r} \equiv \hat{U}_r$.

Recall $U$ and $P$ from~\eqref{U, omega, P, unrotated}. Set
\begin{equation}\label{Pr, def}
P_{\xi, r}(x) \equiv P_r(x)=r^{-2}P(Q_\xi^{\mathsf T}x/r).
\end{equation}
It then holds for $y:=Q_\xi^{\mathsf T}x/r$ that
\begin{align*}
 r^{-1}(\xi\cdot\nabla_x)U_r(x)
   &=r^{-3}Q_\xi(\mathbf{e}_1\cdot\nabla_y)U(y),\\
 \operatorname{div}_x(U_r\otimes U_r)(x)
   &=r^{-3}Q_\xi\operatorname{div}_y(U\otimes U)(y),\\
 \nabla_x P_r(x)&=r^{-3}Q_\xi\nabla_y P(y). 
\end{align*}
Hence, we have the PDE for $U_r \equiv U_{\xi, r}$:
\begin{equation}\label{S2}
-\frac{1}{r }(\xi\cdot\nabla)U_r+
 \operatorname{div}(U_r\otimes U_r)+\nabla P_r=0.
\end{equation} 
Equivalently, $U_r(x-t\xi/r)$ is the traveling dipole with
constant velocity $\xi/r$. Its streamfunction in the moving frame
is $\psi_{\xi,r}(x)+(x\cdot\xi^\perp)/r$, whose perpendicular
gradient is $U_r-\xi/r$.

The vorticity satisfies
$$\operatorname{curl} U_r(x)=\frac{1}{r^2}(\Delta\psi)(Q_\xi^{\mathsf T}x/r)$$ and is supported in $\overline{\B}_r(0)$. For $|x|>r$, it follows from $\psi(y)=-y_2/|y|^2$ that
$$\psi_{\xi,r}(x)=-r \frac{x\cdot\xi^\perp}{|x|^2}.$$ Differentiation leads to
\begin{equation}\label{S3}
U_r(x)=r\nabla\Phi_\xi(x),\qquad
 \Phi_\xi(x)=-\frac{\xi\cdot x}{|x|^2}.
\end{equation}
Compare Equations~\eqref{S2}, \eqref{S3} with \eqref{S2'}, \eqref{S3'}, respectively.

In addition, for a universal constant $C$ we have that
\begin{align}\label{S4}
&\|{U}_r\|_{L^2} \equiv C,\qquad \|{U}_r\|_{L^p} =r^{\frac{2}{p}-1}\|U\|_{L^p},\nonumber\\
& \|\curl \, {U}_r\|_{L^1} \equiv \|\curl\, U\|_{L^1}
\end{align}
for arbitrary $r>0$. The velocity is even and its stream function and exterior potential are odd. 
 
\begin{figure}[htbp]
\centering\includegraphics[width=\linewidth]{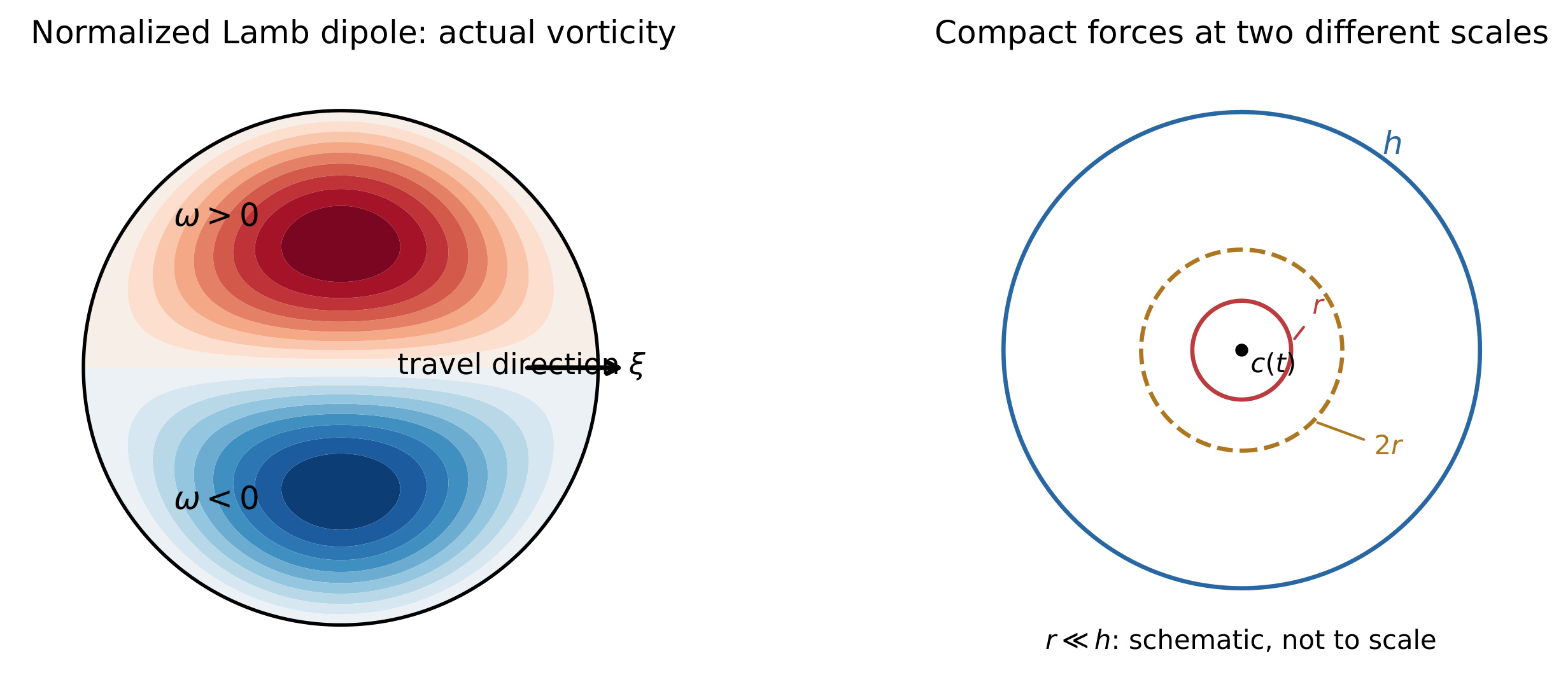}
\caption{Normalised dipole vorticity (left) and schematic scales (right). Only vorticity and forces have the displayed compact supports; the principal velocity has an
irrotational tail.}
\label{fig:dipole-scales}
\end{figure}

In the sequel, we shall refer to $r>0$ as the \emph{core radius}, since inside $\B_r(0) \subset \R^2$ the rescaled velocity field $U_r$ is the desingularised doublet flow. As a foreshadow of the construction in \S\ref{subsec: Traversing Lamb--Chaplygin dipoles} below, we shall centre the rescaled Lamb--Chaplygin dipole $U_r$ at $\cen$ to translate it by solving an ODE for $\cen=\cen(t)$. Another scale $h \gg r$ to be specified shall play a key role. See Figure~\ref{fig:dipole-scales}.

\section{The perturbation lemma}\label{sec: perturbation lemma}

This section is devoted to the proof of the key Lemma~\ref{lem: key}, reproduced below:

\begin{lemma*}
Let $I \subset \R$ be a closed finite interval. Let $(u,p,R)$ satisfy
\begin{align*}
    &\partial_tu+\operatorname{div}(u\otimes u)+\nabla p
       =\operatorname{div}R,\\
    &\operatorname{div}u=0,
\end{align*}
where $u \in \left(\CC^\infty_{t,x} \cap \CC^0_t L^3_x\right) \left(I \times \R^2;\R^2\right)$, $\omega \in \CC^0_t L^1_x\left(I \times \R^2\right)$, and $R \in \left(\CC^\infty_{t,x} \cap\CC^0_t L^1_x\right)\left(I \times \R^2; \sym\right)$. Assume that $u-\srad \in \CC^0_t L^2_x(I \times \R^2;\R^2)$ for some fixed smooth radial vector field $\srad$ and $p(t, 0) \equiv 0$. Also assume that $R(t,\cdot) \equiv 0$ for $t$ in an open neighbourhood of $\p I$, and     that 
\begin{align*}
    \sup_{t \in I} \|R(t,\cdot)\|_{L^1_x} \leq \delta\qquad\text{for some $\delta>0$}.
\end{align*}
Then, for any given $\e_\star>0$, we can find a smooth triple $\left(\widetilde{u}, \widetilde{p}, \widetilde{R}\right)$ satisfying
\begin{align*}
    &\partial_t\widetilde{u}+\operatorname{div}\left(\widetilde{u}\otimes \widetilde{u}\right)+\nabla \widetilde{p}
       =\operatorname{div}\widetilde{R},\\
    &\operatorname{div}\widetilde{u}=0,
\end{align*}
such that for some  universal constant $C_0>0$,
\begin{align*}
&\left(\widetilde{u}, \widetilde{p}, \widetilde{R}\right)(t,\cdot) \equiv (u,p,0)(t,\cdot)\quad \text{ for each $t$ in an open neighbourhood of $\p I$},\\
&\sup_{t\in I}\left\|\widetilde{R}(t,\cdot)\right\|_{L^1_x} <\e_\star,\\
& \sup_{t\in I}\left\|\widetilde{u}(t,\cdot)-u(t,\cdot)\right\|_{L^2_x} +\sup_{t\in I}\left\|\operatorname{curl}\left(\widetilde{u}-u\right)(t,\cdot)\right\|_{L^1_x}
 \leq C_0\sqrt\delta+\e_\star.
\end{align*} 
Moreover, for each $t>0$, the correction term $\left(\widetilde{u}-u\right)(t,\cdot)$ has zero circulation, and $\curl\left(\widetilde{u}-u\right)(t,\cdot)$ has the vanishing vector first moment and scalar second moment, and is compactly supported in space. The updated triple lies in the same regularity classes as for $(u,p,R)$; that is, $\widetilde{u} \in \left(\CC^\infty_{t,x} \cap \CC^0_t L^3_x\right) \left(I \times \R^2;\R^2\right)$, $\curl\widetilde{u} \in \CC^0_t L^1_x\left(I \times \R^2\right)$, and $\widetilde{R} \in\left(\CC^\infty_{t,x} \cap\CC^0_t L^1_x\right)\left(I \times \R^2; \sym\right)$.  
\end{lemma*}

The proof occupies the next eleven subsections.

\subsection{Travelling Lamb--Chaplygin dipoles}\label{subsec: Traversing Lamb--Chaplygin dipoles}

Recall \S\ref{subsec: dipole} for the construction of the rescaled Lamb--Chaplygin dipole. We follow Bru\`e--Colombo--Kumar~\cite[\S 1.2.3 and \S 2.5]{BCK} to vary the amplitude, the core radius, and the centre of the dipole velocity profile $U_r \equiv U_{\xi, r}$, with some unit vector $\xi$ to be specified later: Consider for some smooth functions $\eta=\eta(t)$ and $r=r(t)$ that
\begin{equation}\label{def, V, travelling dipole}
 V(t,x)=\eta(t)U_{r(t)}(x-\cen(t)),
\end{equation} 
where the location $\cen$ of the centre is determined from the ODE:
\begin{equation}\label{c'}
\cen'(t)=\frac{\eta(t)}{r(t)}\xi.
\end{equation}
The unit vector $\xi \in \R^2$ is the direction in which the dipole is travelling. 

We compute that
\begin{align*}
\partial_tV(t,x) &=\eta'(t)U_{r(t)}
\big(x-\cen(t)\big)+\eta (t) r'(t)\partial_rU_{r(t)}\big(x-\cen(t)\big)\\
&\qquad - \eta(t) \cen'(t)\cdot\nabla U_{r(t)}\big(x-\cen(t)\big)\\
&=\eta'(t)U_{r(t)}
\big(x-\cen(t)\big)+\eta (t) r'(t)\partial_rU_{r(t)}\big(x-\cen(t)\big)\\
&\qquad -\frac{\eta^2(t)}{r(t)} \xi \cdot \na U_{r(t)}
\big(x-\cen(t)\big)\\
&= \eta'(t)U_{r(t)}
\big(x-\cen(t)\big)+\eta (t) r'(t)\partial_rU_{r(t)}\big(x-\cen(t)\big)\\
&\qquad - {\rm div}\Big(V(t,x) \otimes V(t,x)\Big) + \na \Big(\eta^2(t) P_r(t,x)\Big).
\end{align*}
The first equality follows from the definition of $V$ in ~\eqref{def, V, travelling dipole}, the second equality from the ODE~\eqref{c'} for $\cen$, and the final one from the PDE~\eqref{S2} for $U_r$, with the rescaled pressure $P_r$ defined in~\eqref{Ur, dipole, def}. We write it more compactly as
\begin{equation} \label{S5}
 \partial_t V+\operatorname{div}(V\otimes V)+\nabla(\eta^2P_r)
       =\eta' U_r+\eta r'\partial_r U_r.
\end{equation}
In the subsequent developments, $V=V(t,x)$ shall serve as the principal building block of the convex integration scheme.

We next introduce a decomposition of the profile $V$, which is a variant of \cite[\S 2.3]{BCK}. Fix a radial cutoff function $\chi=\chi(|x|) \in \CC^\infty([0,\infty[)$ such that $\chi \equiv 0$ on $[0,1]$ and $\chi \equiv 1$ on $[2,\infty[$. Define for any $r>0$:
\begin{equation}\label{Pi-r, def}
\Pi_r(x) := r\chi\left(\frac{|x|}{r}\right)\Phi_\xi(x),
\end{equation}
which is an odd extension and truncation of the exterior velocity potential $\Phi_\xi$. Set
\begin{equation}\label{Wr, def}
W_r(t,x) := U_{r}(t,x)-\nabla\Pi_{r}(t,x).
\end{equation} 
We also introduce the forcing term, allowing that $r$ varies in $t$:
\begin{equation}\label{Fc, def}
 F_c(t,x):=\eta'(t) W_{r(t)}\big(x-\cen(t)\big)+\eta(t) r'(t)\partial_rW_{r(t)}\big(x-\cen(t)\big).
\end{equation}

\begin{lemma}\label{lem: Wr, Fc}
The velocity field $W_r$ is an even vector field compactly supported in $\B_{2r}$ with 
\begin{equation}\label{S7}
 \int_{\R^2} W_{r(t)}(t,x)\,\dd x=2\pi r(t)\,\xi.
\end{equation}

The vector field $F_c$ is the external forcing for $V=V(t,x)$. The PDE for $V$ is
\begin{equation}\label{Euler for V with Fc}
\p_tV + {\rm div}(V \otimes V) + \na p_V = F_c,
\end{equation}
where $p_V$ is a scalar field. Moreover, $F_c$ is compactly supported in $\B_{2r(t)}$ with the properties
\begin{align}
&\|F_c(t,\cdot)\|_{L^1_x}\le C(|\eta'(t)|r(t)+|\eta(t)|\,|r'(t)|),\label{eq:scan-F-bound}\\
 &\int_{\R^2} F_c(t,x)\,\dd x=2\pi\frac{d}{dt}\big(\eta(t) r(t)\big) \, \xi,
 \label{eq:scan-F-mean}\\
 &\int_{\R^2}\big(x-\cen(t)\big)^j (F_c)^k(x)\,\dd x = 0\qquad \text{ for all } j,k \in \{1,2\},\label{S8}
\end{align}
where $C$ is a universal constant.  
\end{lemma} 

In fact, the sign convention in \eqref{S2}--\eqref{S3} is adopted in light of~\eqref{S7}.

\begin{proof}[Proof of Lemma~\ref{lem: Wr, Fc}] 

By the definition of $W_r$ in~\eqref{Wr, def} and that of $U_r$, $P_r$ in~\eqref{def, Ur, psi-r} and \eqref{Pi-r, def}, respectively, we see that ${\rm supp}(W_r) \Subset \B_{2r}$. In particular, the integral in~\eqref{S7} converges. On a large circle the boundary integral for the  stream function $\psi_r$ gives us $\int_{\R^2} U_r(t,x)\,\dd x=\pi r\xi$, whereas the boundary integral for $\Phi_\xi$ gives us 
$\int_{\R^2}\nabla\Pi_r(t,x)\,\dd x=-\pi r\xi$. This proves~\eqref{S7}.

In view of the definition of $W_r$ in~\eqref{Wr, def}, by setting 
$$p_V=\eta^2P_r-\eta'\Pi_r-\eta r'\partial_r\Pi_r$$ we immediately deduce the forced Euler Equation~\eqref{Euler for V with Fc} for $W_r$. The support of the forcing $F_c$ is contained in that of $W_{r(t)}$. Identity~\eqref{eq:scan-F-mean} follows from the definition~\eqref{Fc, def} of $F_c$ and~\eqref{S7} established above. Moreover,  Identity~\eqref{S8} holds since $F_c$ is even with respect to the centre $c(t)$.

Finally, as $W_r(x)=r^{-1}W_1(x/r)$, we have $ \partial_rW_r(x)=-r^{-2}[W_1(z)+z\cdot\nabla W_1(z)]_{z=x/r}$, so that 
$\|W_r\|_{L^1_x}=r\|W_1\|_{L^1_x}$ and
$\|\partial_rW_r\|_{L^1_x}\le\|W_1\|_{L^1_x}+\||x|\nabla W_1\|_{L^1_x}$.  Therefore, \eqref{eq:scan-F-bound} holds by taking $C=C\left(\|W_1\|_{L^1_x}, \||x|\nabla W_1\|_{L^1_x}\right)$. By the definition of $W_r$ in~\eqref{Wr, def}, $C$ depends only on the $\CC^{1,1}$-norm of $U$, and hence on the $\CC^{0,1}$-norm of $\psi$ (see~\eqref{U, omega, P, unrotated}), which is universal.  \end{proof}

As shown in \eqref{eq:scan-F-mean}, the integral over $\R^2$ of the forcing $F_c$ is nonzero in general, so \(F_c\) cannot simply be written as the divergence of a compactly supported tensor. We thus need a further modified version of the forcing $F_c$. This is postponed to \S\ref{subsec: aux force}, for it involves a parameter $h$ detailed later in \S\ref{subsec: parameter D,h}.

\subsection{Decomposition of error; geometric configurations}

In this section, we decompose the Reynolds stress tensor into a linear combination of ``good directions'' $\{\xi_i: i = 1,2,3,4\}$, modulo a harmless term that goes to the pressure. The arguments here essentially follow from \cite[\S 4.2]{BCK}, and the choice of $\xi_i$ is, in turn, taken from Bru\`{e}--Colombo ~\cite{BC}.

Pick four unit vectors:
\begin{align*}
& \xi_1= \begin{bmatrix}
1\\
0
\end{bmatrix},\qquad \xi_2 = \begin{bmatrix}
0\\
1
\end{bmatrix},\\
& \xi_3 = \begin{bmatrix}
\frac{1}{\sqrt{2}}\\
\frac{1}{\sqrt{2}}
\end{bmatrix}, \qquad \xi_4 = \begin{bmatrix}
\frac{1}{\sqrt{2}}\\
-\frac{1}{\sqrt{2}}
\end{bmatrix}.
\end{align*}

Let the time interval $I$, the initial Reynolds stress tensor $R$, and its $L^1_x$-bound $\delta$ be as in the statement of Lemma~\ref{lem: key}. Take a small parameter $\e_0 \in ]0,1[$ as follows: 
\begin{equation}\label{e0}
\e_0 :=  \min\left\{
\frac12,\ \frac{\e_\star}{200},\
\frac{\e_\star}{16\sqrt{3/\pi}\,
  (\|U\|_{L^2_x}+\|\operatorname{curl}(U)\|_{L^1_x})}
\right\},
\end{equation}
where $\e_\star$ is as given in the statement of Lemma~\ref{lem: key}, and $U$ is the velocity field of the (unscaled, unrotated) Lamb--Chaplygin dipole; see~\eqref{U, omega, P, unrotated}. In particular, $\e_0$ is chosen to depend only on $\e_\star$. Observe that the image $R(I)$ is compact in $L^1_x$, so
\begin{equation}\label{tail}
\sup_{t\in I}\int_{|x|>K}|R(t,x)|\,\dd x\longrightarrow 0
\end{equation}
as \(K\to\infty\). Then, take 
\begin{equation}\label{L, choice}
L=L\left(\e_0,\delta,\|R\|_{\CC^0_t L^1_x(I \times \R^2)}, R(I)\right)>1000
\end{equation} 
so that the tail bound in~\eqref{tail} at \(K=L/2\) is at most \(\varepsilon_0\).

We shall fix once and for all a cutoff function $\chi \in \CC^\infty_c(\B_1)$ with $0 \leq \chi \leq 1$, such that
\begin{equation}\label{def, Rc}
R_c(t,x) := \chi\left(\frac{x}{L}\right) R(t,x) 
\end{equation} 
satisfies
\begin{equation}\label{S9}
\begin{cases}
\|R-R_c\|_{\CC^0_t L^1_x (I \times \R^2)} \leq \e_0,\\
 \|R_c\|_{\CC^0_t L^1_x (I \times \R^2)} \leq \delta.
 \end{cases}
\end{equation} 
 
Now let us fix the small positive number 
\begin{equation}\label{floor}
\flf := \frac{\e_0^2}{1024 L^2}.
\end{equation}
With $L$ chosen as above, $\flf$ depends only on $\e_0$, $\delta$, and $\|R\|_{\CC^0_t L^1_x}$. Then set
\begin{equation}\label{S10}
\begin{cases}
b=b(t,x):=16\sqrt{|R_c|^2+{\flf}^2},\\
a_1=\frac{b}{2}-(R_c)_{11}+(R_c)_{12},\\ 
a_2:=\frac{b}{2}-(R_c)_{22}+(R_c)_{12}\\
a_3 :=\frac{b}{2}-2(R_c)_{12},\\
a_4:=\frac{b}{2}.
\end{cases}
\end{equation}
One may directly check that
\begin{equation}\label{pointwise bound for ai}
\sqrt{62}\flf\le a_i\le (8+\sqrt{2})(|R_c|+\flf)\qquad\text{for each } i \in \{1,2,3,4\}.
\end{equation}  
Also,
\begin{equation}\label{ai outside supp Rc}
a_i \equiv 8\flf \quad\text{ on } \R^2 \setminus {\rm supp}(R_c) \text{ for each } i \in \{1,2,3,4\}.
\end{equation}
It is crucial to notice that 
\begin{equation}\label{S11}
-R_c(t,x) = \sum_{i=1}^4\Big(a_i(t,x)\xi_i\otimes\xi_i\Big) - b(t,x){\bf I}_2.
\end{equation}

Next let us fix the four rotated squares of sidelength $4L \gg r(t)$, oriented in $\xi_i$-directions:
\begin{equation}\label{square Qi}
 Q_i :=\Big\{\alpha\xi_i+ \beta\xi_i^\perp:\,|\alpha|<2L,\ |\beta|<2L\Big\}.
\end{equation}
Squares are used only to fix
the geometry conveniently; other rectangles with sidelengths comparable to $L$ also work. Note that $\frac{1}{2}Q_i \Supset \B_L \Supset {\rm supp}(R_c)$. In view of~\eqref{square Qi} and~\eqref{floor}, we have
\begin{equation*}
8\flf\sum_i\leb^2(Q_i)=512\flf L^2<\e_0^2,
\end{equation*}
and hence
\begin{align}\label{S12}
 \int_{Q_i}a_i(t,x)\,\dd x\le (8+\sqrt{2})\delta+ \frac{1}{4}\e_0^2
\end{align} 
by the bound~\eqref{pointwise bound for ai} for $a_i$ and the estimate~\eqref{S9} for $R_c$.



\subsection{Stripes}

Let the time interval $I$, the initial Reynolds stress tensor $R$, and its $L^1_x$-bound $\delta$ be as in the statement of Lemma~\ref{lem: key}. Choose $R_c$ and hence $a_i$, $b$ as in~\eqref{def, Rc} and \eqref{S10} as in the previous subsection.

Note that a corner of the square $Q_i$ is
\begin{equation}\label{corner}
\ci := -2L \left(\xi_i+\xi_i^\perp\right).
\end{equation} 
We take a large number $N_i$ of stripes (to be specified later) inside $\overline{Q_i}$, emanating from equally spaced points on the edge $\left\{\ci + s\xi_i^\perp:\,s\in [0,4L]\right\}$ and traversing along the $\xi$-direction. More precisely, set the spacing
\begin{equation}\label{si, footstep}
s_i:= \frac{4L}{N_i}. 
\end{equation}
The $j^{\text{th}}$ stripe $\mathbf{X}_{ij}$ is then parametrised by 
\begin{equation}\label{parametrisation, stripe}
\stripeij (\lambda):= \ci +\lambda\xi_i+\left(j-\frac{1}{2}\right)s_i \xi_i^\perp,
 \qquad 0\le\lambda\le4L.
\end{equation}

Then, define the integral of $a_i$ along the $j^{\text{th}}$ stripe in $Q_i$ as
\begin{equation}\label{S13}
I_{ij}(t) := \int_0^{4L} a_i\big(t, \stripeij(\lambda)\big)\,\dd\lambda. 
\end{equation}
The ``quadrature''
\begin{equation}\label{Delta-i, quadrature, def}
\Delta_i(t) := s_i \cdot \sum_{j=1}^{N_i} I_{ij}(t)
\end{equation}
approximates the area integral of $a_i$ along $Q_i$. Indeed, 
\begin{align}\label{FTC'}
\left|\Delta_i(t) -\int_{Q_i}a_i(t,x)\,\dd x\right|
 \le s_i\int_{Q_i} \left|{\xi_i^\perp} \cdot \nabla a_i(t,x)\right|\,\dd x \leq s_i \left\|\na a_i\right\|_{L^\infty_t L^1_x(I \times Q_i)}, 
\end{align} 
thanks to the fundamental theorem of calculus.

\begin{figure}[htbp]
\centering\includegraphics[width=0.8\linewidth]{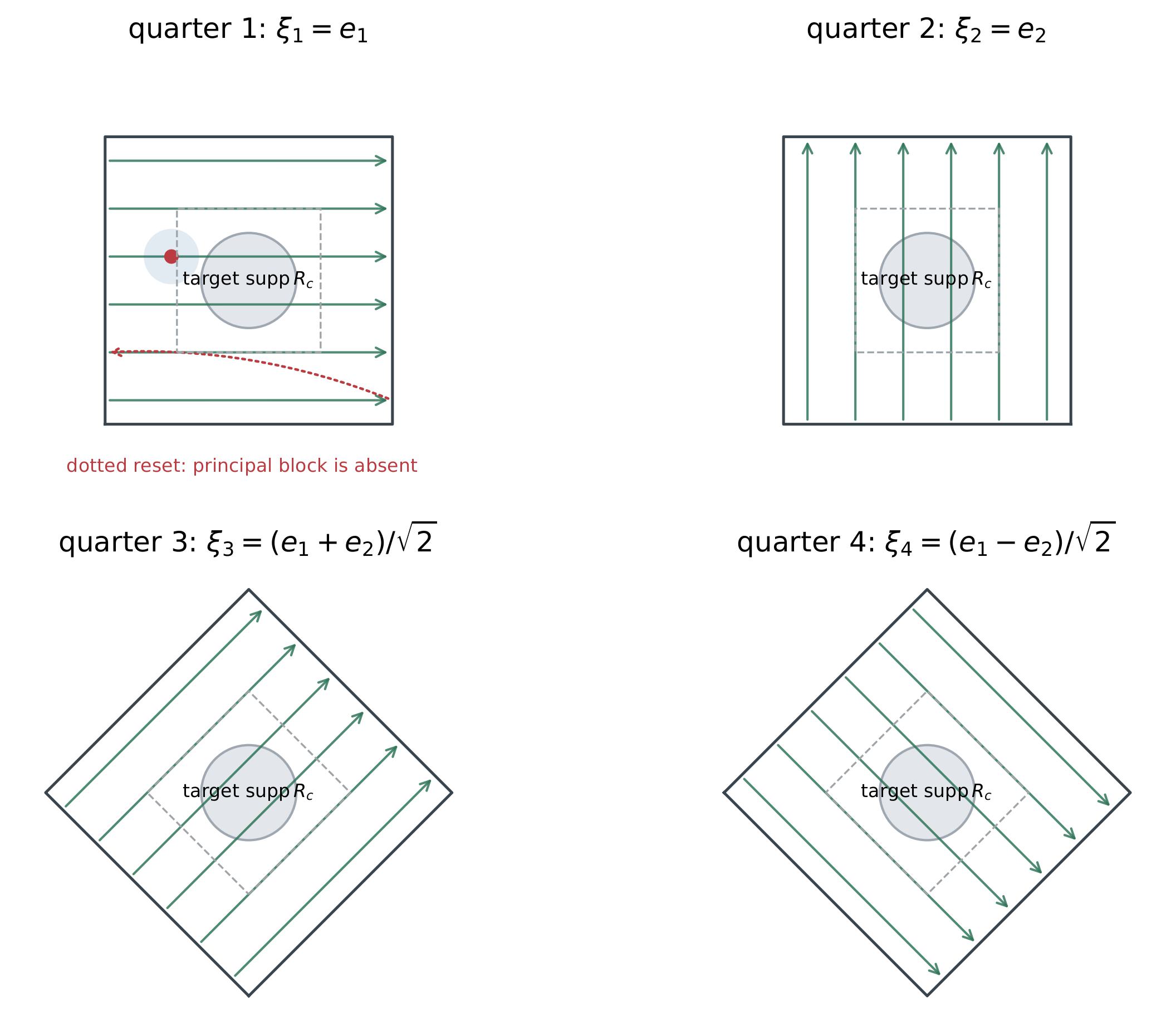}
\caption{The four direction families, used successively. Each square
has side $4L$; the dashed inner square has side $2L$. Arrows indicate
midpoint-line traversals. The dotted connecting path is followed only
while the principal block is absent. The red core and blue auxiliary
bump are schematic and not to scale; only one principal dipole is
active at a time.}
\label{fig:stripe-families}
\end{figure}

See Figure~\ref{fig:stripe-families} for the choice of $\xi_i$, $Q_i$, and $\stripeij$, etc.

\subsection{Choice of parameters}\label{subsec: parameter D,h}

By now, for given $\e_\star$, $\delta$, and the smooth Euler--Reynolds triple $(u,p,R)$ as in the statement of Lemma~\ref{lem: key}, the following parameters have already been chosen:
\begin{align*}
L,\, R_c,\, \e_0,\, \flf,\, b, \, \{a_i\}_{i=1,2,3,4}. 
\end{align*}
In this subsection, let us further specify $D$, $h$, $s_i$, as well as time parameter  $\bigtime$.

\subsubsection{Spatial parameters}
For reasons that will become clear later, we choose the large parameter $D>0$ such that
\begin{equation}\label{D}
D = D\bigg(\e_0,\,\delta,\,L,\,
\|R_c\|_{L^\infty_{t,x}},\,
\|\partial_tR_c\|_{L^\infty_tL^1_x},\,
\|\nabla_x R_c\|_{L^\infty_{t,x}},\,
\|\nabla_xR_c\|_{L^\infty_tL^1_x},\, \|u\|_{L^\infty_t L^3_x}\bigg).
\end{equation} 
Note that $D$ is selected after $L$ and the cutoff stress tensor $R_c$ are chosen. The first four norms in the parenthesis of~\eqref{D} may increase to infinity as $n \nearrow \infty$, where $n$ indexes the time intervals glued together such that $I_n$ shrinks to the initial time $t=0$. We made the choice of $D$ in \eqref{D} on \emph{one fixed time interval} $I$ as in the statement of Lemma~\ref{lem: key}. In what follows, we shall keep careful track of the dependence of $D$ on these indicated parameters.

Now, after $D$ is picked, let us choose another small parameter $h=h(\e,D) >0$ such that
\begin{equation}\label{h, choice}
h := \min\left\{ \frac{1}{128},\,\left(\frac{\e_0}{400 D}\right)^2 \right\}.
\end{equation} 
We then choose $s_i>0$ so small that
\begin{equation}\label{si, choice}
s_i<\min\left\{h,\frac{\e_0^2}{2D},\frac{\e_\star h}{100D}\right\}.
\end{equation}
In view of~\eqref{D}, \eqref{h, choice}--\eqref{si, choice}, $D$, $h$, and $s_i$ all depend on the same set of parameters. Moreover, by~\eqref{si, footstep} the choice of $s_i$ is equivalent to that of $N_i$.

Observe that for $D$ as in \eqref{D}, we may require $\|\na_x a_i\|_{L^\infty_t L^1_x} \leq D$, in view of the expression~\eqref{S10} for $a_i$, the choice of $\flf$ in~\eqref{floor}, and the choice of $L$ in~\eqref{L, choice}. It then follows from~\eqref{S13}, \eqref{FTC'}, \eqref{pointwise bound for ai}, and \eqref{S12} that
\begin{align}\label{S14}
& I_{ij} \geq 4\sqrt{62} \cdot \flf L,\nonumber\\
& \Delta_i \leq (8+\sqrt{2})\delta+ \frac{3}{4}\e_0^2.
\end{align}

The core radius $r=r(t)$ for the Lamb--Chaplygin dipole need not be comparable to $s_i$. Indeed, it is selected
later so that
\begin{equation*}
r(t) \ll s_i <h.
\end{equation*}

\subsubsection{Temporal parameter}

Now let $I = [A,B]$ be the closed finite time interval in Lemma~\ref{lem: key}. By assumption of the lemma, $R(t,\cdot)\equiv 0$ for $t$ in some relatively open neighbourhood of $\p I$. If $R \equiv 0$ on $I \times \R^2$,  one may simply set $\left(\widetilde{u}, \widetilde{p}, \widetilde{R}\right) := (u,p,R)$, and the assertion follows trivially. Thus, in the remaining parts of this section, we assume that $R$ is not constantly zero.

 To fix the idea, label
\begin{equation}
\begin{cases}
t_I^-:=\max\big\{ t>A:\, R(t,\cdot)\equiv 0 \text{ on } [A,t] \big\},\\
t_I^+:=\min\big\{ t<B:\, R(t,\cdot)\equiv 0 \text{ on } [t,B] \big\},\\
J:=\left[t_I^-, t_I^+\right].
\end{cases}
\end{equation}
We then have $t_I^- <t_I^+$.

We shall divide the time interval $J$ into $M$ time cells, each of duration length
\begin{equation}\label{timelength}
\bigtime := \frac{|J|}{M}.
\end{equation} 
Choose $M$ so large that
\begin{equation}\label{Tcell}
\bigtime < \frac{\e_\star}{100 D}; 
\end{equation}
In the $m^{\text{th}}$ cell $$J_m \equiv \left[t_I^- + m\bigtime,\, t_I^- + (m+1)\bigtime\right],$$
we arrange the configurations  $Q_1, Q_2, Q_3, Q_4$ in this order  in each of its four equal quarters. That is, set the ``\emph{stripe durations}'' as 
\begin{equation}\label{duration dij}
d_{ij}^{(m)} := \text{ the time spent on $j^{\text{th}}$ stripe in square $Q_i$ in $m^{\text{th}}$ time cell}.
\end{equation}
Then (recalling the paragraph above~\eqref{si, footstep} for $N_i$) it holds that
\begin{equation}
\sum_{j=1}^{N_i} d^{(m)}_{ij} = \frac{\bigtime}{4}\qquad\text{ for each } i \in \{1,2,3,4\}. 
\end{equation}

\begin{figure}[htbp]
\centering\includegraphics[width=1\linewidth]{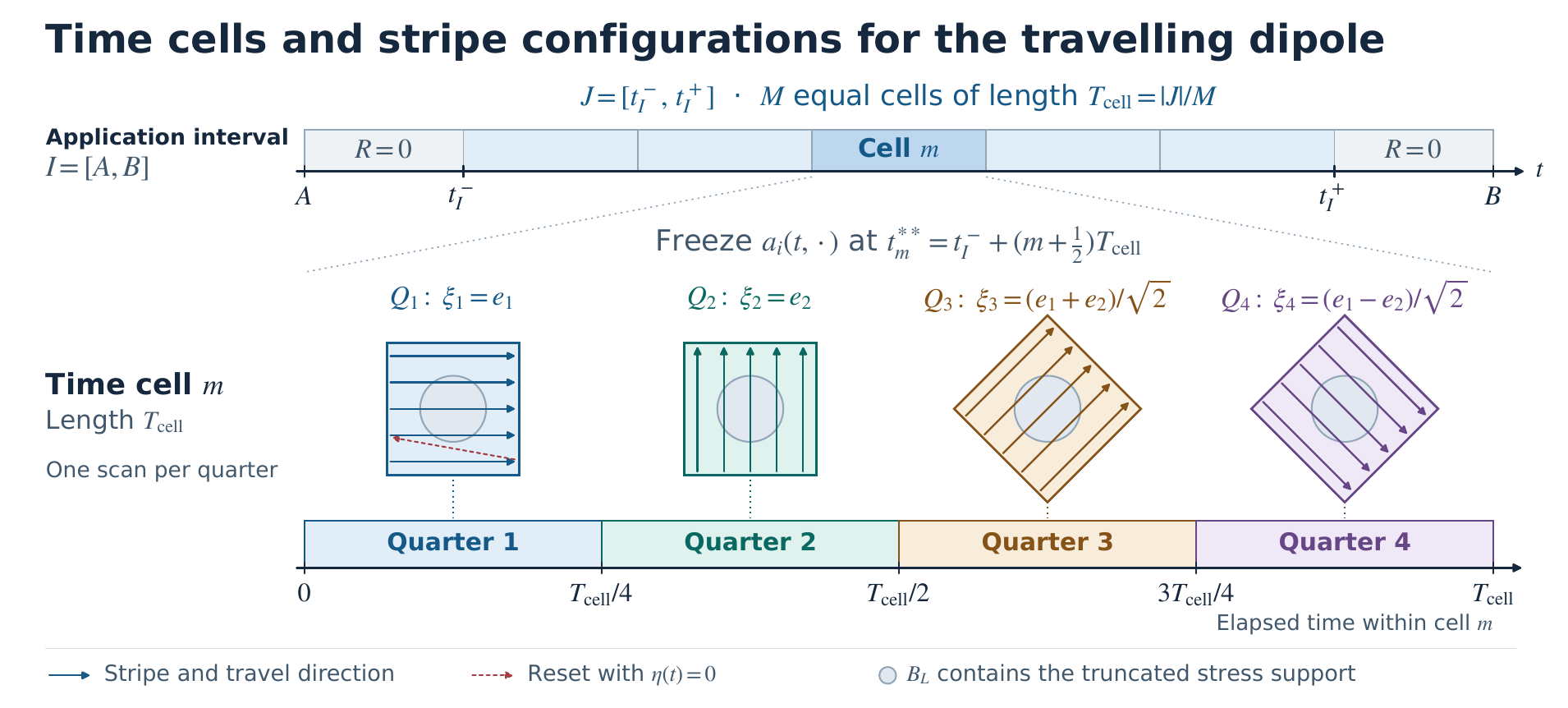}
\caption{Time intervals}
\label{fig:time}
\end{figure}

Among all the parameters introduced by now, only $D$ is yet to be specified.

\subsection{Profiles of the travelling dipole on each stripe}

We now perform a fine-tuning for the profiles of the travelling dipole, which may vary across stripes, directions, and time  cells.  

\subsubsection{Freezing time}

Denote 
\begin{equation}\label{t*, slot}
t^{*}_{m,i,j}=
t_I^-+m\mathsf T_{\mathrm{cell}}
+\frac{(i-1)\bigtime}{4}
+\sum_{\ell=1}^{j-1}d_{i\ell}^{(m)},
\end{equation}
where $m \in \{0,1,\ldots,M-1\}$ indexes time cells, $i \in \{1,2,3,4\}$ labels the directions in which the dipole traverses, and $j \in \{1,2,\ldots, N_i\}$ indexes the stripes within $Q_i$. It is the beginning time for the dipole to travel along the $j^{\text{th}}$ stripe in direction $\xi_i$ within the $m^{\text{th}}$ time cell. Note that 
\begin{equation}\label{+1 in time cell}
t^*_{m,i,j}+d_{ij}^{(m)}=\begin{cases}
t^*_{m,i,j+1},\qquad \text{for } m \in \{0,\ldots, M-1\},\, i \in \{1,2,3,4\},\, j <N_i;\\
t^*_{m,i+1,1}\qquad \text{for } m \in \{0,\ldots,M-1\},\,  i \in \{1,2,3\}, \,j = N_i;\\
t^*_{m+1,1,1}\qquad\text{for } m< M-1, \, i= 4,\,  j=N_4.
\end{cases}
\end{equation} 

Also denote 
\begin{equation}\label{t**, def}
t^{**}_m := t_I^-+\left(m+\frac12\right)\bigtime\qquad\text{for each } m \in \{0,1,\ldots,M-1\}.
\end{equation}
This is the midpoint of the time cell $m$ (regardless of directions and stripes).

\subsubsection{Cutoff $b_\theta$}

We shall also need a cutoff function allows smooth transition near the endpoints from stripe to stripe. Consider a parameter 
\begin{equation}\label{theta, def}
\theta := \sqrt{h} 
\end{equation}
for $h=h(\e_0, D)$ chosen in~\eqref{h, choice}. Note that $0 \leq \theta \leq 1/(8\sqrt{2})$. Take a smooth cutoff $b_\theta:[0,1]\to[0,1]$ such that
\begin{align}\label{b-theta}
b_\theta \equiv 0 \text{  on } [0,\theta]\cup[1-\theta,1],\qquad b_\theta \equiv 1 \text{  on } [2\theta,1-2\theta],\qquad \left\|b_\theta'\right\|_{\CC^0([0,1])} \leq \frac{2}{\theta}. 
\end{align} 
Then, the number
\begin{equation}\label{B-theta, def}
B_\theta :=\int_0^1 b_\theta(y)\,\dd y
\end{equation}
satisfies
\begin{align*}
\frac{1}{2} < 1-4\theta\le B_\theta < 1.
\end{align*}
 
\subsubsection{Amplitude $\eta$} 
We now prescribe the amplitude function $\eta(t)$ on the  time slot for the dipole to traverse the stripe $\mathbf{X}_{ij}$ within the $m^{\text{th}}$  cell; see~\eqref{t*, slot} and~\eqref{duration dij}. The ansatz for the amplitude is
\begin{equation}\label{eta on stripe ij}
\eta(t) = \eta_i^{(m)} b_\theta\left( \frac{t-t^{*}_{m,i,j}}{d_{ij}^{(m)}}\right)\qquad\text{on } \left[t^{*}_{m,i,j},\,t^{*}_{m,i,j} + d_{ij}^{(m)}\right] 
\end{equation}
for some constant $\eta_i^{(m)}$. We choose
\begin{equation}\label{eta-i,m, choice}
\eta_i^{(m)}:=\sqrt{\frac{2\Delta_i(t^{**}_m)}{\pi B_\theta} \,}
\end{equation}
where the freezing time \(t^{**}_m\) is as in~\eqref{t**, def}. Also, in view of $B_\theta>1/2$ and the uniform-in-time bound~\eqref{S14} for $\Delta_i$, we have the upper bound: 
\begin{equation}\label{S18}
\eta_i^{(m)}
\le \overline\eta:=
\left[\frac4\pi\left((8+\sqrt2)\delta+\frac34 \e_0^2\right)\right]^{1/2},
\end{equation}
with $\overline{\eta}$ depending on $\delta$ and $\e_0$ only.

In view of the factor $b_\theta$ therein, we can glue the expression for $\eta$ in~\eqref{eta on stripe ij} across all the slots and extend by zero. Thus, from now on, we view $$\eta \in \CC^\infty_c(I).$$

\subsubsection{Core radius $r$}

Recall that $r(t)$ denotes the radius of the core of the travelling dipole. We propose the following ansatz for $r(t)$: 
\begin{equation}\label{r, ansatz}
r(t) := r_{0i}^{(m)}\,a_i\left(t^{**}_m,\cen(t)\right)\qquad\text{on } 
\left[t^*_{m,i,j},\,t^*_{m,i,j}+d_{ij}^{(m)}\right].
\end{equation} 
Recall~\eqref{+1 in time cell} for the time slot in~\eqref{r, ansatz}, as well as the definition of $a_i$ from~\eqref{S11}. Here we set
\begin{equation}\label{rm0i choice}
r_{0i}^{(m)}
:= \bigtime s_i\sqrt{\frac{B_\theta}{8\pi\,\Delta_i(t^{**}_m)}},
\end{equation}
as well as 
\begin{equation}\label{ODE for centre}
\cen'(t)=\frac{\eta(t)}{r(t)}\xi_i\qquad\text{on } 
\left[t^*_{m,i,j},\,t^*_{m,i,j}+d_{ij}^{(m)}\right],
\end{equation}
where $\eta(t)$ has been defined in~\eqref{eta on stripe ij}. To further unburden the notation, write
\begin{align*}
a_i^{(m)}(x) \equiv a_i\big(t^{**}_m, x\big).
\end{align*}

As with $\eta$ above, we may view $r$ as globally defined on $I$. Indeed, between two consecutive stripes, choose a transition interval \([a,b]\) within their combined off portions, with \(a\) in the old slot and \(b\) in the new one. Here, if the slot starts at \(t^*_{M,i,j}\) and has length \(d=d_{ij}^{(m)}\), the ``off portions'' are  $[t^*_{m,i,j},\,t^*_{m,i,j}+\theta d_{ij}^{(m)}]$ and $[t^*_{m,i,j}+(1-\theta)d_{ij}^{(m)},\,t^*_{m,i,j}+d_{ij}^{(m)}]$.  On these intervals \(b_\theta((t-t^*_{m,i,j})/d_{ij}^{(m)})=0\), so \(\eta=\eta'=0\) thereon. If stripe \((m,i,j)\) is followed by stripe \((m',i',j')\), then we set
\[
\begin{aligned}
r_-&=r_{0i}^{(m)}\,a_i^{(m)}\!\left(\stripeij(4L)\right)
&&\text{(radius at the old stripe's exit)},\\
r_+&=r_{0i'}^{(m')}\,a_{i'}^{(m')}\!\left({\mathsf{x}}_{i'j'}(0)\right)
&&\text{(radius at the new stripe's entrance)}.
\end{aligned}
\]
In passing, note by \eqref{S10} that $r_-=8\flf\,r_{0i}^{(m)}$ and $ r_+=8\flf\,r_{0i'}^{(m')}$. We may then set
\[
r(t)=\left(1-\chi\!\left(\frac{t-a}{b-a}\right)\right)r_-
+\chi\!\left(\frac{t-a}{b-a}\right)r_+,
\]
where \(0\le\chi\le1\) is smooth, zero near \(0\), and one near \(1\).

In this way, we obtain the globally defined
\begin{align}\label{r, global}
r \in \CC^\infty(I; ]0,\infty[)\qquad\text{such that  } \sup_{t \in I} \, 2r(t) \leq h
\end{align}

The justification for our choice of the amplitude $\eta$ and the core radius $r$ in \eqref{eta on stripe ij}, \eqref{r, ansatz} will be given in \S\ref{subsec: justification for eta and r} below. There is still one missing piece of the puzzle: an auxiliary force.

\subsection{Auxiliary force}\label{subsec: aux force}


Recall our principal building block for the convex integration scheme: 
\begin{align*}
V(t,x)=\eta(t)U_{r(t)}(x-\cen(t));
\end{align*}
see~\eqref{def, V, travelling dipole}. It satisfies the Euler equation~\eqref{Euler for V with Fc} with the forcing term $F_c$ defined in~\eqref{Fc, def}:
\begin{align*} F_c(t,x):=\eta'(t) W_{r(t)}\big(x-\cen(t)\big)+\eta(t) r(t)'\partial_rW_{r(t)}\big(x-\cen(t)\big),
\end{align*}
where $W_{r(t)}$ is the compact support cutoff of the velocity profile $U_{r(t)}$, as given in~\eqref{Wr, def}. We have proved in 
Lemma~\ref{lem: Wr, Fc} that $$\int_{\R^2} F_c(t,x)\,\dd x=2\pi\frac{d}{dt}\big(\eta(t) r(t)\big) \, \xi,$$
which is nonzero in general. Thus,  \(F_c\) cannot simply be written as the divergence of a compactly supported tensor.

For a remedy, fix a standard bump function $\Omega: \R^2 \to \R$ with
\begin{align*}
\Omega \in \CC^\infty_c\left(\B_1\right),\quad \int_{\R^2}\Omega(x)\,\dd x =1, \quad \text{ $\Omega$ is nonnegative and radial}.
\end{align*}
Define its scaled version as usual:
\begin{align}\label{Omega-h}
\Omega_h(x):= \frac{1}{h^2}\Omega\left(\frac{x}{h}\right)\qquad\text{for } h>0.
\end{align}
We then define the auxiliary force:
\begin{equation}\label{Fh, def}
F_h(t,x)=2\pi(\eta r)'(t)\,\Omega_h\big(x-\cen(t)\big)\,\xi_i.
\end{equation}

\begin{lemma}\label{lem: Fh}
Assume $h>2r(t)$. The difference $(F_h-F_c)(t,\cdot)$ between the forcing terms is compactly supported in $\B_{h}\left(\cen(t)\right)$, and it has zero mean and zero torques, both about $\cen(t)$ and the origin $\mathbf{0}$. That is, for each $t>0$:
\begin{itemize}
\item
$\int_{\R^2}(F_h-F_c)(t,x)\,\dd x =0$. 

\item
$\int_{\R^2}\big(x-\cen(t)\big) \wedge (F_h-F_c)(t,x)\,\dd x  =0$.

\item
$\int_{\R^2} x \wedge (F_h-F_c)(t,x)\,\dd x  =0$.

\end{itemize}

\end{lemma}

\begin{proof}[Proof of Lemma~\ref{lem: Fh}]
Since $\int_{\R^2}\Omega(x)\,\dd x =1$, the mean of $F_h$ is $2\pi(\eta r)'\xi_i$, which equals the mean of $F_c$ by \eqref{eq:scan-F-mean} in Lemma~\ref{lem: Wr, Fc}. This proves the first point.

Since $\Omega_h$ is radial and hence even, we have $\int_{\R^2}z\Omega_h(z)\,\dd x=0$. Thus $F_h$ has zero intrinsic torque about $\cen(t)$. The same holds for $F_c$ by~\eqref{S8} in Lemma~\ref{lem: Wr, Fc}. This proves the second point.

The third point is an immediate consequence of the previous two. 

Finally, the factor $\Omega_h\big(x-\cen(t)\big)$ in $F_h$ yields $\mathrm{supp}(F_h) \Subset \B_h(\cen(t))$, and by Lemma~\ref{lem: Wr, Fc} we have $\mathrm{supp}(F_c) \Subset \B_{2r(t)}(\cen(t))$. The assertion follows from $h>2r(t)$.  \end{proof}

In view of Lemma~\ref{lem: Fh} above, we may obtain $S_c = {\rm div}^{-1}(F_c-F_h)$ by Lemma~\ref{lemma: IO}. More precisely, for each $t>0$, we can solve for $S_c(t,\cdot) \in L^1_x\left(\sym\right)$ such that
\begin{equation}\label{Sc}
\operatorname{div}S_c=F_c-F_h
\end{equation} 
with
\begin{align}\label{Sc, L1 bound}
\|S_c(t,\cdot)\|_{L^1_x}
 &\le C h\left(\|F_c(t,\cdot)\|_{L^1_x}+\|F_h(t,\cdot)\|_{L^1_x}\right)\nonumber\\
  &\le  Ch(|\eta'(t)|r(t)+|\eta(t)|\,|r'(t)|).
\end{align}
Here $C$ are universal constants, and the second line follows from the bound~\eqref{eq:scan-F-bound} for $F_c$ in Lemma~\ref{lem: Wr, Fc} and the definition for $F_h$ in~\eqref{Fh, def}, where the bump function $\Omega_h$ has integral one.

The terms $\eta$, $\eta'$, $r$, and $r'$ will be estimated after we specify the choice of $\eta$ and $r$ in \S\ref{subsec: justification for eta and r}.

\subsection{Justification for the choice of $\eta$ and $r$}\label{subsec: justification for eta and r}
Now let us justify our preceding choice for the amplitude $\eta$ and the core radius $r$ in \eqref{eta on stripe ij} and \eqref{r, ansatz}. We shall follow three simple rules:
\begin{itemize}
\item[(i)]
The center of the travelling dipole, $\cen(t)$, traverses precisely one stripe ${\mathbf{X}}_{ij}$ in the time slot $\left[t^*_{m,i,j},\,t^*_{m,i,j}+d_{ij}^{(m)}\right]$.

\item[(ii)]
For each $i \in \{1,2,3,4\}$, the dipole takes a quarter of the cell time, $\bigtime/4$, to traverse all the stripes $\left\{{\mathbf{X}}_{ij}: 1\leq j \leq N_i\right\}$ oriented towards $\xi_i$.

\item[(iii)]
The stripes ${\mathbf{X}}_{ij}$ are evenly spaced in $Q_i$, with a separation of $s_i$ between adjacent stripes. In the time-averaged approximation of $R_c$ (the spatially compact cutoff of the Reynolds stress tensor), each stripe should carry the correct weight up to the controlled ``ramp error'' (caused by the portions where $b_\theta \neq 1$).

\end{itemize}

\subsubsection{Rule~$(i)$}
Note by~\eqref{ODE for centre} that the scaled dipole with amplitude \(\eta\) translates at speed \(\eta/r\).  If \(\lambda\) is arclength along the stripe ${\mathbf{X}_{ij}}$ in time cell $m$, then 
\begin{equation*}
\frac{d\lambda}{dt}(t)=\frac{\eta(t)}{r_{0i}^{(m)}a_i(t^{**}_m,\cen(t))}
\end{equation*}
on $\left[t^*_{m,i,j},\,t^*_{m,i,j}+d_{ij}^{(m)}\right]$. Here $t^{**}_m$ is the midpoint of the time cell $m$, defined as in \eqref{t**, def}. We naturally put 
\begin{align*}
\lambda(t^*_{m,i,j})=0, \qquad
\lambda\left(t^*_{m,i,j}+d_{ij}^{(m)}\right)=4L,
\end{align*}
and extend $\lambda$ continuously to $\R_+$ by constants $0$ and $4L$ outside $[0,4L]$.

As in \eqref{parametrisation, stripe}, we denote by $\stripeij(\lambda)$ with $0\le\lambda\le4L$ the arclength parametrisation of the stripe $\mathbf{X}_{ij}$. Then the integral of $a_i(t^{**}_m,\cen(t))$ along the stripe becomes
\begin{align*}
G_{ij}(\lambda) := \int_0^\lambda a_i\big(t^{**}_{m},\, \stripeij(s) \big)\,\dd s,
\end{align*} 
where $\cen(t)= \stripeij\big(\lambda(t)\big)$. Here $\stripeij'(\lambda)=\xi_i$ and $G_{ij}(4L)=I_{ij}$ as in~\eqref{S13}, where we suppress the superscript $m$. As $\eta \equiv 0$ near the endpoints of stripes, we may compute by the chain rule and fundamental theorem of calculus:
\begin{align*}
\frac{d}{dt}G_{ij}\big(\lambda(t)\big)
=a_i\big(t^{**}_m,\cen(t))\frac{d\lambda}{dt}(t)
=\frac{\eta(t)}{r_{0i}^{(m)}} \qquad\text{on } 
\left[t^*_{m,i,j},\,t^*_{m,i,j}+d_{ij}^{(m)}\right].
\end{align*}
Thus, by integration and the change of variables $y=\frac{t-t^*_{m,i,j}}{d_{ij}^{(m)}}$, we arrive at
\begin{align}\label{condition I}
I_{ij}
&=\frac1{r_{0i}^{(m)}}\int_{t^{*}_{m,i,j}}^{t^{*}_{m,i,j}+d^{(m)}_{ij}}\eta(t)\,\dd t
\nonumber\\
&=\frac{\eta_i^{(m)} d^{(m)}_{ij}}{r_{0i}^{(m)}}\int_0^1 b_\theta(y)\,\dd y
=\frac{\eta_i^{(m)} B_\theta d^{(m)}_{ij}}{r_{0i}^{(m)}}.
\end{align}
 
 Condition~\eqref{condition I}  is sufficient and necessary: \(G_{ij}\) is strictly increasing as \(a_i>0\), so integrating the right-most expression in~\eqref{condition I}  and applying \(G_{ij}^{-1}\) gives a trajectory that reaches \(4L\) exactly.
 
 \subsubsection{Rule~$(ii)$} To this end, we simply impose 
 \begin{align}\label{condition II}
 \sum_{j=1}^{N_i} d^{(m)}_{ij} = \frac{\bigtime}{4}\qquad\text{for each } m \in \{0,1,\ldots,M-1\} \text{ and } i \in \{1,2,3,4\}. 
 \end{align}
Here, $d^{(m)}_{ij}$ is the time of the dipole spent on the stripe ${\mathbf{X}_{ij}}$ in the $m^{\text{th}}$ time cell; see~\eqref{duration dij}. Also, $\bigtime$ as in~\eqref{timelength} is the duration of each time cell. 
 
\begin{figure}[htbp]
\centering
\includegraphics[width=1\linewidth]{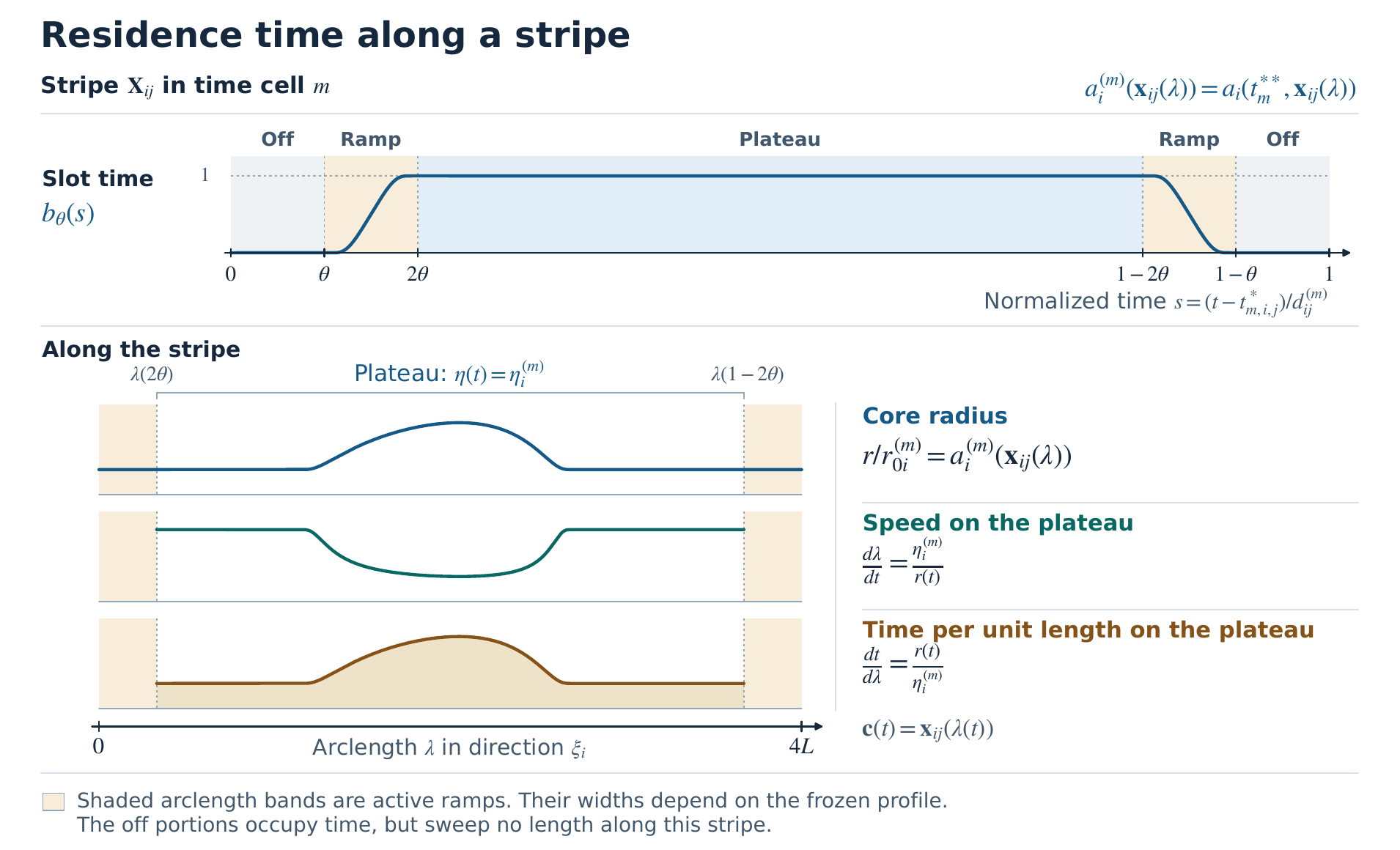}
\caption{Schematic profiles along stripe $\mathbf{X}_{ij}$ in time cell $m$.
The coefficient $a_i^{(m)}(x)=a_i(t_m^{**},x)$ is evaluated at the
moving centre $\cen(t)=\mathbf{x}_{ij}(\lambda(t))$, and
$r(t)=r_{0i}^{(m)}a_i^{(m)}(\cen(t))$.
On the plateau $\eta(t)=\eta_i^{(m)}$, the speed is
$d\lambda/dt=\eta_i^{(m)}/r(t)$ and the time per unit arclength is
$dt/d\lambda=r(t)/\eta_i^{(m)}$.
The shaded bands mark the ramp portions (namely $\{b_\theta \neq 1\}$) schematically.}
\label{fig:residence}
\end{figure}

 \begin{remark}\label{remark: reset}
One caveat is in order:  We interpolate the centre and core radius during the final ``off portion'' of one time slot and the initial off portion of the next. Since \(\eta=\eta'=0\) thereon, this interpolation does not affect the forcing and other relevant perturbations. No additional time is allocated to the ``\emph{reset}'' of the centre of the travelling dipole from one stripe to another. See \S\ref{subsubsec: reset} for details. Here and throughout, by ``\emph{off portion}'' we mean the two subintervals on $\left[t^*_{m,i,j},\,t^*_{m,i,j}+d_{ij}^{(m)}\right]$, one ``initial'' and one ``final'', on which $b_\theta \equiv 0$; recall $b_\theta$ from~\eqref{b-theta}.

On the plateau $\{b_\theta = 1\}$, it holds that 
\begin{align*}
\frac{dt}{d\lambda}=\frac{r_{0i}^{(m)}}{\eta_i^{(m)}}a_i^{(m)}.
\end{align*} 
Thus, the time
spent per unit length is proportional to the frozen coefficients $a^{(m)}_i$ at $t^{**}_m$, the midpoint of the $m^{\text{th}}$ time cell. See Figure~\ref{fig:residence}. 

 \end{remark}

 \subsubsection{Rule~$(iii)$}
This is the most technical part, which connects the forcing terms to the geometric configuration of the stripes. First of all, recall~\eqref{S11}:
\begin{equation*}
-R_c(t,x) = \sum_{i=1}^4\Big(a_i(t,x)\xi_i\otimes\xi_i\Big) - b(t,x){\bf I}_2.
\end{equation*}
Hence, controlling the $L^1_x$-norm of $R_c$ amounts to estimating the integral of $a_i$ modulo a pressure term. Note that $a_i$ is strictly positive: $a_i \geq \sqrt{62}\flf$ by~\eqref{pointwise bound for ai}.

To this end, we fix cell \(m\) and direction \(i\) and consider the coordinates
\begin{equation}\label{Yi, coord}
Y_i(\lambda,z)=\ci+\lambda\xi_i+z\xi_i^\perp,
\qquad 0\le \lambda,z\le 4L,
\end{equation}
where $\ci$ is the corner of square $Q_i$ in~\eqref{corner}. Recall from~\eqref{si, footstep} that   $s_i:= \tfrac{4L}{N_i}$ is the spacing between two consecutive stripes in $Q_i$. Thus, in the $z$-variable in \eqref{Yi, coord}, the $j^{\text{th}}$ ``transverse band'' is $$\band_{ij}=[(j-1)s_i,js_i],$$ whose midpoint is \(z_j=(j-\tfrac12)s_i\) is the natural coordinate for the $j^{\text{th}}$ stripe $\mathbf{X}_{ij}$; \emph{i.e.},
\begin{equation*}
Y_i(\lambda, z_j) \equiv \stripeij(\lambda)\qquad\text{for each } \lambda \in [0,4L].
\end{equation*}

To proceed, notice by~\eqref{Omega-h} that $\Omega_h$ is a standard approximation of identity. We shall thus focus on the estimate for the convolution $\Omega_h \star \left(a_i^{(m)}\mathbf1_{Q_i}\right)$, where $a_i^{(m)}(x) \equiv a_i\left(t^{**}_m, x\right)$ is evaluated at the freezing time $t^{**}_m$, the centre of the $m^{\text{th}}$ time cell. We have
\begin{equation}\label{convolution, equal}
\left(\Omega_h \star \left(a_i^{(m)}\1_{Q_i}\right)\right)(x)
=
\sum_{j=1}^{N_i}\int_{\band_{ij}}\int_0^{4L}
a_i^{(m)}(Y_i(\lambda,z))\,
\Omega_h\big(x-Y_i(\lambda,z)\big)\,\dd \lambda\,\dd z.
\end{equation}
Its quadrature evaluated at $z =z_j$ is then 
\begin{align}\label{widetilde-Omega}
\widetilde{\left(\Omega_h \star \left(a_i^{(m)}\1_{Q_i}\right)\right)}(x) &:= \sum_{j=1}^{N_i}\int_{\band_{ij}}\int_0^{4L}
a_i^{(m)}(Y_i(\lambda,z_j))\,
\Omega_h\big(x-Y_i(\lambda,z_j)\big)\,\dd \lambda\,\dd z \nonumber\\
&= s_i\sum_{j=1}^{N_i}\int_0^{4L}
a_i^{(m)}(\stripeij(\lambda))\,
\Omega_h(x-\stripeij(\lambda))\,\dd\lambda,
\end{align}
thanks to $\int_{\band_{ij}} 1\,\dd z = s_i$ and $Y_i(\lambda, z_j) \equiv \stripeij(\lambda)$.

Now, let us recall the auxiliary force $F_h$ defined in~\eqref{Fh, def}. We compute the time-averaged auxiliary force over the $m^{\text{th}}$ time cell:
\begin{align}\label{Fh-bar, def}
\overline F_h^{(m)}(x)
&:= \frac{1}{\bigtime}\int_{t_I^- + m\bigtime}^{t_I^- + (m+1)\bigtime} F_h(t,x)\,\dd t \nonumber\\
&\equiv \frac{1}{\bigtime} \sum_{i=1}^4 \sum_{j=1}^{N_i}\int_{t^*_{m,i,j}}^{t^*_{m,i,j}+d_{ij}^{(m)}} F_h(t,x)\,\dd t.
\end{align}
It holds that
\begin{align}\label{a1}
\overline F_h^{(m)}(x) &= \frac{2\pi}{\bigtime}\sum_{i=1}^4 \sum_{j=1}^{N_i}\int_{t^*_{m,i,j}}^{t^*_{m,i,j}+d_{ij}^{(m)}}
   \partial_t\!\bigl(\eta(t)r(t)\bigr)
   \Omega_h\!\left(x-\cen(t)\right)\xi_i\,\dd t\nonumber\\
&=-\frac{2\pi}{\bigtime}\sum_{i=1}^4 \sum_{j=1}^{N_i}\int_{t^*_{m,i,j}}^{t^*_{m,i,j}+d_{ij}^{(m)}}
   \eta(t)r(t)\,
   \partial_t\!\left[\Omega_h\!\left(x-\cen(t)\right)\right]
   \xi_i\,\dd t\nonumber\\
&=\frac{2\pi}{\bigtime}\sum_{i=1}^4 \sum_{j=1}^{N_i}\int_{t^*_{m,i,j}}^{t^*_{m,i,j}+d_{ij}^{(m)}}
   \eta(t)r(t)
   \left[\p_t\cen(t)\cdot
   \nabla_x\Omega_h\!\left(x-\cen(t)\right)\right]
   \xi_i\,\dd t\nonumber\\
&=\operatorname{div}_x\!\left[
   \frac{2\pi}{\bigtime}\sum_{i=1}^4 \sum_{j=1}^{N_i}\int_{t^*_{m,i,j}}^{t^*_{m,i,j}+d_{ij}^{(m)}}
   \eta(t)^2\Omega_h\!\left(x-\cen(t)\right)
   \xi_i\otimes\xi_i\,\dd t
   \right].
\end{align}
The boundary term vanishes because $\eta(t)=0$ at both endpoints
of $\left[t^*_{m,i,j},\,t^*_{m,i,j}+d_{ij}^{(m)}\right]$. The last equality uses $r(t)\dot \cen(t)=\eta(t)\xi_i$. 

By passing to the arclength variable, we infer from $r\p_t\lambda=\eta$ that
\begin{align}\label{a2}
\eta(t)^2\,\dd t
&=\eta(t)r(t)\,\dd\lambda\nonumber\\
&=\eta_i^{(m)}r_{0i}^{(m)}\,
b_\theta\!\left(\frac{t-t^*_{m,i,j}}{d_{ij}^{(m)}}\right)
a_i^{(m)}(\stripeij(\lambda))\,\dd\lambda\qquad\text{on } \left[t^*_{m,i,j},\,t^*_{m,i,j}+d_{ij}^{(m)}\right].
\end{align}
See~\eqref{eta on stripe ij}, \eqref{eta-i,m, choice}, \eqref{r, ansatz}, \eqref{rm0i choice}, and \eqref{b-theta} for the terms on the right-most side. Therefore, substituting~\eqref{a2}  into~\eqref{a1} and writing $t=t(\lambda)$ implicitly on $\left[t^*_{m,i,j},\,t^*_{m,i,j}+d_{ij}^{(m)}\right]$, we arrive at
\begin{align}\label{a3}
\overline F_h^{(m)}(x)
 = \operatorname{div}_x\!\left(
\sum_{i=1}^4 A_i^{(m)}\,\xi_i\otimes\xi_i\right),
\end{align}
where
\begin{equation}\label{a4}
A_i^{(m)}(x)
=\frac{2\pi\eta_i^{(m)}r_{0i}^{(m)}}{\mathsf T_{\mathrm{cell}}}
\sum_{j=1}^{N_i}\int_0^{4L}
b_\theta\!\left(\frac{t(\lambda)-t^*_{m,i,j}}{d_{ij}^{(m)}}\right)
a_i^{(m)}(\stripeij(\lambda))\,
\Omega_h\!\left(x-\stripeij(\lambda)\right)\,\dd\lambda.
\end{equation}
To make $t(\lambda)$ well defined, we may further require $b_\theta > 0$ on $]\theta,1-\theta[$.

Comparing~\eqref{a3} and \eqref{a4} with~\eqref{widetilde-Omega} suggests that
\begin{equation}\label{condition III}
\frac{2\pi \eta_i^{(m)}r_{0i}^{(m)}}{\bigtime} = s_i.
\end{equation}


\subsubsection{Summary} Putting together the three conditions~\eqref{condition I}, \eqref{condition II}, and \eqref{condition III}, we obtain the unique choice for the coefficients $\eta_i^{(m)}$, $r_{0i}^{(m)}$, and $d_{ij}^{(m)}$, for each stripe $\mathbf{X}_{ij}$ in $Q_i$ with $j \in \{0,1,\ldots,N_i\}$, each direction $\xi_i; i\in\{1,2,3,4\}$, and each time cell labelled by $m\in \{0,1,\ldots, M-1\}$:
\begin{align}\label{S15}
\eta_i^{(m)}
  &=\sqrt{\frac{2\Delta_i^{(m)}}{\pi B_\theta}\,},\nonumber\\
r_{0i}^{(m)}
  &=\frac{\bigtime\,\eta_i^{(m)}B_\theta}{4 \sum_{j=1}^{N_i}I_{ij}^{(m)}},\\
d_{ij}^{(m)}
  &=\frac{\bigtime}{4}\frac{I_{ij}^{(m)}}{ \sum_{j'=1}^{N_i}I_{ij'}^{(m)}}
   =\frac{r_{0i}^{(m)}I_{ij}^{(m)}}{\eta_i^{(m)}B_\theta}.\nonumber
\end{align}

Here, as in~\eqref{S13}, \(I_{ij}^{(m)}\) is the integral of the frozen coefficient \(a_i^{(m)}(x)=a_i(t_m^{**},x)\) along stripe \(\mathbf X_{ij}\) with respect to arclength:
\begin{equation}\label{Iij,m}
I_{ij}^{(m)}
:=\int_{\mathbf X_{ij}} a_i(t_m^{**},x)\,\dd\mathcal H^1(x)
=\int_0^{4L}
a_i\!\left(t_m^{**},Y_i\!\left(\lambda,(j-\tfrac12)s_i\right)\right)\,\dd\lambda.
\end{equation}

\subsubsection{Reset}\label{subsubsec: reset}
For the sake of bookkeeping, we summarise what happens between the traversal of two consecutive stripes by the travelling dipole in our convex integration construction.

The centre $\cen(t)$ of the principal building block travels along stripes $\mathbf{X}_{ij}$ in the $m^{\text{th}}$ time cell. Its velocity profile is
$$V(t,x)=\eta(t)U_{r(t)}(x-\cen(t)),$$ with $$\cen'(t)=\frac{\eta(t)}{r(t)}\xi_i \qquad \text{for } t \in \left[t^*_{m,i,j},\,t^*_{m,i,j}+d_{ij}^{(m)}\right],$$ as in~\eqref
{def, V, travelling dipole} and \eqref{c'}. Here, $U_{r(t)}$ is the scaled Lamb--Chaplygin dipole introduced in \S\ref{subsec: dipole}.

In view of the structure of the amplitude function $\eta(t)$ in ~\eqref{eta on stripe ij} --- in particular, the presence of the cutoff function $b_\theta$ as in~\eqref{b-theta} --- the profile $V(t,\cdot)$, as well as the forcing terms $F_c$ and $F_h$ given in \eqref{Fc, def} and \eqref{Fh, def}, respectively, are all identically zero near the endpoint of each stripe's time interval $\timeint$.

For the radius function $r(t)$: One may choose an arbitrary smooth connecting path for
the inactive center, agreeing with the endpoint on the previous stripe before the ``reset'' (see Remark~\ref{remark: reset}) and the starting point of the next stripe after the reset. The radius function can be joined smoothly by a positive interpolation between its endpoint values, so that the same radius bound holds during the reset. See the paragraph right above~\eqref{r, global} for the detailed construction. This does not alter $V$, $F_c$, or $F_h$, as $\eta=\eta'=0$ on the connecting interval.

The definition of the center $\cen(t)$ in~\eqref{c'} is not imposed on that connecting path: imposing it there would incorrectly require $\cen'=0$. Indeed, instead of glueing the centre function $\cen(t)$ for one dipole, we use separate centre functions on different stripes $\mathbf{X}_{ij}$ and glue only the resulting velocity building blocks as entailed in the preceding paragraph. Change of directions $\xi_i$ abides by the same convention. There is at most one active principal dipole at a time.

\subsection{Key estimates for forcing and stress quadratures}

Given the choices for $\eta$ and $r$ specified in~\S\ref{subsec: justification for eta and r} above, we are ready to bound the forcing terms $F_c$, $F_h$, as well as $S_c$ and the approximation error or $A_i^{(m)}$, in terms of $D$ and $h$ only. These estimates are summarised in Lemma~\ref{lem: estimates given eta, r} below.

Recall from~\eqref{a3} the expression of the time-averaged auxiliary force:  
 \begin{align*}
\overline F_h^{(m)}(x)
 = \operatorname{div}_x\!\underbrace{\left( \sum_{i=1}^4 A_i^{(m)}\,\xi_i\otimes\xi_i\right)}_{=:\biga}.
\end{align*} 
Recall also from \eqref{Delta-i, quadrature, def} and
\eqref{Iij,m} the quadrature $\Delta_i^{(m)}$ that approximates the area integral of $a_i$ over the square $Q_i$ in $m^{\text{th}}$ time cell: $$\Delta_i^{(m)}=s_i\sum_j I_{ij}^{(m)}
=\Delta_i(t_m^{**}).$$
Furthermore, recall~\eqref{widetilde-Omega}: \begin{align*}
\widetilde{\left(\Omega_h \star \left(a_i^{(m)}\1_{Q_i}\right)\right)}(x) &:= \sum_{j=1}^{N_i}\int_{\band_{ij}}\int_0^{4L}
a_i^{(m)}(Y_i(\lambda,z_j))\,
\Omega_h\big(x-Y_i(\lambda,z_j)\big)\,\dd \lambda\,\dd z \nonumber\\
&= s_i\sum_{j=1}^{N_i}\int_0^{4L}
a_i^{(m)}(\stripeij(\lambda))\,
\Omega_h(x-\stripeij(\lambda))\,\dd\lambda.
\end{align*}

For ease of future reference, we introduce the notation:
\begin{equation}\label{gi, def}
g_i^{(m)}(x):=a_i^{(m)}(x)\,\1_{Q_i}(x)
=a_i\big(t_m^{**},x\big)\,\1_{Q_i}(x).
\end{equation}
The term $A_i^{(m)}$ (see~\eqref{a4}, reproduced below) has been designed to approximate $g_i^{(m)}$:
\begin{align*}
A_i^{(m)}(x)
=\frac{2\pi\eta_i^{(m)}r_{0i}^{(m)}}{\mathsf T_{\mathrm{cell}}}
\sum_{j=1}^{N_i}\int_0^{4L}
b_\theta\!\left(\frac{t(\lambda)-t^*_{m,i,j}}{d_{ij}^{(m)}}\right)
a_i^{(m)}(\stripeij(\lambda))\,
\Omega_h\!\left(x-\stripeij(\lambda)\right)\,\dd\lambda.
\end{align*}

\begin{lemma}\label{lem: estimates given eta, r}
Fix a time cell \(m\), a direction \(i\), and a stripe \(\mathbf{X}_{ij}\).
\begin{itemize}
\item[(i)]
Both $A^{(m)}_i$ and \(\overline F_h^{(m)}\) are smooth and spatially
supported in \(\overline \B_{2\sqrt{2}L+h}(0)\).
\item[(ii)]
The following holds for $\eta$ and $r$:
\begin{align}\label{S19}
&|\eta'(t)|\,r(t)
\le
\frac{2\eta_i^{(m)}r_{0i}^{(m)}
      \left\|a_i^{(m)}\right\|_{L^\infty(Q_i)}}
     {\theta d_{ij}^{(m)}}
=
\frac{4\left(\eta_i^{(m)}\right)^2 B_\theta
      \left\|a_i^{(m)}\right\|_{L^\infty(Q_i)}}
     {\theta I_{ij}^{(m)}},\nonumber
\\
&|\eta(t)|\,|r'(t)|
=
\eta(t)^2
\left|
  \left.\xi_i\cdot\nabla_x\log a_i^{(m)}(x)
  \right|_{x=\cen(t)}
\right|.
\end{align}
In particular, $I_{ij}^{(m)}$ has a uniform positive lower bound.

\item[(iii)]
The forcing terms satisfy that
\begin{equation*}
\|F_c(t,\cdot)\|_{L^1_x} + \|F_h(t,\cdot)\|_{L^1_x} \leq \frac{D}{\theta}
\end{equation*}
for all $t \in J$, provided that $D$ is sufficiently large but keeping the same dependency on parameters as in~\eqref{D}.

\item[(iv)] By further enlarging $D$ in (iii) above, the correction term $S_c$ for the Reynolds stress tensor given in~\eqref{Sc} can be bounded by
\begin{equation}\label{Sc-bound-good}
\|S_c(t,\cdot)\|_{L^1_x} \leq D \theta.
\end{equation}

\item[(v)] $A_i^{(m)}$ nicely approximates $g_i^{(m)}$ in the spatial $L^1$-norm:
\begin{align}\label{S25}
 \left\|A_i^{(m)}-g_i^{(m)}\right\|_{L^1_x}
 &\le 8\theta\,\Delta_i^{(m)}
   +s_i\left\| \xi_i^\perp \cdot \nabla a_i^{(m)}\right\|_{L^1(Q_i)}
   +\frac{s_i}{2h}\left\|\nabla\Omega\right\|_{L^1_x}\left\|g_i^{(m)}\right\|_{L^1_x} \nonumber \\
 &\quad+ h\left(\int |y|\Omega(y)\,\dd y\right)
 \left|D_x g_i^{(m)}\right|(\mathbb R^2).
\end{align}

\end{itemize}
\end{lemma}

As will become transparent along the proof, for Point~$(v)$, any pointwise sampling will be applied to the smooth \(a_i^{(m)}\) on the square $Q_i$. The \(BV\)-estimates, which accounts for the term~$ \left|D_x g_i^{(m)}\right|(\mathbb R^2)$ in~\eqref{S25},  concerns only the spatial convolution
error $\left\|\Omega_h\star g_i^{(m)} - g_i^{(m)}\right\|_{L^1_x}$. In particular, no pointwise evaluation of an arbitrary \(BV\)
representative is required.

\begin{proof}[Proof of Lemma~\ref{lem: estimates given eta, r}]
The first point follows from  \eqref{square Qi},
\eqref{parametrisation, stripe}, and \eqref{Omega-h}. Notice too that the same support
bound holds for \(F_h(t,\cdot)\) at all times in \(J\).

For~$(ii)$, recall from the pointwise lower bound~\eqref{pointwise bound for ai} for $a_i$ and the choice of $\flf$ in~\eqref{floor} that
\begin{equation}\label{Iijm, lower bound}
I_{ij}^{(m)} \geq 4\sqrt{62} \flf L = \frac{\sqrt{62}\e_0^2}{256 L}.
\end{equation}
Also note that on a reset interval, the global centre and radius can be interpolated without satisfying the stripe equations, but both products on the left of~\eqref{S19} vanish there because \(\eta(t)=\eta'(t)=0\). Hence, we only need to argue over the stripes $\mathbf{X}_{ij}$.

Write $y(t):=\frac{t-t^*_{m,i,j}}{d_{ij}^{(m)}}$ as before. Differentiation  yields $\eta'(t)
=\frac{\eta_i^{(m)}}{d_{ij}^{(m)}}\,b_\theta'(y(t))$, so
\[
|\eta'(t)|\,r(t)
=
\frac{\eta_i^{(m)}r_{0i}^{(m)}}{d_{ij}^{(m)}}
\,\big|b_\theta'(y(t))\big|\,a_i^{(m)}\big(\cen(t)\big).
\]
By~\eqref{b-theta}, \(\|b_\theta'\|_{L^\infty(0,1)}\le 2/\theta\). Meanwhile, \(\cen(t)\in\overline{Q_i}\) for $t \in \timeint$, so by
continuity \(a_i^{(m)}(\cen(t))\le \left\|a_i^{(m)}\right\|_{L^\infty(Q_i)}\). These two bounds prove the inequality in the first line of~\eqref{S19}. The equality in this line follows
from the last formula for $d_{ij}^{(m)}$ in \eqref{S15}. For the second line, we compute
\[
\begin{aligned}
r'(t)
&=r_{0i}^{(m)}
  \left.\nabla_x a_i(t_m^{**},x)\right|_{x=\cen(t)}
   \cdot\dot \cen(t)\\
&=\eta(t)\,
  \frac{\left.\xi_i\cdot\nabla_x a_i^{(m)}(x)\right|_{x=\cen(t)}}
       {a_i^{(m)}(\cen(t))}
 =\eta(t)
  \left.\xi_i\cdot\nabla_x\log a_i^{(m)}(x)
  \right|_{x=\cen(t)} .
\end{aligned}
\]
We used the centre equation, the formula for \(r(t)\), and the
positivity of \(a_i^{(m)}\) from \eqref{pointwise bound for ai}.
Taking absolute values and multiplying by \(|\eta(t)|\) proves
the second point.

For the third point, substituting the expression~\eqref{S15} for $\eta_i^{(m)}$, we deduce that
\begin{align*}
\frac{2\left(\eta_i^{(m)}\right)^2 B_\theta
      \left\|a_i^{(m)}\right\|_{L^\infty(Q_i)}}
     {\theta I_{ij}^{(m)}} = \frac{4\Delta_i^{(m)} \left\|a_i^{(m)}\right\|_{L^\infty(Q_i)}}{\pi\theta I^{(m)}_{ij}}.
\end{align*}
With the upper bounds for $\Delta_i^{(m)}$ and $a_i^{(m)}$ given in~\eqref{S14} and~\eqref{pointwise bound for ai}, respectively, as well as the lower bound for $I^{(m)}_{ij}$ given in~\eqref{Iijm, lower bound}, we infer that $|\eta'(t)|r(t) \lesssim \theta^{-1}$ modulo a constant with the same dependency on parameters as for $D$. The same argument holds for $|\eta(t)||r'(t)|$, thanks to the second line in~\eqref{S19}. Notice too that $\|F_c(t,\cdot)\|_{L^1_x}\lesssim (|\eta'(t)|r(t)+|\eta(t)|\,|r'(t)|)$ modulo a universal constant by \eqref{eq:scan-F-bound}. Moreover, for $F_h$, in view of its definition~\eqref{Fh, def} and that the nonnegative kernel $\Omega_h$ verifies $\int_{\R^2}\Omega_h(x)\,\dd x =1$, one has $\|F_h(t,\cdot)\|_{L^1_x} \lesssim (|\eta'(t)|r(t)+|\eta(t)|\,|r'(t)|)$, again modulo a universal constant. This proves~$(iii)$. 

Next, in view of Lemma~\ref{lem: Fh} above, $F_c-F_h$ is compactly supported in the spatial variable and has zero mean and zero torque at every time. This enables us to apply  Lemma~\ref{lemma: IO} to obtain $S_c(t,\cdot) \in L^1_x\left(\sym\right)$ such that $\operatorname{div}S_c=F_c-F_h$ with \begin{align*}
\|S_c(t,\cdot)\|_{L^1_x}
 \lesssim h(|\eta'(t)|r(t)+|\eta(t)|\,|r'(t)|)
\end{align*}
modulo a universal constant; see~\eqref{Sc, L1 bound}. By Point~$(iii)$ above, the right-hand side $\lesssim Dh/\theta$. The proof is complete once we invoke~\eqref{theta, def}, \emph{i.e.}, $\theta =\sqrt{h}$. This proves~$(iv)$.

It remains to show $(v)$. We proceed by triangle inequality to bound
\begin{align*}
\left\|A_i^{(m)}-g_i^{(m)}\right\|_{L^1_x} &\leq \left\|A_i^{(m)} - \widetilde{\left(\Omega_h\star g_i^{(m)}\right)}\right\|_{L^1_x} \\
&\qquad  + \left\|\widetilde{\left(\Omega_h\star g_i^{(m)}\right)} - \left(\Omega_h\star g_i^{(m)}\right)\right\|_{L^1_x} +  \left\|\Omega_h\star g_i^{(m)} - g_i^{(m)}\right\|_{L^1_x}. 
\end{align*}
Here we recall the notation in~\eqref{widetilde-Omega}: for an integral quantity~$\mathscr{I} = \int \{\cdots\}\,\dd z$, we write $\widetilde{\mathscr{I}}$ for the quantity obtained from $\mathscr{I}$ by replacing the $z$-variable in the integrand, which is integrated over $\band_{ij}=[(j-1)s_i,js_i]$, by the midpoint $z_j = (j-1/2)s_j$.

\noindent
\underline{Term~1.} On each time slot $\timeint$, we shift to the scaled variable~$y=\frac{t-t^*_{m,i,j}}{d_{ij}^{(m)}}$. The arclength parameter \(\lambda=\lambda(y)\) is nondecreasing from \(0\) to \(4L\). It follows from ~\eqref{condition III}, \eqref{widetilde-Omega}, the definition of $B_\theta$ in~\eqref{B-theta, def}, and the bounds for $b_\theta$, $B_\theta$ that
\begin{align}\label{first, difference est}
\left\|A_i^{(m)}
  -\widetilde{\left(\Omega_h\star g_i^{(m)}\right)}\right\|_{L^1_x}
=\frac{\Delta_i^{(m)}}{B_\theta}
  \int_0^1 b_\theta(y)(1-b_\theta(y))\,\dd y \le 8\theta\Delta_i^{(m)}.
\end{align}
In particular, $$\left\|A^{(m)}_i\right\|_{L^1_x}\leq  \Delta_i^{(m)}.$$

\noindent
\underline{Term~2.} Recall from~\eqref{Yi, coord} the coordinate functions \(Y_i\)
on $Q_i$. For \(z\in\band_{ij}\), the integrands in \eqref{widetilde-Omega} and
\eqref{convolution, equal} differ by
\begin{align}\label{difference, x}
&a_i^{(m)}\Big(Y_i(\lambda,z_j)\Big)\Omega_h\Big(x-Y_i(\lambda,z_j)\Big)-a_i^{(m)}\Big(Y_i(\lambda,z)\Big)\Omega_h\Big(x- Y_i(\lambda,z) \Big)\nonumber\\
 &\qquad=\bigg[a_i^{(m)}\Big(Y_i(\lambda,z_j)\Big)-a_i^{(m)}\Big(Y_i(\lambda,z)\Big)\bigg]\,\Omega_h\Big(x-Y_i(\lambda,z_j)\Big)\nonumber\\
  &\qquad\qquad   +a_i^{(m)}\Big(Y_i(\lambda,z)\Big)\bigg[\Omega_h\Big(x-Y_i(\lambda,z_j)\Big)-\Omega_h\Big(x-Y_i(\lambda,z)\Big)\bigg],
\end{align}
where $z_j = (j-1/2)s_j$. Using $\int_{\R^2}\Omega_h(x)\,\dd x=1$ and the fundamental theorem of calculus as for the derivation of~\eqref{Delta-i, quadrature, def}, the first term on the right-hand side of~\eqref{difference, x} is controlled by
\[
 s_i\int_{\band_{ij}}
  \left|\partial_{\xi_i^\perp}a_i^{(m)}\big(Y_i(\lambda,z)\big)\right|\,\dd z.
\]
Here it is crucial that the consecutive stripes in $Q_i$ are separated by length $s_i$. On the other hand, standard translation estimates for the convolution kernel $\Omega_h$ leads to
\begin{align*}
&\left\|\Omega_h\Big(\cdot-Y_i(\lambda, z_j)\Big)-\Omega_h\Big(\cdot-Y_i(\lambda, z_j)\Big)\right\|_{L^1_x} \\
&\qquad \le |z-z_j| \|\nabla\Omega_h\|_{L^1_x} \le \frac{s_i}{2h}\|\nabla\Omega\|_{L^1_x},
\end{align*}
which controls the final term in~\eqref{difference, x}. As $J_i$ has Jacobian determinant $1$, we integrate over $\lambda$ and $z$ in each $\band_{ij}$ to conclude that
\begin{align}\label{second, difference est}
\left\|\widetilde{\left(\Omega_h\star g_i^{(m)}\right)} - \left(\Omega_h\star g_i^{(m)}\right)\right\|_{L^1_x} &\leq s_i\left\|\partial_{\xi_i^\perp}a_i^{(m)}\right\|_{L^1(Q_i)}
      +\frac{s_i}{2h}\left\|\nabla\Omega\right\|_{L^1_x}\left\|g_i^{(m)}\right\|_{L^1_x}. 
\end{align}

\noindent
\underline{Term~3.} We shall take into account the jump of \(g_i^{(m)}\) across \(\partial Q_i\).  Note that for any \(g\in BV(\mathbb R^2)\cap L^1(\R^2)\) we have \(\|g(\cdot-y)-g\|_{L^1_x}\le |y|\,|D_xg|(\mathbb R^2)\). Hence, given $\Omega_h$ as in~\eqref{Omega-h}, we deduce that
\begin{align}\label{third, difference est}
\left\|\Omega_h\star g_i^{(m)}-g_i^{(m)}\right\|_{L^1_x}
 \le h\left(\int |y|\Omega(y)\,\dd y\right)
             \left|D_x g_i^{(m)}\right|\left(\mathbb R^2\right).
\end{align}

\noindent
\underline{Conclusion}
The estimate~\eqref{S25}, \emph{i.e.}, Point~(v) of the lemma, follows from \eqref{first, difference est}, \eqref{second, difference est}, and \eqref{third, difference est}.

 The proof is now complete. \end{proof}

\subsection{Time corrector}
\label{subsec:scan-time-corrector}

In the above constructions, on each of the time cells we have approximated the forcing $F_h$ by $\overline{F}^{(m)}_h = \operatorname{div}_x\!\left(
\sum_{i=1}^4 A_i^{(m)}\,\xi_i\otimes\xi_i\right)$, the sum of four symmetric tensors with time-independent coefficients $A_i^{(m)}$ evaluated at frozen time points~\(t_m^{**}\) in \eqref{t**, def}. In this subsection, we shall estimate the temporal fluctuation $F_h(t,\cdot) - \overline{F}^{(m)}_h(\cdot)$ on the $m^{\text{th}}$ time cell: it has zero time average but need not be small at any particular time.

For notational convenience, let us write the $m^{\text{th}}$ time cell as
\begin{align*}
J_m=\left[t_m^-,t_m^+\right],
\end{align*}
where \(\mathsf T_{\mathrm{cell}}=|J|/M\) and
\begin{align*}
 t_m^-=t_I^-+m\mathsf T_{\mathrm{cell}},\qquad
  t_m^+=t_I^-+(m+1)\mathsf T_{\mathrm{cell}}.
\end{align*}

Define the time corrector:
\begin{equation}\label{q, def}
\mathfrak{q}(t,x) := -\int_{t_m^-}^{t}
             \left\{F_h(s,x)-\overline F_h^{(m)}(x)\right\}\,\dd s\qquad\text{for } t \in J_m
\end{equation}
and 
\begin{equation}\label{Q, def}
Q(t,x) := \leray \mathfrak{q}(t,x) \qquad\text{for } t \in J_m,
\end{equation}
where $\leray =\mathrm{Id}-\nabla\Delta^{-1}\operatorname{div}$ is the Leray projection on the plane. The two endpoints of \(J_m\)
have \(\mathfrak{q}=Q=0\): for the right endpoint, it follows from~\eqref{Fh-bar, def}. Thus, we may define $\mathfrak{q}$ and $Q$ on the whole time interval $I$ via extension-by-zero outside $J$. On the other hand, let us keep $\mathfrak{q}$ unchanged across the stripes and quarters
inside the cell $J_m$. (Discarding it would create a jump in the time corrector.)

Furthermore, let us set
\begin{equation}\label{Lambda, def}
\Lambda(t) := \Delta^{-1}\operatorname{div}
\left(F_h(t)-\overline F_h^{(m)}\right)\qquad \text{on } \mathring{J}_m.
\end{equation}
We have the explicit expression~$\Lambda=\sum_{k=1}^2\partial_k
((2\pi)^{-1}\log|\cdot|) \star \left[F_h-\overline F_h^{(m)}\right]^k$, which is locally integrable. From
\eqref{a3} and 
\eqref{q, def}, namely that 
\begin{equation*}
\operatorname{div}_x \left(
\sum_{i=1}^4 A_i^{(m)}(x)\,\xi_i\otimes\xi_i\right)
=\overline F_h^{(m)}(x)
=\frac{1}{\mathsf T_{\mathrm{cell}}}\int_{J_m}F_h(t,x)\,\dd t,
\end{equation*} 
we may deduce the evolution equation for $Q$:
\begin{equation}\label{eq:scan-projected-corrector}
 \partial_t Q(t,x)=-F_h(t,x)+\operatorname{div}_x \biga +\nabla\Lambda(t,x).
\end{equation}
As before, we abbreviate
\begin{equation}\label{A}
\biga(x):= \sum_{i=1}^4 A_i^{(m)}(x)\,\xi_i\otimes\xi_i,
\end{equation}
which is a $\sym$-valued function.

\begin{lemma}\label{lem: estimate for Q}
The time corrector $Q$ satisfies the following estimates:
\begin{align}\label{eq:scan-Q-bound}
&\|Q\|_{L^\infty_t L^2_x}+\left\|\operatorname{curl}Q\right\|_{L^\infty_t L^1_x}
    \le \frac{D\bigtime}{\theta h},\nonumber\\
 & \|Q\|_{L^\infty_t L^{{3}/{2}}_x} \le \frac{D\bigtime}{\theta h^{2/3}}.
\end{align}
\end{lemma}

\begin{proof}[Proof of Lemma~\ref{lem: estimate for Q}]

Fix a cell \(J_m=[t_m^-,t_m^+]\) with
\(t_m^+-t_m^-=\mathsf T_{\mathrm{cell}}\).  Recall from~\eqref{Fh, def} that
\[
 F_h(s,x)
 =2\pi\,\frac{d}{ds}\bigl(\eta(s)r(s)\bigr)\,
       \Omega_h(x-\cen(s))\,\xi_i,
 \qquad \text{for } s\in \left[ t^*_{m,i,j},\,t^*_{m,i,j}+d_{ij}^{(m)}\right]\text{ and } x\in\mathbb R^2.
\]
On the reset intervals this expression is zero, because
\(\eta(s)=\eta'(s)=0\). Thus the changes in direction at the
slot joins do not enter the derivative of \(\eta(s)r(s)\).

By the change of variables \(y=x/h\) in \eqref{Omega-h}, for each $p \in ]1,\infty[$ it holds that
\begin{align*}
&\|\Omega_h\|_{L^p_x}
 =h^{2/p-2}\|\Omega\|_{L^p_x},\\
 &\|\nabla_x\Omega_h\|_{L^1_x}
 =h^{-1}\|\nabla\Omega\|_{L^1_x}.
\end{align*} 
Thus,
\[
 \|F_h(s,\cdot)\|_{L^p_x}
 =2\pi\left|\frac{d}{ds}\bigl(\eta(s)r(s)\bigr)\right|
        h^{2/p-2}\|\Omega\|_{L^p_x}
\]
as well as
\[
 \bigl|\operatorname{curl} F_h(s,x)\bigr|
 \le
 2\pi\left|\frac{d}{ds}\bigl(\eta(s)r(s)\bigr)\right|
       \bigl|\nabla\Omega_h(x-\cen(s))\bigr|.
\]
By the triangle inequality applied to the time average $ \overline F_h^{(m)}(x)
 =\frac{1}{\mathsf T_{\mathrm{cell}}}
   \int_{t_m^-}^{t_m^+}F_h(s,x)\,\dd s$, the same bounds hold for \(\overline F_h^{(m)}\) and its curl.

For fixed \(t\in J_m\), we integrate $\mathfrak{q}(s,x)$ in~\eqref{q, def} over \(s\in[t_m^-,t]\) to deduce
\[
\begin{aligned}
 \|Q(t,\cdot)\|_{L^p_x}
 &\le \|\leray\|_{L^p\to L^p}
     \int_{t_m^-}^{t}
     \left(
       \|F_h(s,\cdot)\|_{L^p_x}
       +\left\|\overline F_h^{(m)}\right\|_{L^p_x}
     \right)\,\dd s \\
 &\le
 4\pi\|\leray\|_{L^p\to L^p}
 \bigtime\,h^{\frac{2}{p}-2}\|\Omega\|_{L^p_x}
 \sup_{s\in J_m}
 \left|\frac{d}{ds}\bigl(\eta(s)r(s)\bigr)\right|\quad\text{for }  p>1.
\end{aligned}
\]
Moreover, \(\operatorname{curl}_x Q=\operatorname{curl}_x \mathfrak{q}\) in the sense of distributions, so one may use the explicit expression for $F_h$ to deduce that
\[
\begin{aligned}
 \|\operatorname{curl}_x Q(t,\cdot)\|_{L^1_x}
 &\le \int_{t_m^-}^{t}
     \left(
       \|\operatorname{curl} F_h(s,\cdot)\|_{L^1_x}
       +\left\|\operatorname{curl}\overline F_h^{(m)}(s,\cdot)\right\|_{L^1_x}
     \right)\,\dd s\\
 &\le
 \frac{4\pi\mathsf T_{\mathrm{cell}}}{h}
 \|\nabla\Omega\|_{L^1_x}
 \sup_{s\in J_m}
 \left|\frac{d}{ds}\bigl(\eta(s)r(s)\bigr)\right|.
\end{aligned}
\]

To estimate $\left|\frac{d}{ds}\bigl(\eta(s)r(s)\bigr)\right|
 \le |\eta'(s)|r(s)+|\eta(s)|\,|r'(s)|$, on each stripe we may control $|\eta'|r$ via \eqref{S19}, \(B_\theta\le 1\) and
\(I_{ij}^{(m)}\ge 4\sqrt{62} L \flf\), and control $|\eta||r'|$ via
\(a_i^{(m)}\ge\sqrt{62}\,\flf\) and
\(\nabla_x\log a_i^{(m)}=(\nabla_x a_i^{(m)})/a_i^{(m)}\). Combining these bounds with \eqref{S18} and \eqref{S10}--\eqref{pointwise bound for ai}, we obtain that
\[
\begin{aligned}
 \sup_{s\in I}
 \left|\frac{d}{ds}\bigl(\eta(s)r(s)\bigr)\right|
 &\le
 \frac{(8+\sqrt{2})\bar\eta^2}
      {\theta\sqrt{62}\,\flf}
 \left(
   \frac{\|R_c\|_{L^\infty_{t,x}(I\times\mathbb R^2)}
         +\flf}{2L}
   +\|\nabla_x R_c\|_{L^\infty_{t,x}(I\times\mathbb R^2)}
 \right),
\end{aligned}
\]
where
$\overline\eta:=
\left[\frac4\pi\left((8+\sqrt2)\delta+\frac34 \e_0^2\right)\right]^{1/2}$ as in~\eqref{S18}.

The proof is complete one we choose a large enough parameter $D$ with parameter dependence indicated in~\eqref{D}.     \end{proof}

Now, with the time corrector $Q$ satisfying~\eqref{eq:scan-projected-corrector} at hand, we may further correct $V$, the velocity profile for the principal building block in our convex integration scheme, by $Q$: Define
\begin{equation}\label{V-corr, def}
\vcorr:=V+Q.
\end{equation}
Recall from~\eqref{S3}, \eqref{Pi-r, def}, \eqref{Wr, def} that the rescaled Lamb--Chaplygin dipole $U_r$ can be decomposed as  
\begin{equation*}
U_r=W_r+\nabla\Pi_r,
\end{equation*} 
where \(W_{r(t)}\) is the profile in direction \(\xi_i\) on
\(\timeint\) with $\leray (W_r)=U_r$ and  $\nabla\Pi_r(x)\in L^2(\mathbb R^2)$. Then, for each \(t\in I\), it holds that
\begin{equation}\label{Vcorr, eq}
\vcorr (t,\cdot)
 =\leray\!\Big[
      \eta(t)W_{r(t)}(\cdot-\cen(t))+\mathfrak{q}(t,\cdot)
   \Big],
\end{equation}
and hence (suppressing the variables $(t,x)$)
\begin{equation*}
 \operatorname{curl} \,\left(\vcorr\right)
 =\operatorname{curl}\,\Big[
      \eta W_r(\cdot-\cen)+\mathfrak{q}\Big].
\end{equation*}

\begin{lemma}\label{lem: V-corr}
Let $\vcorr$ be defined as in~\eqref{V-corr, def}. 
It has zero circulation, zero (vector) first vorticity moment, and zero (scalar) second vorticity moment. More precisely:
\begin{align*}
&\int_{\mathbb R^2}\operatorname{curl}\,\left(\vcorr\right)(t,x)\,\dd x =0,\\
&\int_{\mathbb R^2}x\,\operatorname{curl}\,\left(\vcorr\right)(t,x)\,\dd x = \mathbf{0},\\
&\int_{\mathbb R^2}|x|^2\,\operatorname{curl}\left(\vcorr\right)(t,x)\,\dd x =0.
\end{align*}

\end{lemma}

\begin{proof}[Proof of Lemma~\ref{lem: V-corr}]

First, observe that $\curl\vcorr(t,\cdot)$ is compactly supported in \(\overline{\B}_{2\sqrt{2}L+h}(0)\) for all $t\in I$: This is because the support of $ \eta(t)W_{r(t)}(\cdot-\cen(t))$ is within distance \(2r(t)<h\) of the centre, and $\mathfrak{q}$ is the time integral of source terms supported within distance $h$ of the stripes.   Then, an integration by parts over sufficiently large balls in $\R^2$ leads to the first identity.

Now we turn to the first vorticity moment of $\vcorr$. Integration by parts yields that
\begin{align*}
\int_{\mathbb R^2}x\,\operatorname{curl}(V+Q)(t,x)\,\dd x =\left(\int_{\mathbb R^2} \Big\{\eta(t)W_{r(t)}(x-\cen(t))+\mathfrak{q}(t,x)\Big\}\,\dd x  \right)^\perp.
\end{align*}
To see that the right-hand side vanishes, 
Recall from~\eqref{S7} that
\[
 \int_{\mathbb R^2}
      \eta(t)W_{r(t)}\big(x-\cen(t)\big)\,\dd x
 =2\pi\eta(t)r(t)\xi_i
\]
for each $t\in\timeint$. Also, the spatial integral of \(F_h(s,x)\) is
\(2\pi\frac{d}{ds}(\eta(s)r(s))\,\xi_k\), and that of \(\overline F_h^{(m)}\) is zero. The latter follows from the fact that the integrand is the divergence of a compactly supported tensor; see~\eqref{a3}. Hence, we have
\[
\begin{aligned}
 \int_{\mathbb R^2}\mathfrak{q}(t,x)\,\dd x
 &=-2\pi\sum_{i=1}^4\sum_{j=1}^{N_i}
   \xi_i\int_{\timeint\cap[t_m^-,t]}
          \frac{d}{ds}\bigl(\eta(s)r(s)\bigr)\,\dd s\\
 &=-2\pi\eta(t)r(t)\xi_i,
\end{aligned}
\]
from which it follows that
\begin{equation*}
 \int_{\mathbb R^2}
 \Big\{\eta(t)W_{r(t)}(x-\cen(t))+\mathfrak{q}(t,x)\Big\}\,\dd x=0.
\end{equation*}

Finally, for the second moment, we integrate by parts to deduce that 
\begin{align*}
\int_{\mathbb R^2}|x|^2\,\operatorname{curl}\left(\vcorr\right)(t,x)\,\dd x =-2\int_{\mathbb R^2}
  x\wedge\left[\eta(t)W_{r(t)}\big(x-\cen(t)\big)+\mathfrak{q}(t,x)\right]\,\dd x.
\end{align*}
It thus remains to check that torque of the term in~\eqref{Vcorr, eq} whose Leray projection is $\vcorr$, namely  $\eta(t)W_{r(t)}(\cdot-\cen(t))+\mathfrak{q}(t,\cdot)$, is zero.

To this end, notice that as \(W_{r(t)}\) is even, 
\[
\begin{aligned}
 &\int_{\mathbb R^2}
   \eta(t)W_{r(t)}(x-\cen(t))\otimes x\,\dd x\\
 &\qquad=\eta(t)\int_{\mathbb R^2}
   W_{r(t)}(y)\otimes \big(y+\cen(t)\big)\,\dd y
 =2\pi\eta(t)r(t)\,\xi_i\otimes\cen(t).
\end{aligned}
\]
Similarly, as \(\Omega_h\) is radial, we have
\(\int_{\mathbb R^2}y\,\Omega_h(y)\,\dd y=0\). Therefore, 
\begin{align*}
\int_{\mathbb R^2} F_h(s,x)\otimes x\,\dd x
 =2\pi\,\frac{d}{ds}\bigl(\eta(s)r(s)\bigr)\,
       \xi_k\otimes\cen(s)\qquad\text{ for all } s \in \timeint.
\end{align*}
Moreover, labelling $\biga \equiv \sum_{i=1}^4 A_i^{(m)}\,\xi_i\otimes\xi_i$ as in~\eqref{A}, we have 
\begin{align*}
\int_{\mathbb R^2}
      \left(\overline F_h^{(m)}(x)\right)^\alpha x^\beta\,\dd x
 &=\sum_{\gamma=1}^2
   \int_{\mathbb R^2}
      \partial_{x_\gamma} A^{(m)}_{\alpha\gamma}(x)\,
      x^\beta\,\dd x
 =-\int_{\mathbb R^2} A^{(m)}_{\alpha\beta}(x)\,\dd x.
\end{align*} 

Inserting the above identities into the defining expression~\eqref{q, def} for $\mathfrak{q}$ and integrating by parts in the time variable, one deduces that
\begin{equation}
\label{eq:scan-compact-moment}
\begin{aligned}
 &\int_{\mathbb R^2}
   \Big[\eta(t)W_{r(t)}(x-\cen(t))+\mathfrak{q}(t,x)\Big]
   \otimes x\,\dd x\\
 &\quad=
 2\pi\sum_{i=1}^4\sum_{j=1}^{N_i}
   \int_{\timeint\cap[t_m^-,t]}
      \eta(s)r(s)\,\xi_i\otimes\p_s{\cen}(s)\,\dd s
 -(t-t_m^-)\int_{\mathbb R^2}\biga(x)\,\dd x\\
 &\quad=
 2\pi\sum_{i=1}^4\sum_{j=1}^{N_i}
   \int_{\timeint\cap[t_m^-,t]}
      \eta^2(s)\,\xi_i\otimes\xi_i\,\dd s
 -(t-t_m^-)\int_{\mathbb R^2}\biga(x)\,\dd x.
\end{aligned}
\end{equation}
On resets both integrands
in~\eqref{eq:scan-compact-moment} are zero because \(\eta=0\). The last line follows since
\begin{align*}
\p_s{\cen}(s)=\frac{\eta(s)}{r(s)}\xi_k\qquad\text{for } s \in \timeint.
\end{align*} 
This computation shows that $\int_{\mathbb R^2}
   \Big[\eta(t)W_{r(t)}(x-\cen(t))+ \mathfrak{q}(t,x)\Big]
   \otimes x\,\dd x$ is a symmetric matrix field; nevertheless, the quantity in consideration, \emph{i.e.}, $$\int_{\mathbb R^2}
 x\wedge\left[\eta(t)W_{r(t)}(x-\mathbf c(t))+ \mathfrak{q}(t,x)\right]\,\dd x,$$ is the anti-symmetric part of that matrix field. Hence it must be zero.

The proof is now complete.  \end{proof}

\subsection{The ``raw'' updated Euler--Reynolds triple}
\label{sec:parameter-order}

In this subsection, we shall ``almost'' prove the Key  Lemma~\ref{lem: key}: we construct a triple $(u',p',R')$ from the given $(u,p,R)$ that solves the Euler--Reynolds equation in the sense of distribution, such that the stress $R'$ and the velocity increments $u'-u$ satisfy the desired smallness estimates in Lemma~\ref{lem: key}. However, 
$(u',p',R') \notin \CC^\infty$ in general, so we need a further smoothing step in the next subsection to conclude the proof.

\subsubsection{Preparation}

We first introduce some notation. Set
\begin{equation}\label{s, r-star}
\begin{cases}
s_* := \max_{1\leq i\leq 4} s_i,\\
r_* :=\sup_{t\in I}r(t).
\end{cases}
\end{equation}
Also write
\begin{equation}\label{M3, def}
M_3 := \sup_{t\in I}\|u(t,\cdot)\|_{L^3_x},
\end{equation}
which is assumed to be finite by the hypothesis of Lemma~\ref{lem: key}. Recall from~\eqref{t**, def} that $t^{**}_m$ is the midpoint of the $m^{\text{th}}$ time cell $J_m=\left[t_m^-,t_m^+\right]=\left[t_I^-+m\mathsf T_{\mathrm{cell}}, t_I^-+(m+1)\mathsf T_{\mathrm{cell}}\right]$. Also, define
\begin{equation}\label{def, b0m}
b_0^{(m)}(x)=b(t_m^{**},x)-16\flf,\qquad\text{for any } x \in \R^2.
\end{equation}
Note that by \eqref{S10} and \eqref{def, Rc}, $b\equiv 16\flf$ outside the compact support of $R_c$. Hence, $b_0^{(m)}$ is compactly supported. As will become clear later, we introduce $b_0^{(m)}$ since the $L^1_x$-norm of the tensor field $R+A^{(m)}-b_0^{(m)}{\bf I}_2$ is small.

Let us also recall the following quantities:
\begin{itemize}
\item
In~\eqref{Lambda, def} we have defined 
\begin{equation*}
\Lambda(t,\cdot) := \Delta^{-1}\operatorname{div}
\left(F_h(t,\cdot)-\overline F_h^{(m)}(\cdot)\right)\qquad \qquad\text{for } t \in \mathring{J}_m.
\end{equation*}

\item
$V(t,x)=\eta(t)U_{r(t)}(x-\cen(t))$ is the principal building block; see~\eqref{def, V, travelling dipole}.

\item
$Q$ is the solenoidal time corrector defined in~\eqref{Q, def}:
\begin{equation*}
Q(t,x):= -\leray\left(\int_{t_m^-}^{t}
             \left\{F_h(s,x)-\overline F_h^{(m)}(x)\right\}\,\dd s\right)\qquad\text{for } t \in J_m.
\end{equation*}

\item
$\biga = \sum_{i=1}^4 A_i^{(m)}\,\xi_i\otimes\xi_i$ as in~\eqref{A}. An explicit expression for $A_i^{(m)}=A_i^{(m)}(x)$ is given in~\eqref{a4}. It serves as a good approximation of $a_i$ in the spatial square $Q_i$ and time cell $J_m$; see Lemma~\ref{lem: estimates given eta, r}, $(v)$. Here, $a_i$ are the profiles of $R_c$, the compactly cutoff version of the Reynolds stress tensor $R$, via~\eqref{S11}:
$$-R_c(t,x) = \sum_{i=1}^4\Big(a_i(t,x)\xi_i\otimes\xi_i\Big) - b(t,x){\bf I}_2.$$

\item
$S_c = {\rm div}^{-1}(F_c-F_h)$ as in~\eqref{Sc}, whose $L^1_x$-norm is bounded by $D\sqrt{\theta}$ in Lemma~\ref{lem: estimates given eta, r}, $(iv)$. 
\end{itemize}

\subsubsection{Construction of the ``raw'' triple}

Now we construct the ``raw'' triple $(u',p',R')$ and prove the smallness estimates.

\begin{lemma}\label{lemma: update rule for u', p', R'}
Given any Euler--Reynolds triple $(u,p,R)$ as in Lemma~\ref{lem: key}. Define $(u',p',R'): J_m \times \R^2 \longrightarrow \R^2 \times \R \times \sym$ in each of the time cell $J_m \subset J$ as follows:
\begin{equation}\label{u',p',R', def}
\begin{cases}
u' := u + \vcorr = u+V+Q,\\
p' := p+p_V-\Lambda-b_0^{(m)},\\
R' := R+\biga-b_0^{(m)}\mathbf{I}_2+S_c \\ 
\qquad\qquad+u\otimes V+V\otimes u +Q\otimes(u+V)+(u+V)\otimes Q+Q\otimes Q,
\end{cases}
\end{equation}
where $p_V$ is as in the Euler equation~\eqref{Euler for V with Fc} for $V$. Set $(u',p',R')\equiv (u,p,R)$ on $(I \setminus J)\times \R^2$. Then ${\rm div}\,u'=0$ and $(u',p',R')$ is a Euler--Reynolds triple; \emph{i.e.}, the equation
\begin{align}\label{ER'}
\p_t u' + {\rm div} \left(u'\otimes u'\right) + \na p' = {\rm div}\,R'
\end{align} 
holds in the sense of distributions in $I \times \R^2$.
\end{lemma}

\begin{proof}[Proof of Lemma~\ref{lemma: update rule for u', p', R'}]

The incompressibility of $u'$ follows from that of $u$, $V$, and $Q$. Indeed, by~\eqref{def, V, travelling dipole}, $V$ is $\eta(t)$ times the rescaled Lamb--Chaplygin dipole, which is an exact solution to the Euler equation; also, by~\eqref{Q, def}, $Q$ lies in the image of Leray projection by construction.  

To see that $(u',p',R')$ is an Euler--Reynolds triple, we add up together the equations for $u$, $V$, and $Q$. Indeed, $u$ satisfies the Euler--Reynolds equation $$\p_t u + {\rm div}(u\otimes u) + \na p = {\rm div}\,R,$$ $V$ satisfies~\eqref{Euler for V with Fc}:
\begin{align*}
\p_tV + {\rm div}(V \otimes V) + \na p_V = F_c,
\end{align*}
and $Q$ satisfies~\eqref{eq:scan-projected-corrector}, namely
\begin{equation*}
 \partial_t Q(t,x)=-F_h(t,x)+\operatorname{div} \biga +\nabla\Lambda(t,x).
\end{equation*}
We then conclude by \(\operatorname{div}S_c=F_c-F_h\) from \eqref{Sc}.   \end{proof}

 In passing, recall from~\eqref{Euler for V with Fc} that the explicit formula for $p_V$ is as follows:
\begin{align*}
p_V(t,x):=\eta^2(t) P_r\big(x-\cen(t)\big) -\eta'(t)\Pi_r\big(x-\cen(t)\big) -\eta(t) r'(t)\,(\partial_r\Pi_r)\big(x-\cen(t)\big).
\end{align*}

\begin{proposition}\label{propn: raw estimate}
Under the hypotheses of Lemma~\ref{lem: key}, with the parameters chosen in \S\ref{subsec: parameter D,h}, the Euler--Reynolds triple $(u',p',R')$ satisfies
\begin{align}\label{eq:scan-raw-increment} 
&\operatorname*{ess\,sup}_{t\in I}
     \|R'(t,\cdot)\|_{L^1_x}
< \frac{\varepsilon_\star}{2},\nonumber\\
&\sup_{t\in I}\|(u'-u)(t,\cdot)\|_{L^2_x}
 +\sup_{t\in I}\|\operatorname{curl}(u'-u)(t,\cdot)\|_{L^1_x}
\le C_0\sqrt{\delta}+\frac{\varepsilon_\star}{8},
\end{align}
where $C_0$ is a universal constant.
\end{proposition}

Here, the key is to choose the large number $D$ as in~\eqref{D}:
\begin{align*}
D = D\bigg(\e_0,\,\delta,\,
\|R_c\|_{L^\infty_{t,x}},\,
\|\partial_tR_c\|_{L^\infty_tL^1_x},\,
\|\nabla_x R_c\|_{L^\infty_{t,x}},\,
\|\nabla_xR_c\|_{L^\infty_tL^1_x},\, \underbrace{\|u\|_{L^\infty_t L^3_x}}_{=:\, M_3},\, I\bigg).
\end{align*}
Once we specify $D$, the choice of $h$, $s_i$, $\bigtime$, $\flf$, as well as $L$ and the cutoff Reynolds stress tensor $R_c$, are all fixed. To this end, we recall~\eqref{L, choice}, \eqref{e0}, \eqref{floor}, \eqref{h, choice}, \eqref{si, choice}, and \eqref{Tcell}, and summarise them below for ease of future reference:
\begin{align}\label{parameter, summary}
&L=L\left(\e_0,\delta,\|R\|_{\CC^0_t L^1_x}, R(I)\right)>1000,\qquad  \flf := \frac{\e_0^2}{1024 L^2},\qquad \bigtime < \frac{\e_\star}{100 D},\nonumber\\
&s_i<\min\left\{h,\frac{\e_0^2}{2D},\frac{\e_\star h}{100D}\right\},\qquad h := \min\left\{ \frac{1}{128},\,\left(\frac{\e_0}{400 D}\right)^2 \right\},\nonumber\\
&\e_0 =  \min\bigg\{
\frac12,\ \frac{\e_\star}{200},
\frac{\e_\star}{16\sqrt{3/\pi}\, (\|U\|_{L^2_x}+\|\operatorname{curl}(U)\|_{L^1_x})}
\bigg\}.
\end{align}

\begin{proof}[Proof of Proposition~\ref{propn: raw estimate}]

We first prove the $L^\infty_t L^1_x$-bound for $R'$, while obtaining a specification of $D$ along the way. Then we deduce the desired bounds for $(u'-u)$.

\smallskip
\noindent
\underline{\bf Estimate for $R'$.} To bound $\operatorname*{ess\,sup}_{t\in I}
     \|R'(t,\cdot)\|_{L^1_x}$, recall from Lemma~\ref{lemma: update rule for u', p', R'} that the stress tensor $R'$ is defined as follows:
     \begin{equation}\label{R', final}
     R' := R+\biga-b_0^{(m)}\mathbf{I}_2+S_c +u\otimes V+V\otimes u +Q\otimes(u+V)+(u+V)\otimes Q+Q\otimes Q.
     \end{equation}
We shall treat the nine terms on the right-hand side in four groups.

\smallskip
\noindent
\underline{\textbf{Terms 1--3.}} This is the most technical part of this proposition. Indeed, notice that if we bound these terms individually, at the best we have 
\begin{align*}
 \|R(t,\cdot)\|_{L^1_x}\le\delta,\qquad
 \left\|\biga\right\|_{L^1_x}\le\sum_i\Delta_i^{(m)},\qquad
 \left\|b_0^{(m)}{\bf I}_2\right\|_{L^1_x}\le 16\sqrt2\,\delta.
\end{align*}
The second inequality follows from~\eqref{first, difference est}, where $\Delta_i^{(m)}=s_i\sum_j I_{ij}^{(m)}
=\Delta_i(t_m^{**})$ approximates the area integral of $a_i$ over the square $Q_i$ in the $m^{\text{th}}$ time cell. The third  holds by~\eqref{S9} and $0\le \sqrt{\left|R_c^{(m)}\right|^2+\flf^2}-\flf \le \left|R_c^{(m)}\right|$. As usual, $|\cdot|$ denotes the Hilbert--Schmidt norm of a matrix, so $|{\bf I}_2|=\sqrt{2}$. To obtain smallness, we shall need to group these terms together and exploit cancellation.

To this end, write
\begin{align}\label{split, group 1}
 R(t)+\biga-b_0^{(m)}{\bf I}_2
  = \big[R(t)-R_c(t)\big]+\big[R_c(t)-R_c^{(m)}\big]
       +\big[R_c^{(m)}+\biga-b_0^{(m)}{\bf I}_2\big].
\end{align}
We emphasise that the only terms dependent on the time variable here are $R$ and $R_c$.

\begin{itemize}
\item
By~\eqref{S9},  
\begin{equation}\label{group1, a}
\left\|R(t,\cdot)-R_c(t,\cdot)\right\|_{L^1_x} \leq \e_0.
\end{equation}
\item
Since $R_c$ is smooth, spatially compactly supported with $\partial_tR_c\in L^\infty_tL^1_x$, we deduce from the choice of $t^{**}_m$ in~\eqref{t**, def} that 
\begin{align}\label{group 1,b}
 \left\|R_c(t,\cdot)-R_c^{(m)}(\cdot)\right\|_{L^1_x}
 \le \left|t-t_m^{**}\right|\sup_{s\in I}\left\|\partial_t R_c(t,\cdot)\right\|_{L^1_x}
 \le \tfrac12 D\,\bigtime.
\end{align}

\item
To bound $R_c^{(m)}+\biga-b_0^{(m)}{\bf I}_2$, we estimate the difference between $\sum_{i=1}^4g_i^{(m)}\xi_i\otimes\xi_i$ and $\biga$
using~\eqref{S25} in Lemma~\ref{lem: estimates given eta, r}, where $g_i^{m} = a_i^{(m)} \1_{Q_i}$. Notice that \begin{align*}
 R_c^{(m)}+\sum_i g_i^{(m)}\xi_i\otimes\xi_i-b_0^{(m)}{\bf I}_2
       =8\flf\sum_i\1_{Q_i}\xi_i\otimes\xi_i.
\end{align*}

 To track the constant \(D\), we first give
bounds for each quantity on the right-hand side of \eqref{S25}, reproduced below:
\begin{align}\label{S25'}
 \left\|A_i^{(m)}-g_i^{(m)}\right\|_{L^1_x}
 &\le 8\theta\,\Delta_i^{(m)}
   +s_i\left\| \xi_i^\perp \cdot \nabla a_i^{(m)}\right\|_{L^1(Q_i)}
   +\frac{s_i}{2h}\left\|\nabla\Omega\right\|_{L^1_x}\left\|g_i^{(m)}\right\|_{L^1_x} \nonumber \\
 &\quad+ h\left(\int |y|\Omega(y)\,\dd y\right)
 \left|D_x g_i^{(m)}\right|(\mathbb R^2)
\end{align}
To further unburden the notation, we introduce the symbols:
\begin{equation*}
\begin{cases}
m_0 := (8+\sqrt{2})\delta+\frac{\varepsilon_0^2}{4},\\
G_1 :=(8+\sqrt{2})\left\|\nabla_xR_c\right\|_{L^\infty_t L^1_x(I \times \R^2)}. 
\end{cases}
\end{equation*}
 
We control the four terms on the right-hand side of~\eqref{S25'} separately.

\smallskip
\noindent
\underline{The third term.} Recall that $a_i$ are positive, so by~\eqref{S12} one has
\begin{align}\label{gim, L1 bound}
 \left\|g_i^{(m)}\right\|_{L^1_x}
 =\int_{Q_i}a_i^{(m)}(x)\,\mathrm{d}x
 \le m_0.
\end{align}

\smallskip
\noindent
\underline{The second term.}
By~\eqref{S10}, the common nonlinear part  \(a_i\) is
\(8\sqrt{|R_c|^2+\flf^2}\), which satisfies 
\begin{align*}
 \left|
 \partial_{x_k}\!\left(8\sqrt{|R_c|^2+\flf^2}\right)
 \right|
 =
 \frac{8|R_c:\partial_{x_k}R_c|}
      {\sqrt{|R_c|^2+\flf^2}}
 \le 8\left|\partial_{x_k}R_c\right|\qquad\text{for } k \in\{1,2\}.
\end{align*} 
The  linear parts of \(a_1,a_2,a_3,a_4\), evaluated on a symmetric
matrix \(H\), are respectively
\(-H_{11}+H_{12}\), \(-H_{22}+H_{12}\), \(-2H_{12}\), and \(0\). Their operator norms are at most \(\sqrt{2}\) for the Hilbert--Schmidt norm. We thus have the pointwise derivative bound:
\begin{equation*}
\left|\nabla_x a_i\right|\le(8+\sqrt{2})\left|\nabla_x R_c\right|.
\end{equation*}
Evaluating at \(t_m^{**}\) and using \(|\xi_i^\perp|=1\), we obtain that
\begin{align}\label{nabla ai, bound}
 \left\|\xi_i^\perp\cdot\nabla_x a_i^{(m)}\right\|_{L^1(Q_i)}
 \le \left\|\nabla_x a_i^{(m)}\right\|_{L^1(Q_i)}
 \le G_1.
\end{align}

\smallskip
\noindent
\underline{The first term.} 
We next bound the ``quadrature mass'' \(\Delta_i^{(m)}\). Recall that
\[
 \Delta_i^{(m)}=s_i\sum_{j=1}^{N_i}I_{ij}^{(m)},
 \qquad
 I_{ij}^{(m)}
   =\int_0^{4L}a_i^{(m)}(Y_i(\lambda,z_j))\,\mathrm{d}\lambda,
\]
where \(z_j=(j-\tfrac12)s_i\) is the midpoint of \(\band_{ij}=[(j-1)s_i,js_i]\). Since each band has length
\(s_i\), we can write
\[
\begin{aligned}
 &\Delta_i^{(m)}-\int_{Q_i}a_i^{(m)}(x)\,\mathrm{d}x\\
 &\quad=\sum_{j=1}^{N_i}\int_0^{4L}\int_{\band_{ij}}
 \left[
 a_i^{(m)}(Y_i(\lambda,z_j))-a_i^{(m)}(Y_i(\lambda,z))
 \right]\,\mathrm{d}z\,\mathrm{d}\lambda .
\end{aligned}
\]
Thus, for each fixed \(\lambda\in[0,4L]\) and \(z\in \band_{ij}\), the
fundamental theorem of calculus and
\(\partial_zY_i=\xi_i^\perp\) lead to 
\begin{align*}
\left|a_i^{(m)}(Y_i(\lambda,z_j))
           -a_i^{(m)}(Y_i(\lambda,z))\right| &\le
 \int_{\min\{z,z_j\}}^{\max\{z,z_j\}}
 \bigl|(\xi_i^\perp\cdot\nabla_x a_i^{(m)})
                      (Y_i(\lambda,\zeta))\bigr|\,\mathrm{d}\zeta\\
                      &\leq \int_{\band_{ij}}
 \bigl|(\xi_i^\perp\cdot\nabla_x a_i^{(m)})
                      (Y_i(\lambda,\zeta))\bigr|\,\mathrm{d}\zeta .
\end{align*} 
The last bound is independent of \(z\). A further integration in \(z\)
over \(\band_{ij}\) gives us a good factor \(s_i\). Then, integrating in \(\lambda\)
and summing over the bands, we obtain that
\[
\begin{aligned}
 \left|\Delta_i^{(m)}-\|g_i^{(m)}\|_{L^1_x}\right|
 &\le s_i\sum_{j=1}^{N_i}\int_0^{4L}\int_{\band_{ij}}
   \bigl|(\xi_i^\perp\cdot\nabla_x a_i^{(m)})
                        (Y_i(\lambda,\zeta))\bigr|
                    \,\mathrm{d}\zeta\,\mathrm{d}\lambda\\
 &=s_i\left\|\xi_i^\perp\cdot\nabla_x a_i^{(m)}\right\|_{L^1(Q_i)}\\
 &  \le s_i G_1 .
\end{aligned}
\]
By choosing $D \geq G_1$ and recalling \(s_i<\varepsilon_0^2/(2D)\) from~\eqref{parameter, summary}, we conclude
that
\begin{equation}\label{Delta im bound}
  \Delta_i^{(m)}
 \le m_0+s_iG_1
 \le m_0+\frac{\varepsilon_0^2}{2}.
\end{equation}

\smallskip
\noindent
\underline{The second term.} For the last term of \eqref{S25'}, we need to control the total variation
of \(D_x g_i^{(m)}\) on \(\mathbb R^2\), including the jump at
\(\partial Q_i\). Denote by \(\nu_{Q_i}\) the outward unit normal to \(Q_i\). For every
\(\varphi\in C^\infty_c(\mathbb R^2;\mathbb R^2)\), integration
by parts on \(Q_i\) yields that
\[
\begin{aligned}
 \bra D_x g_i^{(m)},\varphi \ket
 &=-\int_{Q_i}a_i^{(m)}(x)
                   \operatorname{div}_x\varphi(x)\,\mathrm{d}x\\
 &=\int_{Q_i}\nabla_x a_i^{(m)}(x)\cdot\varphi(x)\,\mathrm{d}x -\int_{\partial Q_i}a_i^{(m)}(x)
              \varphi(x)\cdot\nu_{Q_i}(x)\,\mathrm{d}\mathcal H^1(x).
\end{aligned}
\]
As ${\rm supp}(R_c) \subset \B_L(0) \Subset Q_i$, we have \(a_i^{(m)} \equiv 8\flf\)
near \(\partial Q_i\); see~\eqref{ai outside supp Rc}. Hence,
\begin{equation}\label{Dxgi}
 D_x g_i^{(m)}
 =\mathbf{1}_{Q_i}\nabla_x a_i^{(m)}\,\mathrm{d}x
  -8\flf\,\nu_{Q_i}\,
                 \mathcal H^1\!\mres{\partial Q_i}.
\end{equation}
As \(\mathcal H^1(\partial Q_i)=16L\), we deduce that
\begin{align}\label{Dxgi, bound}
\left|D_x g_i^{(m)}\right|(\mathbb R^2) &
 =\int_{Q_i}\left|\nabla_x a_i^{(m)}(x)\right|\,\mathrm{d}x +128\flf L\nonumber\\
&\le G_1+128\flf L.
\end{align}

\smallskip
Substituting the bounds~\eqref{gim, L1 bound}, \eqref{nabla ai, bound}, \eqref{Delta im bound}, and \eqref{Dxgi, bound}  into \eqref{S25}, we deduce that
\[
\begin{aligned}
 \sum_{i=1}^4\left\|A_i^{(m)}-g_i^{(m)}\right\|_{L^1_x}
 &\le 8\theta\sum_{i=1}^4\Delta_i^{(m)}
       +G_1\sum_{i=1}^4s_i
       +\frac{m_0\|\nabla\Omega\|_{L^1_x}}{2h}
                            \sum_{i=1}^4s_i\\
 &\quad+4h(G_1+128\flf L)
                     \int_{\mathbb R^2}|y|\Omega(y)\,\mathrm{d}y\\
 &\le 32\left(m_0+\frac{\varepsilon_0^2}{2}\right)\theta
       +4\left[
           G_1+(G_1+128\underline f L)
                    \int_{\mathbb R^2}|y|\Omega(y)\,\mathrm{d}y
         \right]h\\
 &\quad+2m_0\|\nabla\Omega\|_{L^1_x}\frac{s_*}{h}.
\end{aligned}
\]
In the last line we used $s_* < h$; see \eqref{s, r-star}.  Thus, when choosing \(D\) in \eqref{D}, it suffices for this
estimate to require
\[
 D\ \ge\
 \max\left\{
  32 \left(m_0+\frac{\varepsilon_0^2}{2}\right),\,
  4\left[
     G_1+(G_1+128\flf L)
             \int_{\mathbb R^2}|y|\Omega(y)\,\mathrm{d}y
    \right],\,
  2m_0\|\nabla\Omega\|_{L^1_x}
 \right\}.
\]
This also ensures \(D\ge G_1\) as used above. The quantities in
this maximum are fixed before \(h\) and \(s_i\) are chosen.  We conclude that uniformly in cell index \(m\),
\begin{equation}\label{group 1, c}
 \left\|\biga
       -\sum_{i=1}^4g_i^{(m)}\xi_i\otimes\xi_i\right\|_{L^1_x}
 \le D\left(\theta+h+\frac{s_*}{h}\right).
\end{equation} 
\end{itemize}

Putting together \eqref{group1, a}, \eqref{group 1,b}, and \eqref{group 1, c}, and recalling $\theta = \sqrt{h}$ from \eqref{theta, def}, we arrive at 
\begin{equation}\label{Terms 1-3}
 \left\|R(t,\cdot)+\biga-b_0^{(m)}{\bf I}_2\right\|_{L^1_x}
 \le \varepsilon_0+\frac{ D\bigtime}{2}
       +D\left(\sqrt{h}+h+\frac{s_*}{h}\right)+\frac{\varepsilon_0^2}{2}.
\end{equation}

\smallskip
\noindent
\underline{\textbf{Term~4.}} This term has already been estimated in 
Lemma~\ref{lem: estimates given eta, r}, $(iv)$: we have
\begin{equation}\label{Term 4}
\|S_c\|_{L^\infty_t L^1_x} \leq D\sqrt{h}.
\end{equation}

\smallskip
\noindent
\underline{\textbf{Terms~5 $\&$ 6.}} We now control $u \otimes V + V \otimes u$. By definition of $V$ in~\eqref{def, V, travelling dipole} and the basic estimates~\eqref{S3} and \eqref{S4} for the rescaled Lamb--Chaplygin dipoles $U_{r(t)}$, we have
\begin{align*}
& \|V(t,\cdot)\|_{L^{3/2}_x}=|\eta(t)|\,r(t)^{1/3}\|U\|_{L^{3/2}_x},\\
&
 \|V(t)\|_{L^2_x}=|\eta(t)|\,\|U\|_{L^2_x}.
\end{align*}
In particular, as $|U(x)| \lesssim |x|^{-2}$ near spatial infinity, we have $U \in L^{3/2}(\R^2) \cap L^2(\R^2)$. Thus, by H\"{o}lder's inequality we deduce that
\begin{equation}\label{Terms 5,6}
\left\|\big(u \otimes V + V \otimes u\big)(t,\cdot)\right\|_{L^1_x} \leq 2M_3\bar\eta\,r_*^{1/3}\|U\|_{L^{3/2}_x}.
\end{equation}
Here $r_* = \sup_{I} r$ as in~\eqref{s, r-star}, $M_3 = \|u\|_{L^\infty_t L^3_x}$ as in~\eqref{M3, def}, and $ \overline\eta=
\left[\frac4\pi\left((8+\sqrt2)\delta+\frac34 \e_0^2\right)\right]^{1/2}$ as in~\eqref{S18}, which is an upper bound for $\eta$.

\smallskip
\noindent
\underline{\textbf{Terms~7--9.}} By H\"{o}lder's inequality, the bound~\eqref{eq:scan-Q-bound} for $Q$, and  \(\theta=\sqrt h\) as in \eqref{theta, def}, we have
\begin{align}\label{Terms 7-9}
&\left\|\Big(Q\otimes (u+V) + (u+V)\otimes Q + Q \otimes Q\Big)(t,\cdot)\right\|_{L^1_x} \nonumber\\
&\qquad \leq 2\|Q(t,\cdot)\|_{L^{3/2}_x}\|u(t,\cdot)\|_{L^3_x}+ 2 \|Q(t,\cdot)\|_{L^2_x}\|V(t,\cdot)\|_{L^2_x} + \|Q(t,\cdot)\|^2_{L^2_x}\nonumber\\
&\qquad \le 2M_3\frac{D\mathsf T_{\rm cell}}{\theta h^{2/3}}
 +2\bar\eta\|U\|_2\frac{D\mathsf T_{\rm cell}}{\theta h}
 +\left(\frac{D\mathsf T_{\rm cell}}{\theta h}\right)^2\nonumber\\
&\qquad = 2M_3\frac{D\mathsf T_{\rm cell}}{h^{7/6}}
       +2\bar\eta\|U\|_2\frac{D\mathsf T_{\rm cell}}{h^{3/2}}
       +\frac{D^2\mathsf T_{\rm cell}^2}{h^3}.
\end{align}

\smallskip
Therefore, in view of the decomposition for $R'$ in~\eqref{R', final}, as well as the bounds in~\eqref{Terms 1-3}, \eqref{Term 4}, \eqref{Terms 5,6}, and \eqref{Terms 7-9}, we conclude that
\begin{align}\label{R' estimate, pre}
\|R'\|_{L^\infty_tL^1_x} &\leq \varepsilon_0+\frac{ D\bigtime}{2}
       +D\left(\sqrt{h}+h+\frac{s_*}{h}\right)+\frac{\varepsilon_0^2}{2} + D\sqrt{h}+2M_3\bar\eta\,r_*^{1/3}\|U \|_{L^{3/2}_x} \nonumber\\
       &\qquad  + 2M_3 \frac{D\bigtime}{h^{7/6}} + 2\bar\eta\|U\|_{L^2_x}\,\frac{D\mathsf T_{\mathrm{cell}}}{h^{3/2}} + \frac{D^2 \mathsf{T}_{\mathrm{cell}}^2}{h^3}. 
\end{align}

\smallskip
\noindent
\underline{\bf Choice of $D$ and $\bigtime$.}  Along with the previous arguments, requirements have been imposed on $D$ according to~\eqref{D}, \emph{i.e.}, depending only on 
$M_3$, $\|R_c\|_{L^\infty_{t,x}}$, $\|\partial_tR_c\|_{L^\infty_tL^1_x}$, $\|\nabla_x R_c\|_{L^\infty_{t,x}}$, $\|\nabla_xR_c\|_{L^\infty_tL^1_x}$, as well as $\delta$ and $\e_0$. Also recall from~\eqref{parameter, summary} that $\e_0$ depends only $\e_\star$ and the norms of $U$, $\curl\, U$, which are universal constants since $U$ is a fixed vector field.

We can make an explicit choice of $D$, although completely unnecessary: Set 
$m_0=(8+\sqrt2)\delta+\varepsilon_0^2/4$ and hence $\bar\eta=\sqrt{4(m_0+\varepsilon_0^2/2)/\pi}$ as defined in~\eqref{S18}. Then put 
\begin{align*}
F_*=\left(\|W_1\|_{L^1_x}
       +\bigl\||x|\,|\nabla_x W_1|\bigr\|_{L^1_x}+2\pi\right)\bar\eta^2
 \left(\frac{(8+\sqrt2)\bigl(\|R_c\|_{L^\infty_{t,x}}+\flf\bigr)}{2L\sqrt{62}\flf}+\frac{(8+\sqrt2)\|\nabla_xR_c\|_{L^\infty_{t,x}}}{\sqrt{62}\flf}\right),
\end{align*}
where $W_1$ is the velocity profile $U_1(x) + \na \big[\chi(|x|)\frac{\xi\cdot x}{|x|^2}\big]$ given in~\eqref{Wr, def}, \eqref{Pi-r, def}, and \eqref{S3}, which becomes fixed once the cutoff function $\chi$ is given. In addition, let $C_{\mathrm{IO}}$ be the universal constant in Lemma~\ref{lemma: IO} following Isett--Oh~\cite{IO}. Also, as before, $G_1=(8+\sqrt{2})\left\|\nabla_xR_c\right\|_{L^\infty_t L^1_x}$. With these symbols at hand, together with those already introduced in~\eqref{parameter, summary}, for any given convolution kernel \(\Omega\) satisfying~\eqref{Omega-h}, we may choose $D$ as follows: 
\begin{align}\label{eq:scan-choice-D} 
 D=\max\Bigg\{&
    1,\ M_3,\ \|\partial_t R_c\|_{L^\infty_tL^1_x},\
    G_1,\ \frac{\bar\eta (8+\sqrt2)\bigl(\|R_c\|_{L^\infty_{t,x}}+\flf\bigr)}{16L\sqrt{62} \flf},\ 32\left(m_0+\varepsilon_0^2/2\right),\nonumber\\
 & 4\Bigl(G_1+\bigl(G_1+128\flf L\bigr)
                       \int_{\R^2} |y|\Omega(y)\,\dd y\Bigr),\
    2\|\nabla\Omega\|_{L^1} m_0, (1+C_{\mathrm{IO}}) F_*,\nonumber\\
 & 2F_*\bigl(\|\Omega\|_{L^2}+\|\nabla\Omega\|_{L^1}\bigr),\
   2F_*\|\leray\|_{L^{3/2}\to L^{3/2}}\|\Omega\|_{L^{3/2}_x}
 \Bigg\}. 
\end{align}

To proceed, let us take \(h,s_i\) exactly as in \eqref{h, choice} and \eqref{si, choice}. Then 
\begin{align*}
D(2\sqrt{h}+h) < \frac{3\varepsilon_0}{400},\qquad \frac{D s_*}{h}<\frac{\varepsilon_\star}{100}.
\end{align*}
With $D$ chosen as above, we specify $\bigtime := \frac{|J|}{M}$ by taking the total number of time cells $M$ so large that 
\begin{align}\label{choice of T-cell}
D\mathsf T_{\mathrm{cell}}<\min\Bigg\{&
     \frac h4,\ \frac{\min_i s_i}{100},\
     \frac{\varepsilon_\star}{100},\
      \left(\frac{\varepsilon_\star}{100(1+2M_3\bar\eta\|U\|_{3/2})}\right)^3,\nonumber\\
 &\qquad \,\frac{h^{3/2} \varepsilon_\star}{100(1+2(M_3+\bar\eta\|U\|_2))},\
        h^{3/2}\sqrt{\frac{\varepsilon_\star}{100}}\Biggr\}.
\end{align}
The parameters \(t_m^{**}\), \(I_{ij}^{(m)}\), and \(\Delta_i^{(m)}\) are then computable from \eqref{t**, def}, \eqref{Iij,m}, and \eqref{S15}. In particular, we note that 
\begin{equation}
2r_*<h\qquad\text{and}\qquad r_*<\min_{1\leq i \leq 4} \frac{s_i}{100},
\end{equation}
where $s_i$ is the width between two consecutive stripes in $Q_i$ determined by~\eqref{condition III}. With all these choices, we infer from~\eqref{R' estimate, pre} that
\begin{equation}\label{R' estimate, final}
\left\|R'\right\|_{L^\infty_tL^1_x}<\varepsilon_0+\frac{\varepsilon_0^2}{2}+\frac{3\varepsilon_0}{400}+\frac{5\varepsilon_\star}{100}
< \frac{\varepsilon_\star}{2}.
\end{equation}
This is the desired bound for $R'$.

\smallskip
\noindent
\underline{\bf Estimate for $u'-u$.} The bound for $\vcorr := u'-u \equiv V+Q$ can be readily deduced. Indeed, by \eqref{S4}, \eqref{eta on stripe ij}, and \eqref{S18}, we have 
\begin{align*}
\sup_{t\in I}\|V(t,\cdot)\|_{L^2_x}+\sup_{t\in I} \|\operatorname{curl}V(t,\cdot)\|_{L^1_x}
 &\le \bar\eta(\|U\|_{L^2_x}+\|\operatorname{curl}U\|_{L^1_x})\\
 &\le C_0\sqrt{\delta}
       +\sqrt{3/\pi}(\|U\|_{L^2_x}+\|\operatorname{curl}U\|_{L^1_x})\varepsilon_0
\\
 & \le C_0\sqrt{\delta}+\varepsilon_\star/16,
\end{align*} 
where $$C_0=(2/\sqrt\pi)\sqrt{8+\sqrt2}
\Big(\|U\|_{L^2}+\|\operatorname{curl}U\|_{L^1}\Big)$$ is universal, and the last
inequality is \eqref{e0}. 

The bound for $Q$ has been proved in~\eqref{eq:scan-Q-bound}, reproduced here:
\begin{align*} 
\|Q\|_{L^\infty_t L^2_x}+\left\|\operatorname{curl}Q\right\|_{L^\infty_t L^1_x}
    \le \frac{D\bigtime}{\theta h},\qquad \|Q\|_{L^\infty_t L^{{3}/{2}}_x} \le \frac{D\bigtime}{\theta h^{2/3}}.
\end{align*}
In the above choice for $\bigtime$, namely~\eqref{choice of T-cell}, we have $$ \frac{D\bigtime}{\theta h} + \frac{D\bigtime}{\theta h^{2/3}} \leq \frac{\e_\star}{100}.$$ This proves
\begin{equation}\label{w estimate, final}
\sup_{t\in I}\left\|\vcorr(t,\cdot)\right\|_{L^2_x} + \sup_{t\in I}\left\|\operatorname{curl}\vcorr(t,\cdot)\right\|_{L^1_x}
                 \le C_0\sqrt\delta+\frac{\varepsilon_\star}{8}. 
\end{equation}

In view of~\eqref{R' estimate, final} and \eqref{w estimate, final}, the proof of Proposition~\ref{propn: raw estimate} is complete.   \end{proof}

\subsubsection{Necessity for further modification}\label{subsubsec: nonsmooth}

The ``raw'' triple \((u',p',R')\) constructed in Lemma~\ref{lemma: update rule for u', p', R'} satisfies the Euler--Reynolds equation~\eqref{ER'} in the distributional sense, as well as the smallness estimates in Proposition~\ref{propn: raw estimate}. However, we cannot conclude the proof of Lemma~\ref{lem: key} directly using $\left(\widetilde{u}, \widetilde{p}, \widetilde{R}\right) = (u',p',R')$, since they are not $\CC^\infty$. Indeed, in the proof of Main Theorem~\ref{mainthm}, we shall iterate the construction $(u,p,R) \rightsquigarrow \left(\widetilde{u}, \widetilde{p}, \widetilde{R}\right)$ in infinitely many stages, so we must return to $\CC^\infty$-triples after each stage.

The failure of smoothness of the triple $(u',p',R')$ occurs at the joins of the glueing construction, \emph{i.e.}, at the endpoints of the time cells. For instance, consider the time corrector $Q$ in~\eqref{Q, def}. By construction, it vanishes at $\p J_m = \left\{t^-_m, t^+_m\right\}$, but its time derivative jumps there:
\begin{align*}
\partial_t Q(t_m^+-)=\leray\overline F_h^{(m)},
\qquad
\partial_t Q(t_m^++ )=\leray \overline F_h^{(m+1)}.
\end{align*}
There is an analogous issue for the stress: In~\eqref{R', final}, the terms involving \(V\) and \(S_c\) vanish near $\p J_m$, and the terms involving \(Q\) tend to zero in \(L^1_x\). Thus,
\begin{align*}
\left(R'-R\right)(t_m^+-)=\biga-b_0^{(m)}{\bf I}_2,
\qquad
\left(R'-R\right)(t_m^++)=\mathbf{A}^{(m+1)}-b_0^{(m+1)}{\bf I}_2.
\end{align*}
These need not agree.

Nevertheless, some regularity properties still pertain for the raw triple $(u',p',R')$. Note that $V$ contains only finitely many profiles (\emph{i.e.}, with traversals of finitely many stripes) of the rescaled Lamb--Chaplygin dipoles $U_{r(t)}$ weighted by $\eta(t)$, where $|U(x)| \lesssim |x|^{-2}$ near  infinity. Moreover, both $V$ and its curl vanish on open ``reset'' intervals near the endpoints of stripes. Hence,  $V \in \CC^0_t \big(L^2_x \cap L^{3/2}_x\big)$ and $\operatorname{curl}V\in \CC^0_t L^1_x$.  
In addition, $Q$ lies in the same function class, thanks to~\eqref{eq:scan-Q-bound} and that $Q=0$ for $t \in \p J_m$; see below~\eqref{Q, def}. In particular, it holds that
\begin{align*}
\vcorr \in \CC^0_t \big(L^2_x \cap L^{3/2}_x\big),\qquad \operatorname{curl}\,\vcorr\in \CC^0_t L^1_x.
\end{align*}

On the other hand, the stress $R'$ and pressure $p'$ may jump at $\p J_m$, since in general they have nonzero trace at the endpoints of the joins. The pressure increment is locally integrable: by construction for $\Lambda$ and $b_0^{(m)}$, $p_V$ is locally bounded.

\subsection{Conclusion}
\label{subsec:scan-final-smoothing}

We now conclude the proof of Lemma~\ref{lem: key} by resolving the nonsmoothness issue of the ``raw'' triple $(u',p',R')$ constructed in Lemma~\ref{lemma: update rule for u', p', R'} via a regularisation argument. 

As a caveat, if one directly mollifies $u'$ and $p'$, the resulting smooth velocity and pressure would change the prescribed Euler flow near the endpoints of the time cells $J_m$, as convolution does not commute with tensor products. Instead, we choose to regularise the \emph{increments}.

 \begin{proof}[Conclusion of the proof of Lemma~\ref{lem: key}]
 
Our arguments are divided into eight steps below.

Let $(u,p,R) \in \CC^\infty$ be the Euler--Reynolds triple as given in Lemma~\ref{lem: key}.  

\smallskip
\noindent
{\bf Step~1: Preparation.}  
Without loss of generality, we may fix $$p(t,0) \equiv 0,$$ \emph{e.g.}, by replacing $p$ with $p_1(t,x):=\int_0^1 x\cdot\na_x p(t,sx)\,\dd s$. We relabel $p_1\equiv p$ in the sequel.

As before, take $$J=[t_I^-,t_I^+]\Subset I=[A,B]$$ and $$\vcorr \equiv V+Q = u'-u$$ as in \eqref{V-corr, def} and \eqref{u',p',R', def}. 
The amplitude $\eta(t)$ in \eqref{eta on stripe ij}
vanishes near the endpoints of each stripe by \eqref{b-theta}, so \(V\) vanishes there too. The time corrector $\mathfrak{q}$ defined in~\eqref{q, def} vanishes at both ends of every cell by \eqref{Fh-bar, def}; hence so does $Q$ in~\eqref{Q, def}. Thus, $\vcorr$ has zero traces at cell joins $\p J_m$ and at $\p J$.

In what follows, set
\begin{align}\label{extension to I}
&\vcorr =0,\qquad p'=p,\qquad R'=R={\bf 0} \qquad \text{for each } t \in \R \setminus I;\nonumber\\
&u(t,\cdot) = u(A,\cdot)\quad \text{ for each } t <A,\qquad u(t,\cdot) = u(B,\cdot)\quad \text{ for each } t > B.
\end{align}
Thus, the estimates in Proposition~\ref{propn: joint L1 convergence of difference in Reynolds stress} and the Euler--Reynolds equation~\eqref{ER'} holds for all time $t \in \R$. Recall from \S\ref{subsubsec: nonsmooth} that
\begin{align*}
\vcorr \in \CC^0_t \big(L^2_x \cap L^{3/2}_x\big)(I \times \R^2),\quad \operatorname{curl}\,\vcorr\in \CC^0_t L^1_x(I \times \R^2),\quad p'-p \in L^1_{\loc}(I \times \R^2).
\end{align*}
  
 Subtracting the  Euler--Reynolds equations for $(u',p',R')$ and $(u,p,R)$, we obtain
\begin{align}\label{PDE for Vcorr}
\p_t\vcorr+\na_x(p'-p)
  +\operatorname{div}
      \big(u\otimes\vcorr+\vcorr\otimes u+\vcorr\otimes\vcorr\big)
   =\operatorname{div}(R'-R).
\end{align}
Note that $\vcorr$, $p'-p$, and $R'-R$ are supported in time in $J$. In light of the extension in~\eqref{extension to I}, this PDE~\eqref{PDE for Vcorr}  holds in $\mathcal{D}'\left(\R\times\R^2\right)$, despite that the equation for $u$ only holds for $t \in I$. 

\smallskip
\noindent
{\bf Step~2.} We now define $\widetilde{u}-u$ and $\widetilde{p}-p$ by mollification.

Take any $\alpha_\e\in\CC_c^\infty\left(\R_t\right)$ and $\beta_\e\in\CC_c^\infty\left(\R^2_x\right)$ nonnegative, supported respectively in $]-\e,\e[$ and $\B_\e(0)$, and satisfying \[
 \int_{\R}\alpha_\e(s)\,\dd s=1,\qquad
 \int_{\R^2}\beta_\e(y)\,\dd y=1.
\] 
Throughout the proof, we require 
\begin{align}\label{first requirement, epsilon}
0<\e<\frac{\operatorname{dist}(J,\p I)}{3}.
\end{align}

Then, for any $(t,x)\in I\times\R^2$, we define the increments:
\begin{align}\label{def for u-tilde, p-tilde}
&(\widetilde u-u)(t,x)
 :=
 \int_{\R}\int_{\R^2}\rho_\e(s,y)
                    \vcorr(t-s,x-y)\,\dd y\,\dd s,\nonumber\\
&(\widetilde p-p)(t,x)
 =
 \int_{\R}\int_{\R^2}\rho_\e(s,y)
                    (p'-p)(t-s,x-y)\,\dd y\,\dd s.
\end{align}
Here and hereafter, we take the convolutions kernel:
\begin{equation*}
\rho_\e(s,y):=\alpha_\e(s)\beta_\e(y).
\end{equation*}

Throughout the remaining parts of the paper, we write 
\(\rho_\e\star\) for the spacetime convolution, \emph{i.e.}, in $(t,x)$-variables. It commutes with distributional derivatives.

\smallskip
\noindent
{\bf Step~3.} By convolving~\eqref{PDE for Vcorr} and adding the equation for $(u,p,R)$, we obtain
\begin{align*}
&\p_t\widetilde u+\na_x\widetilde p
   +\operatorname{div}_x(u\otimes u)\nonumber\\
 &\qquad+\operatorname{div}_x\Bigl[
   \rho_\e\star(u\otimes\vcorr)+\rho_\e\star(\vcorr\otimes u)
                       +\rho_\e\star(\vcorr\otimes\vcorr)\Bigr]\nonumber\\
 &\qquad\qquad=\operatorname{div}_x\bigl[R+\rho_\e\star(R'-R)\bigr]\qquad \text{ in } I \times \R^2.
\end{align*}
This is tantamount to
\begin{equation}\label{tilde-eq}
\p_t\widetilde u+\operatorname{div}_x
(\widetilde u\otimes\widetilde u)+\na_x\widetilde p
=\operatorname{div}_x\widetilde R,
\end{equation}
where
\begin{align}\label{eq:scan-smoothed-stress} 
 \widetilde R &:= R+\rho_\e\star(R'-R) +u\otimes(\widetilde u-u)-\rho_\e\star(u\otimes\vcorr)\nonumber\\
&\qquad+(\widetilde u-u)\otimes u-\rho_\e\star(\vcorr\otimes u) +(\widetilde u-u)\otimes(\widetilde u-u) - \rho_\e\star(\vcorr\otimes\vcorr).
\end{align}
In particular, $\widetilde{R} \in \CC^\infty\left(\R \times \R^2; \sym \right)$.

The regularised velocity \(\widetilde u\) is solenoidal: as both $V$ and $Q$ are divergence-free in the sense of distributions, we have $\operatorname{div}(\widetilde u-u)
   =\rho_\e\star\operatorname{div}(V+Q)=0$, where $u$ is solenoidal by assumption.

For \(t\in[A,t_I^--\e]\cup[t_I^++\e,B]\), we have
\(t-s\notin J\) whenever \(|s|<\e\). The convolutions of
\(\vcorr\), \(p'-p\), \(R'-R\), and the products involving
\(\vcorr\) all vanish there. The terms containing \(\widetilde u-u\)
then vanish as well. Since \(R=0\) on these intervals, we obtain that
\begin{equation}
\left(\widetilde u,\widetilde p,\widetilde R\right)(t,\cdot)
                  =(u,p,0)(t,\cdot)
 \quad\text{for }t\in\left[A,t_I^--\e\right]\cup \left[t_I^++\e,B\right].
\end{equation} 
Both intervals contain a relatively open collar of the corresponding endpoint of \(I\).

\smallskip
\noindent
{\bf Step~4.} Now let us prove that
\begin{equation}\label{XX}
\sup_{t\in I}\left\|\widetilde R(t,\cdot)\right\|_{L^1_x}<\e_\star.
\end{equation}
For this purpose, we control each of the eight terms on the right-hand side of~\eqref{eq:scan-smoothed-stress}.

\smallskip
\noindent
\underline{Terms~1--2.} Write $R+\rho_\e\star(R'-R)
=\rho_\e\star R'+(R-\rho_\e\star R)$.

For the first term, as \(\rho_\e\) has integral~$1$ and is nonnegative, it holds for every \(t\in I\) that
\[
\begin{aligned}
 \|(\rho_\e\star R')(t,\cdot)\|_{L^1_x}
 &\le\int_{\R}\int_{\R^2}\rho_\e(s,y)
                    \|R'(t-s,\cdot-y)\|_{L^1_x}\,\dd y\,\dd s\\
 &\le\operatorname*{ess\,sup}_{s\in I}\|R'(s,\cdot)\|_{L^1_x}
   <\frac{\e_\star}{2}.
\end{aligned}
\]
The final bound follows from Proposition~\ref{propn: raw estimate}. 

For the second term, consider
\begin{align*}
R(t,x)-(\rho_\e\star R)(t,x)=\int_{\R}\int_{\R^2}\rho_\e(s,y)
                  [R(t,x)-R(t-s,x-y)]\,\dd y\,\dd s .
\end{align*} 
Its \(L^1_x\)-norm is bounded by $$\sup_{\substack{t\in I\\|s|,|y|\le\e}}
           \|R(t,\cdot)-R(t-s,\cdot-y)\|_{L^1_x},$$ which tends to zero as $\e \searrow 0$.

\smallskip
\noindent
\underline{Terms~3--6.} In view of the expression~\eqref{def for u-tilde, p-tilde} for $\widetilde u-u$, we may express these terms by
\begin{align*}
&\bigg[u\otimes(\widetilde u-u)-\rho_\e\star(u\otimes\vcorr) +(\widetilde u-u)\otimes u-\rho_\e\star(\vcorr\otimes u)\bigg](t,x)\\
&\quad =  \int_{\R}\int_{\R^2}\rho_\e(s,y)
       [u(t,x)-u(t-s,x-y)]\otimes\vcorr(t-s,x-y)\,\dd y\,\dd s \\
 &\qquad + \int_{\R}\int_{\R^2}\rho_\e(s,y)
       \vcorr(t-s,x-y)\otimes[u(t,x)-u(t-s,x-y)]\,\dd y\,\dd s.
\end{align*}
By H\"{o}lder's inequality, each of the double integrals the right-hand side can be bounded by
\[
\begin{aligned}
 &\int_{\R}\int_{\R^2}\rho_\e(s,y)
        \|u(t,\cdot)-u(t-s,\cdot-y)\|_{L^3_x}
                       \|\vcorr(t-s,\cdot-y)\|_{L^{3/2}_x}
                         \,\dd y\,\dd s\\
 &\quad\le
 \sup_{\substack{t\in I\\|s|,|y|\le\e}}
                  \|u(t,\cdot)-u(t-s,\cdot-y)\|_{L^3_x}
       \ \sup_{t\in I}\|\vcorr(t,\cdot)\|_{L^{3/2}_x}.
\end{aligned}
\]
The second factor is finite; see Step~1 in this proof. The first term tends to zero as $\e \searrow 0$: indeed, \(u\), \(R\), and \(\vcorr\)
are uniformly continuous in time in \(L^3_x\), \(L^1_x\), and \(L^2_x\), respectively. 

\smallskip
\noindent
\underline{Terms~7--8.} To control these two symmetric tensor products, by adding and subtracting $\vcorr\otimes\vcorr$, one may express
\begin{align*}
&(\widetilde u-u)\otimes(\widetilde u-u) - \rho_\e\star(\vcorr\otimes\vcorr)\nonumber\\
&\quad = 
(\widetilde u-u-\vcorr)\otimes(\widetilde u-u) +
\vcorr\otimes(\widetilde u-u-\vcorr)\nonumber\\
&\qquad + \vcorr\otimes\vcorr - \rho_\e \star \left(\vcorr\otimes\vcorr\right)
\end{align*}
For the first line on the right-hand side, we proceed by 
\begin{align}\label{L2 contraction}
 &\sup_{t\in I}\|(\widetilde u-u)(t,\cdot)\|_{L^2_x}
 \le \sup_{t\in I}\|\vcorr(t,\cdot)\|_{L^2_x},\\
 & \sup_{t\in I}\|(\widetilde u-u-\vcorr)(t,\cdot)\|_{L^2_x}
 \le\sup_{\substack{t\in I\\|s|,|y|\le\e}}
                  \|\vcorr(t-s,\cdot-y)-\vcorr(t,\cdot)\|_{L^2_x}
 \longrightarrow 0\nonumber.
\end{align} 
For the second line, inside  the integral defining
that difference evaluated at $(t,x)$ we add and subtract
$\vcorr(t,x)\otimes\vcorr(t-s,x-y)$, finding that the $L^1_x$-norm of this line is   bounded by
\[
\begin{aligned}
 &\int_{\R}\int_{\R^2}\rho_\e(s,y)
       \bigl(\|\vcorr(t-s,\cdot-y)\|_{L^2_x}
                           +\|\vcorr(t,\cdot)\|_{L^2_x}\bigr)\\
 &\hspace{7em}{}\times
         \|\vcorr(t-s,\cdot-y)-\vcorr(t,\cdot)\|_{L^2_x}\,
                                                      \dd y\,\dd s\\
 &\quad\le 2\sup_{t\in I}\|\vcorr(t,\cdot)\|_{L^2_x}
            \sup_{\substack{t\in I\\|s|,|y|\le\e}}
                \|\vcorr(t-s,\cdot-y)-\vcorr(t,\cdot)\|_{L^2_x}.
\end{aligned}
\]
This tends to zero as well.

\smallskip
\noindent
\underline{Summary.} By shrinking $\e$ in the range~\eqref{first requirement, epsilon}, we may require that the contribution from Terms~3--8 on the right-hand side of \eqref{eq:scan-smoothed-stress} is less than $\e_\star / 2$. This together with the contibution from Terms~1--2 implies the desired estimate~\eqref{XX}.

\smallskip
\noindent
{\bf Step~5.}  In \eqref{L2 contraction} we have proved the  \(L^2_x\)-contraction for \(\widetilde u-u\). On the other hand, by commuting curl and convolution we obtain that
\begin{align*}
\curl(\widetilde u-u)(t,x)
  =\int_{\R}\int_{\R^2}\rho_\e(s,y)
                   (\curl\vcorr)(t-s,x-y)\,\dd y\,\dd s.
\end{align*}
Hence,
\begin{align*}
\|\curl_x(\widetilde u-u)(t,\cdot)\|_{L^1_x}
 &\le\int_{\R}\int_{\R^2}\rho_\e(s,y)
         \|(\curl_x\vcorr)(t-s,\cdot-y)\|_{L^1_x}\,\dd y\,\dd s\\
 &\le\sup_{t\in I}\|(\curl_x\vcorr)(t,\cdot)\|_{L^1_x}.
\end{align*} 
We take supremum in $t$ 
 and use~\eqref{eq:scan-raw-increment} in Proposition~\ref{propn: raw estimate} to conclude that
 \begin{equation}
 \sup_{t\in I}\|(\widetilde u-u)(t,\cdot)\|_{L^2_x}
 +\sup_{t\in I}\|\curl_x(\widetilde u-u)(t,\cdot)\|_{L^1_x}\le C_0\sqrt{\delta}+\frac{\e_\star}{8},
 \end{equation}
where $C_0$ is the universal constant as in~\eqref{eq:scan-raw-increment}.

\smallskip
\noindent
{\bf Step~6.} As \(\curl_x\leray=\curl_x\), we have
\begin{align*}
\curl\,\vcorr(t,\cdot) = \curl\Big( \eta(t)W_{r(t)}(\cdot-\cen(t))+\mathfrak{q}(t,\cdot)\Big).
\end{align*} 
By \eqref{square Qi} and \eqref{parametrisation, stripe}, $\cen(t) \1_{\{\eta(t)\neq 0\}}$ lies in $\overline{\B}_{2\sqrt{2}L}(0)$, so $\curl\left( \eta(t)W_{r(t)}(\cdot-\cen(t))\right)$ is supported within a distance of \(2r(t)\le h\) of $\cen(t) \1_{\{\eta(t)\neq 0\}}$, thanks to Lemma~\ref{lem: Wr, Fc}
and \eqref{r, global}. On the other hand, $\mathfrak{q}$ is defined in~\eqref{q, def} as an integral of \(F_h\) and \(\overline F_h^{(m)}\), both of which
are supported in \(\overline{\B}_{2\sqrt{2}L+h}(0)\), thanks to \eqref{Fh, def} and Lemma~\ref{lem: estimates given eta, r}~$(i)$.  Therefore, 
\[
 \operatorname{supp}(\curl \vcorr)(t,\cdot)
                     \subset\overline{\B}_{2\sqrt{2}L+h}(0)
\]
and hence $${\rm supp}\left[\curl (\widetilde u-u)(t,\cdot)\right] \subset \overline{\B}_{2\sqrt{2}L+h+\e}(0)$$ uniformly in \(t\), taking into account the mollification parameter $\e$ in the construction of $\widetilde{u}$.

\smallskip
\noindent
{\bf Step~7.} We now verify that $\left(\widetilde{u}-u\right)(t,\cdot)$ has zero circulation, and that $\curl\left(\widetilde{u}-u\right)(t,\cdot)$ has the vanishing first and second moments. To this end, recall from Lemma~\ref{lem: V-corr} that all  the three vorticity moments of $u'-u$ vanish. Indeed, this holds for all $t \in\R$ in view of the extension of increments in~\eqref{extension to I}.

 These identities pass on to $\widetilde{u}-u$ by convolution and Fubini's theorem, using the change of variables \(z=x-y\). Observe that \begin{align*}
&\int_{\R^2}|z+y|^j|(\curl_x\vcorr)(t-s,z)|\,\dd z\\
 &\qquad\le(2\sqrt{2}L+h+\e)^j
             \sup_{\tau\in I}\|(\curl_x\vcorr)(\tau,\cdot)\|_{L^1_x}\qquad\text{for $|y|<\e$ and $j \in \{1,2\}$},
\end{align*} 
which is finite.
 
 \begin{itemize}
 \item
 For the
circulation, it holds that 
\begin{align*}
\int_{\R^2}\curl(\widetilde u-u)(t,x)\,\dd x =\int_{\R}\int_{\R^2}\rho_\e(s,y)
         \left[\int_{\R^2}(\curl_x\vcorr)(t-s,z)\,\dd z\right]
         \,\dd y\,\dd s=0.
\end{align*}

\item
 For the (vector) first moment, it holds by the above point that
\begin{align*}
 &\int_{\R^2}(z+y)(\curl\vcorr)(t-s,z)\,\dd z\\
 &\qquad =\int_{\R^2}z(\curl\vcorr)(t-s,z)\,\dd z
           +y\int_{\R^2}(\curl\vcorr)(t-s,z)\,\dd z=\mathbf{0}.
\end{align*}
Integrating against \(\rho_\e(s,y)\) yields that
\begin{align*}
\int_{\R^2} x\, \curl(\widetilde u-u)(t,x)\,\dd x = \mathbf{0}.
\end{align*}
 
 \item
 For the (scalar) second moment, it holds by the above two points that
 \begin{align*}
&\int_{\R^2}|z+y|^2(\curl\vcorr)(t-s,z)\,\dd z  =\int_{\R^2}|z|^2(\curl\vcorr)(t-s,z)\,\dd z\\
&\quad+2y\cdot\int_{\R^2}z(\curl\vcorr)(t-s,z)\,\dd z +|y|^2\int_{\R^2}(\curl\vcorr)(t-s,z)\,\dd z=0.
 \end{align*}
Integrating against \(\rho_\e(s,y)\) yields that
\begin{align*}
\int_{\R^2} |x|^2 \curl(\widetilde u-u)(t,x)\,\dd x = 0.
\end{align*}
\end{itemize}

\smallskip
\noindent
{\bf Step~8.} Finally, we check that the updated triple $\left(\widetilde{u},\widetilde{p},\widetilde{R}\right)$ is continuous in time.           

For $p \in ]1,\infty[$, by the definition of $\widetilde u-u$ and that the kernel $\beta_\e$ has integral $1$, we have 
\begin{align*}
& \|(\widetilde u-u)(t,\cdot)\|_{L^p_x}
 \le \int_{\R}\alpha_\e(s)\|\vcorr(t-s,\cdot)\|_{L^p_x}\,\dd s
 \le\sup_{s\in I}\|\vcorr(s,\cdot)\|_{L^p_x},\\
&\|(\widetilde u-u)(t,\cdot)
                 -(\widetilde u-u)(t',\cdot)\|_{L^p_x}\le\int_{\R}\alpha_\e(s)
          \|\vcorr(t-s,\cdot)-\vcorr(t'-s,\cdot)\|_{L^p_x}\,\dd s.          
\end{align*}
The last line tends to zero as $t' \to t$ for any fixed $\e$ since $\vcorr$ is continuous in time; see \S\ref{subsubsec: nonsmooth}. We take $p=2$ and $3/2$ here. In addition, by Young's inequality we have that
\begin{align*}
 \|(\widetilde u-u)(t,\cdot)\|_{L^3_x}
 &\le\|\beta_\e\|_{L^{6/5}_x}
        \int_{\R}\alpha_\e(s)\|\vcorr(t-s,\cdot)\|_{L^2_x}\,\dd s\\
 &\le\|\beta_\e\|_{L^{6/5}_x}
                          \sup_{s\in I}\|\vcorr(s,\cdot)\|_{L^2_x}
\end{align*}
as well as
\begin{align*}
&\|(\widetilde u-u)(t,\cdot)
                 -(\widetilde u-u)(t',\cdot)\|_{L^3_x}\\
 &\quad\le \|\beta_\e\|_{L^{6/5}_x}\int_{\R}\alpha_\e(s)
         \|\vcorr(t-s,\cdot)-\vcorr(t'-s,\cdot)\|_{L^2_x}\,\dd s.
\end{align*}
Again, the last line tends to zero as $t' \to t$ for any fixed $\e$. We thus conclude that $\widetilde u\in\CC^0_tL^3_x$, $\widetilde u-\srad\in\CC^0_tL^2_x$, and $\curl_x\widetilde u\in\CC^0_tL^1_x$.

To see that \(\widetilde R\in\CC^0_tL^1_x\), we recall the decomposition~\eqref{eq:scan-smoothed-stress}, namely
\begin{align*}
 \widetilde R &:= R+\rho_\e\star(R'-R) +u\otimes(\widetilde u-u)-\rho_\e\star(u\otimes\vcorr)\nonumber\\
&\qquad+(\widetilde u-u)\otimes u-\rho_\e\star(\vcorr\otimes u) +(\widetilde u-u)\otimes(\widetilde u-u) - \rho_\e\star(\vcorr\otimes\vcorr).
\end{align*}

\noindent
\underline{Terms~2, 4, 6, 8.}
For the second term, we bound by Fubini and $\int \alpha_\e =1$ that 
\[
\begin{aligned}
 &\|(\rho_\e\star(R'-R))(t,\cdot)
            -(\rho_\e\star(R'-R))(t',\cdot)\|_{L^1_x}\\
 &\quad\le
 \left(\operatorname*{ess\,sup}_{s\in I}
                         \|(R'-R)(s,\cdot)\|_{L^1_x}\right)
             \int_{\R}|\alpha_\e(t-s)-\alpha_\e(t'-s)|\,\dd s\\
 &\quad\le
 \left(\operatorname*{ess\,sup}_{s\in I}
                         \|(R'-R)(s,\cdot)\|_{L^1_x}\right)
                  |t-t'|\int_{\R}|\alpha_\e'(s)|\,\dd s.
\end{aligned}
\]
An analogous computation for Terms 4, 6, and 8 yields that
\begin{align*}
&\bigg\|\left[\rho_\e\star\left( u\otimes\vcorr  +  \vcorr \otimes u + \vcorr \otimes \vcorr\right)\right](t,\cdot) \\
&\qquad-\left[\rho_\e\star\left( u\otimes\vcorr  +  \vcorr \otimes u + \vcorr \otimes \vcorr\right)\right](t',\cdot)  \bigg\|_{L^1_x} \lesssim \left|t-t'\right|
\end{align*}

\noindent
\underline{Terms~3,5,7.} By H\"{o}lder's inequality, we bound
\begin{align*}
&\bigg\|\Big[u\otimes(\widetilde u-u)+(\widetilde u-u)\otimes u\Big](t,\cdot) - \Big[u\otimes(\widetilde u-u)+(\widetilde u-u)\otimes u\Big](t',\cdot)\bigg\|_{L^1_x}\\
&\qquad \leq \|u(t,\cdot)-u(t',\cdot)\|_{L^3_x}
                      \|(\widetilde u-u)(t,\cdot)\|_{L^{3/2}_x} + \|u(t',\cdot)\|_{L^3_x}
       \|(\widetilde u-u)(t,\cdot)
                -(\widetilde u-u)(t',\cdot)\|_{L^{3/2}_x}
\end{align*}
as well as
\begin{align*}
&\bigg\|\Big[(\widetilde u-u)\otimes(\widetilde u-u) \Big](t,\cdot) - \Big[(\widetilde u-u)\otimes(\widetilde u-u)\Big](t',\cdot)\bigg\|_{L^1_x}\\
&\qquad \leq\bigg\{\|(\widetilde u-u)(t,\cdot)\|_{L^2_x}
               +\|(\widetilde u-u)(t',\cdot)\|_{L^2_x}\bigg\}
    \|(\widetilde u-u)(t,\cdot)-(\widetilde u-u)(t',\cdot)\|_{L^2_x}.
\end{align*}
Both of the right-hand sides tend to zero as $t'\to t$, thanks to the arguments at the beginning of this step. Since \(R\) is continuous into \(L^1_x\), we may now conclude the asserted
continuity of \(\widetilde R\).

In view of {\bf Steps~1--8} above, the proof of Lemma~\ref{lem: key} is now complete.  \end{proof}

\begin{remark}
In the proof of Lemma~\ref{lem: key} above, we have selected six families of parameters in the following non-interchangeable order:
\begin{enumerate}
 \item
$\e_0$;
\item
$L,R_c,\underline{\mathtt f},b,a_i$;
\item
$D$;
\item
$h,\theta,N_i,s_i$;
\item
$M,\mathsf T_{\rm cell}$;
\item
$t_m^{**}, I_{ij}^{(m)},\eta_i^{(m)},r_{0i}^{(m)},d_{ij}^{(m)}$. \end{enumerate}  
\end{remark}

\section{Velocity field initialising the iteration scheme: glueing construction}\label{sec: glue}

\subsection{The background radial velocity field}

We begin by choosing the azimuthal velocity field $\srad$ in Theorem~\ref{mainthm}, following the identical construction as in DiPerna--Majda~\cite{dm2, dm3} and Delort~\cite{Delort}. 

Given the initial vorticity $\omega_0\in\mcomp(\R^2)\cap H^{-1}(\R^2)$, let $b \in \CC^\infty_c(\R^2)$ be a radial scalar field such that 
\begin{align}\label{b, def}
    \int_{\R^2} b(x)\,\dd x = \G := \omega_0(\R^2).
\end{align}
Then set 
\begin{equation}\label{srad, def}
    \srad := \kbs \star b,
\end{equation}
where $\kbs$ is the 2D Biot--Savart kernel.

The following has been established in~\cite{dm2, dm3}:

\begin{lemma}\label{lem: srad}
The vector field $\srad$ in \eqref{srad, def} is a stationary solution to the 2D Euler equations, such that $\srad \in \left(\CC^\infty \cap L^p \cap \dot{W}^{1,\infty}\right)(\R^2)$ for some $p>2$, and $|\sigma(x)|\lesssim \min\left\{|x|, |x|^{-1} \right\}$ on $\R^2$ modulo a constant depending only on $\G$. 
\end{lemma}

\subsection{Regularisation of the initial vorticity}

In what follows, we write $$u_0 := \kbs \star \omega_0$$ for the initial velocity, and set $$v_0:=u_0-\srad.$$  We shall also write  $\left|\omega_0\right|(\R^2) \equiv \|\omega_0\|_{\M(\R^2)}$ for the total variation norm of the signed measure $\omega_0$.

We introduce a regularisation of the initial vorticity measure $\omega_0 \in \mcomp(\R^2)$. Instead of simply mollifying $\omega_0$, we approximate it by a smooth scalar field $\gdata$ whose up to the second moments coincide with those of $\omega_0$. The invariance of second moment under regularisation will play a crucial role in Lemma~\ref{lemma: R-initial} below: It ensures that the difference between the Euler velocities with initial vorticity $\gdata$ and $\gdataa$ for arbitrary $\e,\eta>0$ is always torque-free, which enables us to invert the matrix divergence and obtain a Reynolds stress tensor by Lemma~\ref{lemma: IO}.

\begin{lemma}\label{lem: approx of omega-0}
    There exists $\gdata \in \CC^\infty_c(\R^2)$ such that 
    \begin{align*}
    &\left\|\gdata\right\|_{L^1(\R^2)} \leq \left|\omega_0\right|(\R^2) + C_\G\e^2,\\
        &\int_{\R^2} \gdata(x)\,\dd x = \G,\quad \int_{\R^2}x\gdata(x)\,\dd x = \int_{\R^2}x\,\dd\omega_0(x),\\
        &\int_{\R^2}|x|^2\gdata(x)\,\dd x = \int_{\R^2} |y|^2\,\dd\omega_0(y),\\
        &\kbs\star\gdata - \srad \longrightarrow  v_0 \qquad\text{ strongly in } L^2(\R^2;\R^2) \text{ as $\e \to 0$}.
    \end{align*}
Here $C_\G$ is a uniform constant which depends only on the circulation $\G$. 
\end{lemma}

\begin{proof}[Proof of Lemma~\ref{lem: approx of omega-0}]
Let $\varrho \in \CC^\infty_c(\R^2)$ be any nonnegative radial kernel with $\int_{\R^2}\varrho(x)\,\dd x =1$, and write $\varrho_\e(x):= \e^{-2} \varrho(x/\e)$. This shall be fixed once and for all, so we suppress the dependence of variables on $\varrho$ in the sequel.

 Define
\begin{equation}\label{def of omega,0,e,good}
    \gdata := \varrho_\e\star\omega_0 - \frac{1}{4} \e^2 \G \left(\int_{\R^2}|x|^2\varrho(x)\,\dd x\right) (\Delta\varrho).
\end{equation}
Since $\int_{\R^2}\Delta\varrho\,\dd x = 0$ and $\int_{\R^2}x \Delta\varrho (x)\,\dd x = {\bf 0}$, the circulation and the first moment of $\gdata$ equal those of $\omega_0$. Also, it is clear from \eqref{def of omega,0,e,good} that $\left\|\gdata\right\|_{L^1(\R^2)} \leq \left|\omega_0\right|(\R^2) + C_\G\e^2$ with $C_\G:= 4^{-1}|\G|\left(\int_{\R^2}|x|^2\varrho(x)\,\dd x\right)\|\Delta\varrho\|_{L^\infty_x}$.

For the second moment, using the identity $\int_{\R^2}|x|^2 \Delta\varrho(x)\,\dd x = 4$, one deduces that
\begin{align*}
    \int_{\R^2}|x|^2\gdata(x)\,\dd x &= \int_{\R^2} |x|^2(\varrho_\e \star\omega_0)\,\dd x - \G \e^2 \left(\int_{\R^2}|x|^2\varrho(x)\,\dd x\right)\\
    &=\e^2 \iint_{\R^2 \times \R^2} \left(|z|^2 - \left|z-\frac{y}{\e}\right|^2\right) \varrho\left(z-\frac{y}{\e}\right)\,\dd\omega_0(y)\,\dd z\\
    &= \iint_{\R^2 \times \R^2} |y|^2\varrho\left(z-\frac{y}{\e}\right)\,\dd\omega_0(y)\,\dd z\\
    &=\int_{\R^2} |y|^2\,\dd\omega_0(y).
\end{align*}

Finally, we have 
\begin{align*}
    \kbs\star\gdata - \srad = \kbs \star \varrho_\e \star\omega_0 -\srad - \frac{1}{4}\e^2 \G \left(\int_{\R^2}|x|^2\varrho(x)\,\dd x\right) \na^\perp\varrho.
\end{align*}
Here $ \kbs \star \varrho_\e \star\omega_0 -\srad  \to \kbs\star\omega_0 - \srad \equiv v_0$ strongly in $L^2$, while the other term equals $\e^2$ times a smooth, compactly supported function. To see that \(v_0 \in L^2\), note that for \(\nu = \omega_0 - b\), compact support and \(\nu(\mathbb{R}^2) = 0\) lead to \(|\widehat{\nu}(\xi)| \le |\xi| \int |x| \,\mathrm{d}|\nu|(x)\). In some neighbourhood of zero, this expression bounds the Fourier transform of \(K_{\mathrm{BS}} \star \nu\), and at high frequencies its \(L^2\)-bound follows from \(\nu \in H^{-1}\). Thus the assertion follows.  \end{proof}

\subsection{Smooth Euler flow with regularised initial vorticity}

Now, we show that the smooth Euler flow $\ugood$ whose initial vorticity is $\gdata$ introduced in the previous subsection has nice convergence properties. See Lemma~\ref{lemma:joint trace} and Proposition~\ref{propn: joint L1 convergence of difference in Reynolds stress} below.

Let $\gsol$ be the global smooth 2D Euler vorticity with initial datum $\gdata$:
\begin{equation}\label{eq: gsol}
\begin{cases}
    \p_t\gsol + \left(\kbs \star \gsol\right) \cdot\na \gsol = 0,\\
    \gsol\big|_{t=0} = \gdata,
    \end{cases}
\end{equation}
and set
\begin{equation}\label{def, v-e}
    \ugood := \kbs \star \gsol,\qquad v_\e := \ugood - \srad. 
\end{equation}
Since $\ugood$ is a 2D Euler flow and, as before, $\srad$ is a stationary 2D Euler flow, we obtain the equation for $v^\e$: there exists a scalar field $p_\e$ such that
\begin{equation}\label{eq for v-e}
    \begin{cases}
        \p_t v_\e + v_\e\cdot \na v_\e + \srad\cdot\na v_\e + v_\e \cdot \na \srad + \na p_\e = 0,\\
        v_\e \big|_{t=0} = v_{\e,0}:=\kbs\star\gdata - \srad.
    \end{cases}
\end{equation}

\begin{lemma}\label{lemma:joint trace}
Jointly as $\e,t \searrow 0$, it holds that
\begin{enumerate}
    \item 
    $v_{\e}(t,\cdot) \to v_0$ strongly in $L^2_x$;
    \item 
    $\gsol \weakstar \omega_0$ as Radon measures; and
    \item 
    $\left\{|\gsol(t,\cdot)|\right\}$ are uniformly tight. 
\end{enumerate} 
\end{lemma}

\begin{remark}
    Here and hereafter, by a slight abuse of notation, we shall always identify $\gsol(t,\cdot)$ with the measure $\gsol(t,\cdot)\leb^2$, and similarly for any other locally integrable functions. 
\end{remark}

\begin{proof}[Proof of Lemma~\ref{lemma:joint trace}]

We first argue for Point~(2). By definition of $\gdata$ in~\eqref{def of omega,0,e,good}, we have $\gdata \weakstar \omega_0$ as Radon measures as $\e \searrow 0$. Also, as $t \searrow 0$,  we have $\gsol(t,\cdot) \to \gdata$ in $\CC^k(\R^2)$ for any $k \in \mathbb{N}$, by the initial trace for smooth 2D Euler flows. To show the joint convergence in $\e,t \searrow 0$, let us prove the precompactness of $\left\{\gsol\right\}$ in $L^1_t \M_x\left([0,T] \times \R^2\right)$, which would follow from the $L^1_{t,x}$-boundedness of $\left\{\gsol\right\}$ and the $L^1_t W^{-2,1}_x$-boundedness of $\left\{\p_t \gsol\right\}$, thanks to the Aubin--Lions lemma. 

To this end, take an arbitrary test function $\varphi = \varphi(x) \in W^{2,\infty}_c(\R^2)$. In view of~\eqref{eq: gsol}, integration by parts, symmetrisation over $x$ and $y$, as well as $\kbs(z) = \frac{z^\perp}{2\pi|z|^2}$, we obtain that
\begin{align*}
    \left|\left\langle \p_t\gsol(t,\cdot),\varphi \right\rangle \right| &= \left| \iint_{\R^2 \times \R^2} \omega(t,x)\omega(t,y)   \left[\frac{1}{2} \kbs (x-y)\cdot \Big(\na\varphi(x) - \na\varphi(y)\Big) \right] \,\dd y\,\dd x \right| \\
    &\lesssim \left\|\gsol(t,\cdot)\right\|_{L^1_x}^2 \|\varphi\|_{W^{2,\infty}_x}
\end{align*}
modulo a universal constant. We may now conclude (2).

Next, for Point~(1), we apply standard energy estimates to \eqref{eq for v-e} to deduce that 
\begin{equation*}
\left\|v_\e(t,\cdot)\right\|^2_{L^2_x} \leq e^{2 \|\na \srad(t,\cdot)\|_{L^\infty_x} \, t} \|v_{\e,0}\|^2_{L^2_x}.
\end{equation*}
The boundedness of the $L^\infty_x$-norm of $\na \srad$ follows from Lemma~\ref{lem: srad}, and $v_{\e,0} \to v_0$ strongly in $L^2$ by Lemma~\ref{lem: approx of omega-0}. Hence, 
\begin{equation}
\label{norm bound}
\limsup_{\varepsilon,t\to0}\|v_\varepsilon(t,\cdot)\|_{L^2_x}\le\|v_0\|_{L^2_x}.
\end{equation}
Moreover, by Point~(2) proved above, any weak \(L^2_x\)-subsequential limit of $\{v_\varepsilon(t,\cdot)\}$ has the same curl and divergence as \(v_0\), so their difference must be zero. This together with~\eqref{norm bound} yields the strong convergence claimed in Point~(1).

Finally, as the smooth Euler flow transports vorticity by a volume-preserving flow, $$\left\|\omega^{\rm good}_{0,\varepsilon}\right\|_{L^1_x}
\le |\omega_0|(\mathbb R^2)+C_\Gamma\varepsilon^2.$$ Integration over time yields that
\begin{align*}
\left|\bra \gsol(t,\cdot) - \gsol(0,\cdot),\varphi\ket\right|
\le \frac{t}{4\pi}\left\|D^2\varphi\right\|_{L^\infty_x}
        \left\|\omega^{\rm good}_{0,\varepsilon}\right\|_{L^1_x}^{\,2}\qquad\text{for each $t>0$ and $\varphi\in \CC_c^\infty(\mathbb R^2)$}.
\end{align*}
The right-hand side tends to zero uniformly in $\e$.  Combining this with Point~(2) above and the uniform mass bound proves the joint weak-star convergence, and hence
\begin{align*}
\int_{\R^2} \left|\gsol(t,x)\right|\,\dd x \longrightarrow \left|\omega_0\right|(\R^2)
\end{align*}
thanks to the lower semicontinuity for the weak-star convergence of measures.

Now, given any $\eta>0$, choose $\psi\in \CC_c^\infty(\R^2)$ such that $|\psi|\le 1$ and with
$\int_{\R^2}\psi\,\dd \omega_0> \left|\omega_0\right|(\R^2)-\eta$. We then have
\begin{align*}
    \int_{\R^2\setminus {\rm supp}(\psi)} \left|\gsol(t,x)\right|\,\dd x &\leq \left\|\gsol(t,\cdot)\right\|_{L^1_x(\R^2)} - \int_{\R^2}\gsol(t,x)\psi(x)\,\dd x \\
    &\leq \eta + C_\G \e^2 + \int_{\R^2}\psi \,\dd\left[\gsol(t,\cdot) - \omega_0 \right],
\end{align*}
where the last term tends to zero by Point~(2). By  choosing first $\eta$ and then  $(\e,t)$ sufficiently small, we conclude  the uniform tightness of $\left\{\gdata(t,\cdot)\right\}$, \emph{i.e.}, Point~(3).        \end{proof}

Recall that $\ugood$ is the global smooth 2D Euler flow whose initial vorticity is $\gdata$. We now prove that the sequence of tensors $\left\{\ugood \otimes \ugood\right\}$ is Cauchy in the strong $L^1$-topology modulo tensors whose divergences are gradients. More precisely:

\begin{proposition}\label{propn: joint L1 convergence of difference in Reynolds stress} Given any $\e,\eta>0$, there exists a tensor field $\mathcal{T}_{\e,\eta} \in L^1\left([0,\infty [\times \R^2; \, \sym\right)$ with ${\rm Div}\,\mathcal{T}_{\e,\eta} = \na \pi$ for some scalar field $\pi$, such that 
\begin{equation}\label{dee, def}
    \dee := \ugood \otimes \ugood - u^{\rm good}_\eta \otimes u^{\rm good}_\eta + \mathcal{T}_{\e,\eta} \longrightarrow 0 \qquad \text{in $L^1$ uniformly jointly as $\e,\eta, t \searrow 0$}.
\end{equation}
\end{proposition}

\begin{proof}[Proof of Proposition~\ref{propn: joint L1 convergence of difference in Reynolds stress}]

Recall from~\eqref{def, v-e} that $v_\e=\ugood - \srad$, where $\srad$ is a background radial vector field independent of $\e$. We shall work with $v_\e$ instead of $\ugood$, since $\{v_\e(t,\cdot)\}$ is strongly convergent in $L^2$ as $\e,t\searrow 0$. See Lemma~\ref{lemma:joint trace} above.

We shall exhibit $$\mathcal{T}_{\e,\eta} \equiv \mathcal{T}'\left[\gsol-\gsoll\right] - \srad \otimes (\kbs\star\nu) + (\kbs\star\nu)\otimes \srad,$$ where $\mathcal{T}'$ maps finite signed Radon measures on $\R^2$ to $L^1(\R^2)$, such that
\begin{itemize}
    \item[(i)]
    $\mathcal{T}'$ maps uniformly tight, weakly convergent sequences of signed Radon measures to strongly convergent $L^1$-tensors;
    \item[(ii)] For some scalar field $\pi$, it holds that
    \begin{equation}
   {\rm div}\left(\mathcal{T}'\left[\nu\right]\right) = {\rm div}\left[\srad \otimes (\kbs\star\nu) + (\kbs\star\nu)\otimes \srad\right] + \na \pi.     
    \end{equation}
\end{itemize}
The assertion then follows from the above properties of $\mathcal{T}_{\e,\eta}$ and Lemma~\ref{lemma:joint trace}. 

For this purpose, choose any cutoff function~$\chi \in \CC^\infty([0,\infty[)$ such that $\chi \equiv 0$ on $[0,1/2]$, $\chi \equiv 1$ on $[2,\infty[$, and $0 \leq \chi \leq 1$ everywhere. Define
\begin{align}\label{T', def}
    \mathcal{T}'[\nu](x) &:= \int_{\R^2} \Bigg\{\srad(x) \otimes \kbs(x-y) + \kbs(x-y)\otimes \srad(x)  \nonumber\\
    &\qquad \qquad \quad -\chi\left(\frac{|x|}{|y|} \right) \Big[ \srad(x)\otimes \kbs(x) + \kbs(x)\otimes \sigma(x) \Big] \Bigg\}\,\dd\nu(y).
\end{align}
For ease of future reference, write the integrand as $F_y(x)$ for the moment.

To verify Point~(i) above, it suffices to check that $y \mapsto F_y$ is continuous into $L^1$ with $\sup_{y \in \R^2}\|F_y\|_{L^1_x(\R^2)} < \infty$. Consider three cases separately:

\noindent
\underline{Case~1: $|x| \leq |y|/2$.} In this case, the cutoff function $\chi\left(\frac{|x|}{|y|} \right) \equiv 0$, $|\srad(x)| \lesssim \min\left\{|y|, |x|^{-1} \right\}$, and $\left|\kbs(x-y)\right| \lesssim |y|^{-1}$. Thus, 
\begin{align*}
\int_{\{|x|< |y|/2\}} |F_y(x)|\,\dd x &\lesssim \frac{1}{|y|}\int_{\{|x|< |y|/2\}} \min\left\{|y|, |x|^{-1} \right\}\,\dd x \\
&\lesssim \min \left\{|y|^2,1\right\}.
\end{align*}

\noindent
\underline{Case~2: $|y|/2 < |x| < 2|y|$.} In this case, $\left|\kbs(x-y)\right| \lesssim |x-y|^{-1}$ and $|\srad(x)| \lesssim \min\{|y|, |y|^{-1}\} \approx  \min\{|x|, |x|^{-1}\}$. As $|x-y|^{-1}$ is locally integrable on 2D, we thus have
\begin{align*}
\int_{\{|y|/2<|x|<2|y|\}} |F_y(x)|\,\dd x &\lesssim \int_{\{|y|/2<|x|<2|y|\}} \left(\frac{\min\left\{|x|, |x|^{-1} \right\}}{|x-y|} +\frac{\min\left\{|x|, |x|^{-1} \right\}}{|x|}\right) \,\dd x\\
&\lesssim \min \left\{|y|^2,1\right\}.
\end{align*}

\noindent
\underline{Case~3: $|x| \geq 2 |y|$.} Note that $\left|\kbs(x-y)-\kbs(x)\right| \lesssim \frac{|y|}{|x|^2}$ and $|\srad(x)|\lesssim \min\{|x|, |x|^{-1}\}$, so we may estimate in polar coordinates:
\begin{align*}
    \int_{\{|x|\geq 2|y|\}} |F_y(x)|\,\dd x \lesssim |y|\int_{2|y|}^\infty \frac{\min\{s,s^{-1}\}}{s^2} s\,\dd s\lesssim \min \left\{|y|,1\right\}.
\end{align*}

By the three cases above, $\|F_y\|_{L^1_x}$ is uniformly bounded for $|y| \leq C$ with any fixed constant $C>0$, and $\|F_y\|_{L^1_x} \lesssim |y|^2 \to 0$ as $y \to 0$. In particular, $y \mapsto F_y$ is continuous into $L^1$ as $y \to 0$. In view of the expression~\eqref{T', def} for $F_y$, we observe that the above estimates together with the $L^1_\loc$-continuity of the translation $\kbs(\cdot - y)$ imply the continuity of $y \mapsto F_y$ at each $y \neq 0$. Hence Point~(i) follows.

Now let us turn to Point~(ii). Observe that $\srad$ and $\kbs$ are radial vector fields whose value at each point $x \in \p\B_1(0)$ is perpendicular to $x$. By radial symmetry, $\srad$ and $\kbs$ are tangential to $\p\B_r(0)$ for each $r>0$. Thus, in view of the definition of $\mathcal{T}'[\nu]$ in~\eqref{T', def}, in polar coordinates $(\rho,\theta)$ one may express
\begin{align*}
 \Big(   \mathcal{T}'[\nu]  - \srad \otimes (\kbs\star\nu) - (\kbs\star\nu)\otimes \srad\Big) (x) = A(\rho)\, {\bf e}_\theta \otimes {\bf e}_\theta.
\end{align*}
It follows that
\begin{align*}
    &{\rm Div}\Big( \mathcal{T}'[\nu] - \srad \otimes (\kbs\star\nu) - (\kbs\star\nu)\otimes \srad\Big)(x) \\
    &\qquad \qquad = -\frac{A(\rho)}{\rho}{\bf e}_r = \na \left\{ -\int_1^\rho \frac{A(\rho')}{\rho'}\,\dd\rho' \right\}.
\end{align*}
The proof is now complete.    \end{proof}

\begin{lemma}\label{lemma: R-initial}
    Given any $0<\e,\eta\ll 1$, there exists a time-independent tensor field $\rinit \in L^1\left(\R^2; \sym\right)$ such that
    \begin{align*}
        & {\rm div}\,\rinit = \kbs\star \left(\gdata - \gdataa\right),\\
        & \left\|\rinit\right\|_{L^1_x} \leq C\max\{\e,\eta\},\\
        & \ugood(t,x) - \ugoodd(t,x) = {\rm div}\left\{ \rinit - \int_0^t \dee(s,x)\,\dd s \right\} + \na \pi_{\e,\eta}(t,x).
    \end{align*}
    Here, $C>0$ is a universal constant, and $\pi_{\e,\eta}$ is a scalar field.
\end{lemma}

For future references, we introduce:
\begin{equation}\label{B, def, reynold}
    \mathscr{B}_{\e, \eta}(t,x) := \rinit - \int_0^t \dee(s,x)\,\dd s.
\end{equation}
It is crucial that the constant $C$ here is universal. It depends only on the Biot--Savart kernel.

\begin{proof}[Proof of Lemma~\ref{lemma: R-initial}]
Consider the vector field $$\mathfrak{f}_0:=\kbs\star \left(\gdata - \gdataa\right).$$

\begin{itemize}
    \item 
    It is smooth and compactly supported in a fixed ball. Indeed, by the definition of $\gdata$ and $\gdataa$ in~\eqref{def of omega,0,e,good}, it holds that $\mathfrak{f}_0 = \kbs\star \left[(\varrho_\e - \varrho_\eta) \star \omega_0\right]$. Thanks to the mean value property for harmonic functions, we have $\kbs\star \varrho_\e \equiv \kbs$ outside the support of $\varrho_\e$, and similarly for $\varrho_\eta$. Hence, as $\|\kbs\star\varrho_\e - \kbs\|_{L^1_x} \lesssim \e$ and similarly for $\varrho_\eta$, one has $\left\| \mathfrak{f}_0\right\|_{L^1_x} \lesssim \max \{\e,\eta\}$  
modulo a universal constant. 
\item 
It is divergence-free and compactly supported, so it has zero mean: As ${\rm div}\,\mathfrak{f}_0 = 0$, for each $j \in \{1,2\}$ we have ${\rm div}(x^j \mathfrak{f}_0) = \mathfrak{f}_0^j$. But the vector field $x^j \mathfrak{f}_0$ is compactly supported, so $\int_{\R^2}\mathfrak{f}_0(x)\,\dd x = 0$ by Gauss's theorem.

\item 
It has zero torque. Indeed, by integration by parts, $\curl^{-1} = \kbs\star$, and the fact that the second moments of $\gdata$ and $\omega_0$ coincide for any $\e>0$ (Lemma~\ref{lem: approx of omega-0}), we deduce:
\begin{align*}
    -2\,{\rm Torque}(\mathfrak{f}_0) &= 2\int_{\R^2}\Big[ -x^1 \mathfrak{f}_0^2(x) + x^2 \mathfrak{f}_0^1(x) \Big]\,\dd x\\
    &= -2\int_{\R^2} x^\perp\cdot \mathfrak{f}_0(x)\,\dd x\\
    &= \int_{\R^2} |x|^2 \curl\,\mathfrak{f}_0(x)\,\dd x\\
    &= \int_{\R^2} |x|^2 \left(\gdata (x) - \gdataa (x) \right)\,\dd x\\
    &=0.
\end{align*}

\end{itemize}

By virtue of Lemma~\ref{lemma: IO}, we can find compactly supported $\rinit \in L^1\left(\R^2; \sym\right)$ such that ${\rm div}\,\rinit = \mathfrak{f}_0$ and $\left\|\rinit\right\|_{L^1_x} \lesssim \|\mathfrak{f}_0\|_{L^1_x}$. Subtracting the Euler equations for $\ugood$ and $\ugoodd$, we deduce the PDE for $\mathfrak{f}:=\ugood - \ugoodd$ as follows:
\begin{align*}
    0 &= \p_t \mathfrak{f} + {\rm div}\,\Big(\ugood \otimes \ugood - \ugoodd\otimes \ugoodd\Big)+ \na \pi\\
    &\qquad = \p_t \mathfrak{f} + {\rm div}\, \left(\dee\right) + \na \pi'
\end{align*}
for some scalar fields $\pi$ and $\pi'$.  See Proposition~\ref{propn: joint L1 convergence of difference in Reynolds stress} for the second line above. Also note that $\mathfrak{f}\big|_{t=0} = \mathfrak{f}_0$, where the trace is taken with respect to the strong $L^1$-topology; see Lemma~\ref{lem: approx of omega-0}. Integrating the above PDE over time, we conclude that $$\ugood(t,x) - \ugoodd(t,x) = {\rm div}\left\{ \rinit - \int_0^t \dee(s,x)\,\dd s \right\} + \na \pi_{\e,\eta}(t,x)$$ for some scalar field $\pi_{\e,\eta}$.

The proof is complete.  \end{proof}

\subsection{Temporal glueing}
Now let us describe the glueing scheme in time. 

Given any $\e_\star>0$ as in the statement of Theorem~\ref{mainthm}. Fix a small number $\kappa \in ]0,1]$ to be specified. (In fact, $\kappa$ will be chosen to depend only on $\e_\star$; see~\eqref{kappa, choice} below.)

 Set for each $n \in \mathbb{N}$:
\begin{equation}\label{Dn, def}
    D_n := \frac{\kappa^2}{4^n}.
\end{equation}
Let $\{a_n\} \searrow 0$ be any sequence decreasing to zero; \textit{e.g.}, $$a_n = \frac{1}{2^n}$$ will do the work. Throughout, $n$ indexes the time intervals used in our glueing construction.

Applying Lemma~\ref{lemma:joint trace} and Proposition~\ref{propn: joint L1 convergence of difference in Reynolds stress}, for each $n \in \mathbb{N}$ we obtain some $h_n>0$ such that
\begin{align}\label{xx}
    & \left\|v_\e(t,\cdot) - v_0\right\|_{L^2_x} \leq a_n, \nonumber\\
    &\left\|\ugood(t,\cdot) - \ugoodd(t,\cdot) \right\|_{L^2_x} \leq c_{\rm Rey} \sqrt{D_n}, \nonumber\\
    &\left\|\dee(t,\cdot)\right\|_{L^1_x} \leq c_{\rm Rey} D_n \qquad \text{whenever } \e,\eta,t \in ]0,h_n]. 
\end{align}
Here $c_{\rm Rey}>0$ is a universal constant independent of $n$, which is also to be determined later.

We choose points of times $t_n$ such that
\begin{equation}\label{choice of tn}
t_1 < \frac{1}{2},\qquad t_{n+1}< \frac{t_n}{4},\qquad \text{and }  t_n < \frac{h_n}{2}.
\end{equation}
Denote the $n^{\text{th}}$-epoch by
\begin{equation*}
    I_n:= ]t_n,2t_n[.
\end{equation*}

Next, we choose $\e_n \searrow 0$, which act as the approximation parameters for $\omega_0$. That is, $\e_n$ and $\e_{n+1}$ shall play the role of $\e$ and $\eta$ in $\gdata$ and $\gdataa$, respectively. Let us pick
\begin{equation}\label{choice of en}
    \begin{cases}
        \e_1 \leq \min\left\{ \frac{\sqrt{\kappa}}{C_\G}, h_1,c_{\rm Rey} t_1 D_1  \right\};\\
        \e_n \leq \min\left\{h_n, h_{n-1}, c_{\rm Rey}t_n D_n\right\}\qquad\text{for each } n \geq 2.
    \end{cases}
\end{equation}
Then, let $\zz{n}$ be the global smooth 2D Euler flow whose initial vorticity is $\gdatan$; see Lemma~\ref{lem: approx of omega-0} and formula~\eqref{def of omega,0,e,good} for its definition. $C_\G$ is the  universal constant as in Lemma~\ref{lem: key}.

\begin{figure}[htbp]
\centering\includegraphics[width=\linewidth]{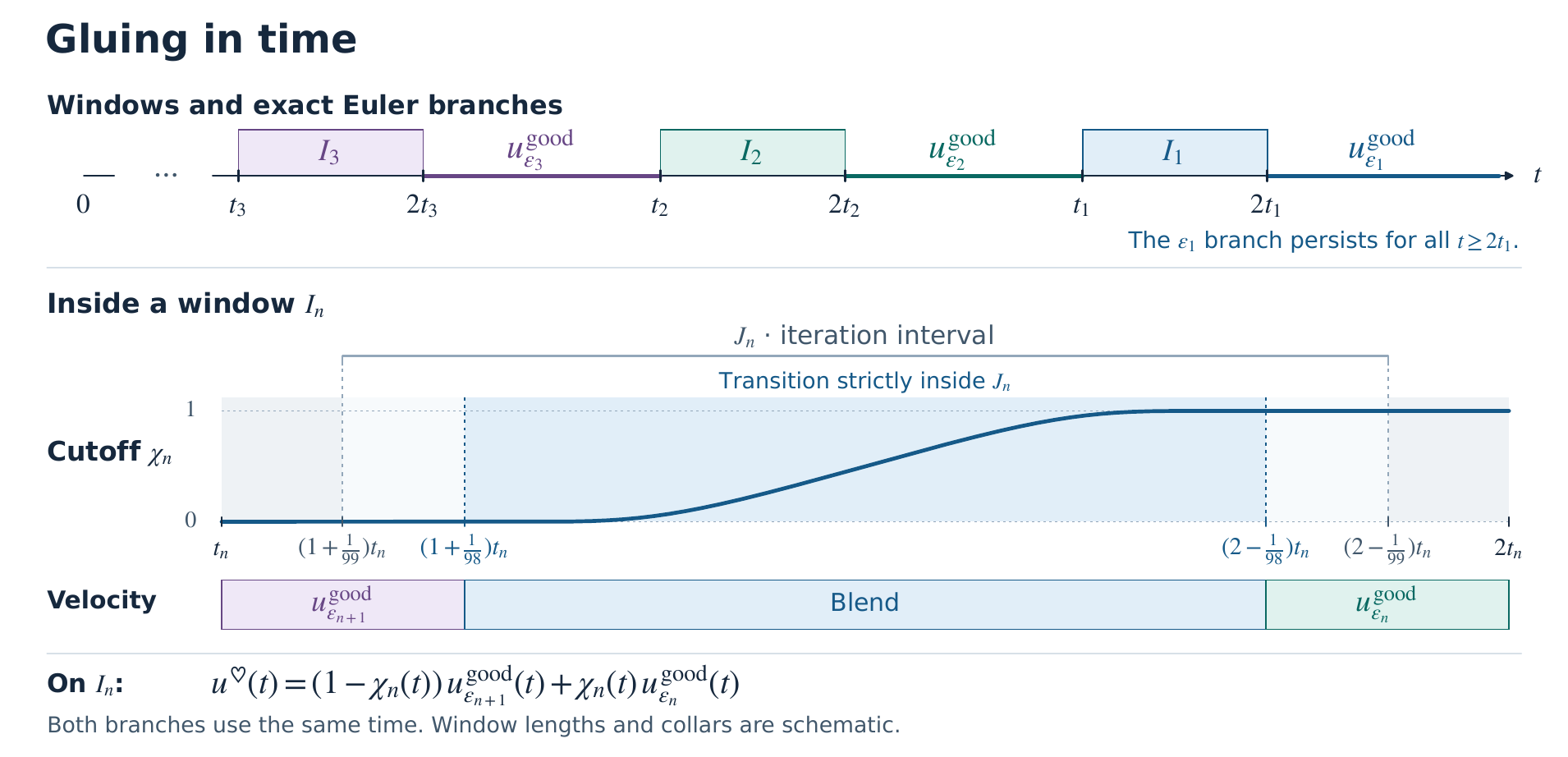}
\caption{Schematic presentation for glueing in time.}
\label{fig:glue}
\end{figure}

We are now at the stage of defining the glued velocity $\glue$. Roughly speaking, we construct $\glue = \zz{n}$ on $I_n$ except at the beginning of the epoch $I_n$. See Figure~\ref{fig:glue}. In the subsequent~\S\ref{sec: proof}, this velocity field $\glue$ will serve as the starting point of the iterative constructions.

More precisely, choose (temporal) cutoff functions $\chi_n = \chi_n(t)$ such that 
\begin{align}\label{def, chi-n}
&\chi_n \in \CC^\infty([0,\infty[),\qquad \chi_n \equiv 0 \text{ on } \left[0,\left(1+\frac{1}{100}\right)t_n\right],\nonumber\\
&\chi_n \equiv 1 \text{ on } \left[ \left(2-\frac{1}{100}\right)t_{n},\infty \right[,\qquad \chi_n' \leq \frac{400}{t_n}\qquad\text{for each $n \in \mathbb{N}$}.
\end{align}
Then define 
\begin{equation}\label{u-heart, key def, glue}
    \glue (t,x):= \sum_{n=1}^\infty \bigg\{ \chi_n(t) \zz{n}(t,x) + \Big[1-\chi_n(t)\Big] \zz{n+1}(t,x) \bigg\}\1_{I_n}(t)
\end{equation}
for each $t >0$ and $x \in \R^2$. 
\begin{proposition}\label{propn: glue}
The glued velocity satisfies the following:
\begin{enumerate}
    \item 
    $\glue(t,\cdot)\equiv \zz{n}(t,\cdot)$ for each $t \in \left[ \left(1+\frac{1}{100}\right)t_n,2t_n\right] \subset I_n$ for each $n \in \mathbb{N}$;
    \item 
    $\glue \in \CC^\infty\left(]0,\infty[ \times \R^2\right)$;
    \item 
    $\glue\big|_{t=0} = u_0 :=\kbs \star \omega_0$ in the strong $L^2$-topology. 
    
    \item 
$\glue \in \CC^0_t L^3_x\left(I_n \times \R^2\right)$ for each $n$.
    
\item 
$\omega^\heartsuit := \curl\,\glue$ has the spatial $L^1$-bounds: for each $n$,
\begin{equation}\label{eq:s4-base-tv}
\|\omega^\heartsuit(t,\cdot)\|_{L^1_x}
 \leq
 \begin{cases}
\left|\omega_0\right|(\R^2) +C_{\G}\e_n^2, & t\in\overline {I_n} = [t_n,2t_n],\\
\left|\omega_0\right|(\R^2) +C_{\G}\e_{n+1}^2, & t \in [2t_{n+1}, t_n],\\
\left|\omega_0\right|(\R^2) +C_{\G}\e_1^2, & t\geq 2t_1,
 \end{cases}
\end{equation}
where $C_\G$ is the constant from Lemma~\ref{lem: approx of omega-0}.
    
    \item 
    Let $\mathscr{B}_{\e_n, \e_{n+1}}$ be the tensor field defined as in \eqref{B, def, reynold} above. Reynolds stress tensor $\rey$ for the glued velocity field satisfying the equation
    \begin{equation}\label{Euler eq for heart variables}
        \p_t\glue + \glue \cdot \na \glue + \na \pi^\heartsuit = -{\rm Div}\,\rey
    \end{equation}
    has the following expression:
    \begin{align}\label{expression for R-heart}
        \rey (t,x) &=\chi_n'(t) \mathscr{B}_{\e_n, \e_{n+1}}(t,x) \nonumber\\
        &\quad- \chi_n(t)\Big[1-\chi_n(t)\Big]\left(\zz{n}(t,x)-\zz{n+1}(t,x)\right) \otimes \left(\zz{n}(t,x)-\zz{n+1}(t,x)\right).
    \end{align}
    The support of $\rey$ lies in $\bigsqcup_{n=1}^\infty \left[\left(1+\frac{1}{98}\right)t_n, \left(2-\frac{1}{98}\right)t_n \right] $.

    \item 
      By choosing $\crey$ to be a sufficiently small universal constant, we have that 
    \begin{equation*}
        \left\|\rey\right\|_{L^\infty_t L^1_x\left(I_n \times \R^2\right)} \leq D_n\qquad\text{ for each } n \in \mathbb{N}.
    \end{equation*}
\end{enumerate} 
\end{proposition}

\begin{proof}[Proof of Proposition~\ref{propn: glue}]
The first point follows from the definition of $\glue$ in \eqref{u-heart, key def, glue} and the choice of cutoff functions $\chi_n$ in \eqref{def, chi-n}.

The second point holds since $\ugood$ for any $\e>0$ is smooth for positive time, and the glueing in \eqref{u-heart, key def, glue} is also smooth.

The third point follows from Lemma~\ref{lemma:joint trace}~(1). 

For the fourth point, note that on each epoch $I_n$, thanks to the structure of the  Biot--Savart kernel and the  uniform $L^1$-bound of the regularised vorticity, $\left|\glue (\cdot,x)\right| \leq \mathcal{O}(|x|^{-1})$ outside twice a ball
containing the support of its vorticity, which is uniformly small in $L^3(\R^2 \setminus \B_{R})$ as $R\to\infty$. Also, by the choice of temporal cutoff functions in \eqref{def, chi-n}, the glued velocity is continuous in time.

For the fifth point, the estimate  on $\overline{I_n}$ follows from the
conservation of $\|\omega_\e^{\rm good}(t,\cdot)\|_{L^1_x}$ and $\e_{n+1}\leq \e_n$. For $t\in [2t_{n+1},2t_n]$, notice from~\eqref{xx} that 
$\|\glue(t)-u_0\|_{L^2_x}\leq a_n$. In particular, $\glue\big|_{t=0}=u_0$ defines an extension with
$\glue-\srad\in \CC^0_t L^2_x([0,\infty[\times \R^2)$.

The sixth point follows from the PDE for $\ugood$ with $\e= \e_n$ and $\e_{n+1}$, the choice of $\chi_n$ in \eqref{def, chi-n}, as well as the definition of $\mathscr{B}_{\e_n, \e_{n+1}}$ in \eqref{B, def, reynold} and the identity $$ \zz{n} - \zz{n+1} = {\rm div} \mathscr{B}_{\e_n, \e_{n+1}} + \na \pi_{\e_n,\e_{n+1}} $$ for some scalar field $\pi_{\e_n,\e_{n+1}}$; see Lemma~\ref{lemma: R-initial}.   

It remains to prove the last point. 
Fix $n \in \mathbb{N}$. In view of~\eqref{expression for R-heart} and $0 \leq \chi_n \leq 1$, we have
\begin{align*}
\left\|\rey\right\|_{L^\infty_t L^1_x\left(I_n \times \R^2\right)} &\leq \left\|\chi_n'\right\|_{L^\infty_t(I_n)} \left\|\mathscr{B}_{\e_n, \e_{n+1}}\right\|_{L^\infty_t L^1_x\left(I_n \times \R^2\right)} + \left\|\zz{n}-\zz{n+1}\right\|_{L^\infty_t L^2_x\left(I_n \times \R^2\right)}^2.
\end{align*}

\begin{itemize}
    \item 
For the first term, note by~\eqref{def, chi-n} that $\left\|\chi_n'\right\|_{L^\infty_t(I_n)}  \leq \frac{400}{t_n}$. In addition,
\begin{align*}
    \left\|\mathscr{B}_{\e_n, \e_{n+1}}\right\|_{L^\infty_t L^1_x\left(I_n \times \R^2\right)} &\leq \left\|\mathscr{R}^{\mathrm{initial}}_{\e_n,\e_{n+1}}\right\|_{L^1_x} + 2t_n \left\|\mathscr{DR}_{\e_n,\e_{n+1}}\right\|_{L^\infty_t L^1_x\left(I_n \times \R^2\right)} \\
    &\leq C\e_n + 2t_n \crey D_n. 
\end{align*}
Here, the first inequality follows from the definition of $\mathscr{B}_{\e_n, \e_{n+1}}$ in~\eqref{B, def, reynold} and $I_n = ]t_n,2t_n[$, and the second inequality holds by the bound for $\mathscr{R}^{\mathrm{initial}}_{\e_n,\e_{n+1}}$ in Lemma~\ref{lemma: R-initial}, as well as the estimate for  $\mathscr{DR}_{\e_n,\e_{n+1}}$ in~\eqref{xx} and $t_n<h_n/2$ in~\eqref{choice of tn}. The constant $C$ is the universal one in Lemma~\ref{lemma: R-initial}. Hence, as $\e_n \leq \crey t_n D_n$ by \eqref{choice of en}, we obtain that
\begin{align*}
\left\|\chi_n'\right\|_{L^\infty_t(I_n)} \left\|\mathscr{B}_{\e_n, \e_{n+1}}\right\|_{L^\infty_t L^1_x\left(I_n \times \R^2\right)} \leq (400C+800)\crey D_n.
\end{align*}

\item 
For the second term, since $t<2t_n<h$ for $t \in I_n$ by~\eqref{choice of tn}, it follows from~\eqref{xx} that 
\begin{align*}
    \left\|\zz{n}-\zz{n+1}\right\|_{L^\infty_t L^2_x\left(I_n \times \R^2\right)}^2 \leq (\crey)^2 D_n.
\end{align*}
\end{itemize}
We thus arrive at 
\begin{align*}
    \left\|\rey\right\|_{L^\infty_t L^1_x\left(I_n \times \R^2\right)} \leq \big[(400C+800)\crey  + (\crey)^2\big]D_n. 
\end{align*}
The proof is complete by choosing any $\crey \leq \sqrt{(200C+400)^2+1} -(200C+400)$.  
\end{proof}

\section{Proof of the main theorem; convex integration}\label{sec: proof}

In this section, we shall conclude the proof of Theorem~\ref{mainthm}.

 For this purpose, we start with the glued velocity $\glue$ and, through a sequence of stages indexed by $q$, insert the Lamb--Chaplygin dipole-based perturbations iteratively to decrease the $L^1_x$-norm of the Reynolds stress tensor in each epoch $I_n$; see Lemma~\ref{lem: key}. In the limit $q \nearrow \infty$, we arrive at the global classical weak solution $u$ to the 2D Euler equations with initial vorticity measure $\omega_0$. The attainability of the initial data is proved by taking $n \nearrow \infty$ due to the temporal glueing construction in \S\ref{sec: glue}.

\subsection{Iteration scheme} \label{subsec: iteration scheme}

Fix $n \in \mathbb{N}_+$. We choose once and for all
\[
 J_n=\left[\left(1+\frac{1}{99}\right)t_n, \left(2-\frac{1}{99}\right)t_n \right]\Subset ]t_n,2t_n[ =:I_n.
\]
The insertion of Lamb--Chaplygin dipoles will be carried out in each $J_n$, \textit{i.e.}, away from the endpoints of the epochs $I_n$. 
Observe that $R^{\heartsuit}$ vanishes on $[0,\infty[ \setminus\bigsqcup_n J_n$: In fact, by Proposition~\ref{propn: glue}~(6), $R^{\heartsuit}$ also vanishes on relative neighbourhoods in $J_n$ of both its endpoints.

We introduce the iteration scheme as follows:

\smallskip
\noindent
\underline{Initialisation.} The iteration begins at Stage 0: On $J_n$ we take
\begin{align}\label{stage 0}
 &(u,p,R) := \left(\glue, \pi^{\heartsuit}, R^{\heartsuit}\right) \qquad\text{at } q=0,\nonumber\\
 &\delta := \delta_{n,0} \equiv D_n,
\end{align}
and return to $(U_1,p_1,R_1)$ with $\e \equiv \delta_{n,1}$ in Lemma~\ref{lem: key}, where $p_1$ is a pressure field  associated with $(U_1, R_1)$. Recall~\eqref{Dn, def} for the choice of $D_n$.

\smallskip
\noindent
\underline{Inductive step.}  We inductively construct smooth triples $(U_q,p_q,R_q)$ for $t>0$ via applications of Lemma~\ref{lem: key}: Set for each $q \in \mathbb{N}$  that
\begin{equation}\label{delta,n,q: def}
    \delta_{n,q} := \frac{D_n}{4^q} = \frac{\kappa^2}{4^{n+q}}.
\end{equation}
At each stage~$q$, we apply Lemma~\ref{lem: key} with 
\begin{align}\label{inductive step}
    I=J_n,\qquad \delta=\delta_{n,q},\qquad\text{and } \e = \delta_{n,q+1}.
\end{align} 
If $R_q=0$ for some $q$, then we simply set $\left(U_{q+1}, p_{q+1}, R_{q+1}\right) := (U_q, p_q, R_q)$.

The triples $(U_q, p_q, R_q)$ thus constructed have the following properties. Throughout this paper, we refer to $(U_q, p_q, R_q)$ as a \emph{Euler--Reynolds triple} if they satisfy the following equations:
\begin{align}\label{ER triple}
\begin{cases}
    \partial_t U_q+\operatorname{div}(U_q\otimes U_q)+\nabla p_q
       =\operatorname{div}R_q,\\
    \operatorname{div}U_q=0 \qquad \text{ in } ]0,\infty[ \times \R^2.
    \end{cases}
\end{align}

\begin{lemma}\label{lem: induction} For  each $q \geq 0$, it holds that: 
\begin{enumerate}
 \item $(U_q, p_q, R_q) \equiv \left(\glue, \pi^\heartsuit, \rey\right)$ on $]0,\infty[ \setminus \bigsqcup_{n=1}^\infty J_n$, as well as on some relative open neighbourhoods of the endpoints of each $J_n$.  
 
 \item On each epoch, $U_q\in \CC^0_tL^3_x(J_n\times \R^2)$,
 $U_q-\srad\in \CC^0_tL^2_x(J_n\times \R^2)$, and $\operatorname{curl}U_q\in \CC^0_tL^1_x(J_n\times \R^2)$. In addition, $\curl\,U_q$ has a common compact spatial support for each $t\in J_n$.
 \item On each epoch, $R_q\in \CC^0_t L^1_x(J_n \times \sym)$, with
 $\sup_{t\in J_n}\|R_q(t,\cdot)\|_{L^1_x}\leq\delta_{n,q}$.
 \item The circulation, the first moment, and the second moment of  $\curl\,U_q$ equal those of $\omega_0$, the unregularised initial vorticity measure, at every $t>0$.
\end{enumerate}    
\end{lemma} 

This is a direct consequence of Lemma~\ref{lem: key}, so the proof is omitted.  As a caveat, the radius of the support of $\curl\, U_q$ is allowed to depend on $n$ and $q$.

\subsection{Bounds for limiting velocity and vorticity}

Denote the velocity increments  as
\begin{equation}\label{wnq, increment}
\begin{cases}
    w_{n,q}(t,x) :=  \Big(U_{q+1} (t,x) - U_q (t,x)\Big)\1_{J_n}(t),\\
    w_q:=\sum_{n=1}^\infty w_{n,q}.
    \end{cases}
\end{equation}
The preceding construction in \S\ref{subsec: iteration scheme} yields that
\begin{align}
 \sup_{t\in J_n}\|R_{q+1}(t,\cdot)\|_{L^1_x}
     &\leq\delta_{n,q+1},\label{eq:s4-stress-step}\\
 \|w_{n,q}\|_{L^\infty_tL^2_x}
 +\|\operatorname{curl}w_{n,q}\|_{L^\infty_tL^1_x}
     &\leq C_0\kappa\,2^{-n-q}
                         +\kappa^2 4^{-n-q-1}.\label{eq:s4-step}
\end{align}
Here, $C_0$ is the universal constant in Lemma~\ref{lem: key}. In particular, these bounds do not depend on the radii of the supports of $U_q$, nor on the lengths or positions of $J_n$.  In particular, \[
\|w_q(t,\cdot)\|_{L^2_x}+\|\operatorname{curl}w_q(t,\cdot)\|_{L^1_x}
 \leq C_0\kappa\,2^{-1-q}+\kappa^2 4^{-2-q}\qquad\text{for any } t>0 
\]
and for every fixed integer $Q\geq 0$, 
\begin{align}
 \sum_{q=Q}^{\infty}
 \bigl(\|w_{n,q}\|_{L^\infty_tL^2_x}
       +\|\operatorname{curl}w_{n,q}\|_{L^\infty_tL^1_x}\bigr)
 &\leq 2C_0\kappa\,2^{-n-Q}
             +\tfrac13\kappa^2 4^{-n-Q}\notag\\
 &\leq \underbrace{\left(2C_0+\frac{1}{3}\kappa\right)}_{=:\, C_1}\,2^{-n-Q}.\label{sumestimate}
\end{align}
In the above, $w_{n,q}$ and $w_q$ as well as their curl are continuous in the indicated norms across the endpoints of time intervals $J_n$.

It is crucial to observe that the right-hand side of \eqref{eq:s4-step} tends to zero uniformly in
$t\in J_n$ for fixed $q$ as $n\to\infty$. Thus,  $U_q-\srad$ converges uniformly in $L^2_x$ on each
compact time interval in $[0,\infty[$, and that
$\omega_q=\operatorname{curl}U_q$ converges uniformly in $L^1_x$ for $t \in ]0,\infty[$. In particular,
\begin{align}\label{def of u, omega as limits}
&\lim_{q \to \infty} \left(U_q-\srad\right) =: u - \srad \in C^0_t L^2_x \left([0,\infty[ \times \R^2\right),\nonumber\\
&\lim_{q \to \infty} \omega_q = \lim_{q \to \infty} \curl\, U_q =: \omega \in C^0_t L^1_x \left(]0,\infty[ \times \R^2\right).
\end{align}

Fix $N \in \mathbb{N}$. For $0<t\leq 2t_N$, using \eqref{sumestimate} with $Q=0$, the trace estimate for $\glue$ in Proposition~\ref{propn: glue}~(3), as well as the monotonicity of $a_n$ as in~\eqref{xx}, we obtain that
\begin{equation}\label{eq:s4-initial-trace}
 \sup_{0<t\leq2t_N}\|U_q(t,\cdot)-u_0(\cdot)\|_{L^2_x} \leq a_N+C_1\kappa 2^{-N}.
\end{equation}
This bound is uniform in $q \in \mathbb{N}$, and it remains valid for $u$ defined in~\eqref{def of u, omega as limits}. 

Moreover, in view of the bound for the $L^1$-norm of $\omega^\heartsuit$ in Proposition~\ref{propn: glue}, Point~(5) and the estimate~\eqref{sumestimate}, it holds for each $N \in \mathbb{N}$ that
\begin{align}
 \sup_{q\geq 0}\sup_{0<t\leq 2t_N}\|\omega_q(t,\cdot)\|_{L^1_x} \leq \left|\omega_0\right|(\R^2)+C_\G\e_N^2+C_1\kappa 2^{-N},\label{eq:s4-near-zero-tv}
\end{align}
By further taking $N=1$ in~\eqref{eq:s4-near-zero-tv}, we arrive at
\begin{equation}
 \sup_{q\in \mathbb{N}}\sup_{t>0}\|\omega_q(t,\cdot)\|_{L^1_x}
       \leq\left|\omega_0\right|(\R^2)+C_\G\e_1^2+\frac{C_1}{2}\kappa.\label{eq:s4-global-tv}
\end{equation}
Both bounds~\eqref{eq:s4-near-zero-tv} and \eqref{eq:s4-global-tv} remain valid for $\omega(t,\cdot)$ in place of $\omega_q(t,\cdot)$ at each $t>0$. Thus, by virtue of the choice for $\e_1$ in~\eqref{choice of en}, it holds that $C_0\e_1^2\leq\kappa$. Thus, given any fixed $\e_\star>0$, we pick $\kappa>0$ such that 
\begin{equation}\label{kappa, choice}
\kappa\leq\min\left\{1,\,\frac{\varepsilon_\star}{{C_0+7/6}} \right\}
\end{equation}
to obtain the desired bound in Theorem~\ref{mainthm}:
\begin{align}
\sup_{t>0}\|\omega(t,\cdot)\|_{L^1_x}\le |\omega_0|(\R^2)+\e_\star.
\end{align}

\begin{remark}
Let us comment that for $\omega_0\ne0$, one may  ensure  $\omega(t,\cdot)\neq 0$ at every $t>0$ by tuning the parameters $\e_1,t_1,\kappa$.  Choose
$\varphi_*\in \CC_c^\infty(\R^2)$ with $|\varphi_*|\leq1$ and
$\langle\omega_0,\varphi_*\rangle>\frac{3}{4}|\omega_0|(\R^2)$. By Lemma~\ref{lemma:joint trace},
there exists $\eta_*>0$ such that
\[
 \left\langle \gsol(t,\cdot),\varphi_*\right\rangle>\frac{|\omega_0|(\R^2)}{2}
 \qquad\text{if }0<\e<\eta_*,\quad 0\leq t<\eta_*.
\]
Impose also $\e_1<\eta_*$, $2t_1<\eta_*$, and
$C_1\kappa\leq |\omega_0|(\R^2)/4$. Then  \eqref{sumestimate}
and $\|\varphi_*\|_\infty\leq1$ lead to
\[
 \langle\omega_q(t,\cdot),\varphi_*\rangle\geq \frac{|\omega_0|(\R^2)}{4}
 \qquad \text{for each $q \in \mathbb{N}$ and $t \in ]0,2t_1]$}.
\]
The same holds for $\omega(t,\cdot)$.

Meanwhile, we check that the limiting velocity $u(t,\cdot)$ cannot support on a Lebesgue null set in $\R^2$ for any $t>0$. Indeed, let $\B$ be any fixed ball containing
$\operatorname{supp}\varphi_*$. Via integration by parts we see that
\[
 \|U_q(t,\cdot)\|_{L^2(\B)},\ \|u(t,\cdot)\|_{L^2(\B)}
 \geq \frac{|\omega_0|(\R^2)}{4\left\|\nabla\varphi_*\right\|_{L^2_x}}.
\]
Suppose that  $U_q(t,\cdot)$ vanishes on $\B\setminus E_q$ for
sets $E_q$ with $\leb^2(E_q)\to0$ for some $t>0$. Then
\[
 \|u(t,\cdot)\|_{L^2_x(\B)}
 \leq \|u(t,\cdot)-U_q(t,\cdot)\|_{L^2_x(\B)}
      +\|u(t,\cdot )\1_{\B\cap E_q}\|_{L^2_x}\longrightarrow 0.
\]
Contradiction.

Also note that for any  $t\geq2t_1$, all velocities coincide with $\zz{1}$, which is nonzero at $t=2t_1$ by the preceding paragraph. 

\end{remark}

\subsection{Classical weak solution}

We check that the limiting velocity field $u$ obtained in~\eqref{def of u, omega as limits} is a classical weak solution. The argument is divided into four steps below. 

\smallskip
\noindent
{\bf Step~1}. Clearly,  $\operatorname{div}u(t,\cdot)=0$ for all $t\geq0$ by the convergence of $u(t,\cdot)$ in $\mathcal{D}'$.

\smallskip
\noindent
{\bf Step~2}. We next show that for any $\phi \in \CC^\infty_c\left([0,\infty [\times \R^2\right)$ with ${\rm div}\phi=0$,
\begin{equation}\label{eq:s4-weak-velocity}
 \int_0^\infty\!\!\int_{\R^2}
     \bigl(u\cdot\partial_t\phi+(u\otimes u):\nabla\phi\bigr)\,\dd x\,\dd t
       +\int_{\R^2}u_0(x)\cdot\phi(0,x)\,\dd x=0.
\end{equation}
 That is, $u$ satisfies the 2D Euler equations in the sense of distributions.

Indeed, by virtue of the bound~\eqref{eq:s4-stress-step} and that $R_q(t,\cdot)\equiv 0$ for $t \in ]0,\infty[ \setminus \bigsqcup_{n=1}^\infty J_n$, we have \begin{equation}\label{eq:s4-global-stress}
 \sup_{t>0}\|R_q(t,\cdot)\|_{L^1_x}
 \leq \sup_{n\geq1}\delta_{n,q}=\kappa^2 4^{-1-q}.
\end{equation}
 Let $\phi$ be a compactly
supported smooth solenoidal test field on $[0,\infty)\times\R^2$, with ${\rm supp}(\phi) \Subset [0,T [\times \B$ for some bounded ball $\B$. for any $0<s<T$, an integration by parts applied to the Euler--Reynolds triple $(U, p_q,R_q)$ in~\eqref{ER triple} yields that
\[
 \int_s^T\!\!\int_{\R^2}
 \bigl(U_q\cdot\partial_t\phi+
       (U_q\otimes U_q-R_q):\nabla\phi\bigr)\,\dd x\,\dd t
       +\int_{\R^2}U_q(s,x)\cdot\phi(s,x)\,\dd x=0.
\]

\smallskip
\noindent
\underline{Send $s \searrow 0$:} For the first term on the left-hand side, the integral on $[0,s]$ vanishes as $s\searrow 0$ by the dominated convergence theorem. Indeed,  the $L^1$-norm of $R_q$ is bounded by~\eqref{eq:s4-global-stress}, while the $L^2$-norm of $U_q$ is bounded by $\left\|\glue-\srad\right\|_{L^2_{t,x}([0,T]\times\B)}$ by Proposition~\ref{propn: glue}, as well as the summable velocity increments as in~\eqref{sumestimate}.  The second term
converges to $\int_{\R^2} u_0(x)\cdot\phi(0,x)\,\dd x$ by \eqref{eq:s4-initial-trace}
and the smoothness of $\phi$.

\smallskip
\noindent
\underline{Send $q \nearrow \infty$:} 
For the quadratic term, we bound
\begin{align*}
&\|U_q(t,\cdot)\otimes U_q(t,\cdot)-u(t,\cdot)\otimes u(t,\cdot)\|_{L^1_x(\B)} \\
&\qquad\leq \bigl(\|U_q(t,\cdot)\|_{L^2_x(\B)}+\|u(t,\cdot)\|_{L^2_x(\B)}\bigr)
 \|U_q(t,\cdot)-u(t,\cdot)\|_{L^2_x(\B)}.
\end{align*}
The right hand side tends uniformly to zero. Moreover, by~\eqref{eq:s4-global-stress},
\begin{align*}
\left|\int_0^T\int_{\R^2} R_q:\nabla\phi\,\dd x\,\dd t \right| \leq T\|\nabla\phi\|_{L^\infty_{t,x}}\kappa^2 4^{-1-q} \longrightarrow 0 \quad\text{as } q \nearrow \infty.
\end{align*}
This proves~\eqref{eq:s4-weak-velocity}.


\smallskip
\noindent
{\bf Step~3}. We show that $u$ is continuous in $t$ and $u(t,\cdot)$ lies in $H^{-L}_\loc$ for each $t$. This is part of the definition for the classical weak solution \textit{\`{a} la} DiPerna--Majda~\cite{dm1, dm2, dm3}. 

\begin{lemma}\label{lemma:dm-time-regularity}
Let $\srad$ be the fixed smooth, bounded, radial stationary Euler velocity field as before. Let $u$ be a distributional solution to the 2D Euler equations, \emph{i.e.}, \eqref{eq:s4-weak-velocity}. If $
 v:=u-\srad \in \CC^0_{t}L^2_x\left([0,\infty)\times\R^2\right)$, then $v \in \CC^1_t H^{-L}_x \left([0,\infty)\times\R^2\right)$ for every $L>2$. As a consequence, $ u\in\operatorname{Lip}\left([0,T];H^{-L}_{\mathrm{loc}}(\R^2)\right)$ for any $0<T<\infty$.  \end{lemma}

The constant depends on $\sup_{0\le t\le T}\|v(t,\cdot)\|_{L^2_x}$. On the whole plane, one must exclude a harmonic pressure term before applying the Leray projection. Also, $|\srad| \simeq |x|^{-1}$ at spatial infinity, so $\srad\notin H^{-L}(\R^2)$. Hence, we may only conclude $u(t,\cdot) \in H^{-L}_\loc$ from $v(t,\cdot) \in H^{-L}$.

\begin{proof}[Proof of Lemma~\ref{lemma:dm-time-regularity}]

Let $\leray$ be the Leray projection on $\R^2$; \emph{i.e.}, the $L^2$-orthogonal projection on divergence-free parts of vector fields with Fourier symbol $\xi \mapsto {\mathrm{Id}} - {(\xi \otimes \xi)} / {|\xi|^2}$. For any tensor field $M \in L^1(\R^2; \R^2 \otimes \R^2)$, by computations on the Fourier side we have that
\begin{align*}
\left\| \lerayy M\right\|^2_{H^{-L}_x} \leq \frac{\|M\|^2_{L^1_x}}{4\pi^2}  \int_{\R^2}\frac{|\xi|^2}{(1+|\xi|^2)^L}\,\dd \xi =\frac{\|M\|_{L^1_x}^2}{2\pi} \int_0^\infty\frac{r^3}{(1+r^2)^L}\,\dd r,
\end{align*}
which is finite whenever $L>2$. On the other hand, 
For $M\in L^2(\R^2; \R^2\otimes \R^2)$ and any $L>1$, as $|\xi|^2/(1+|\xi|^2)^L\le1$, Plancherel leads to $\|\lerayy M\|_{H^{-L}_x}\le\|M\|_{L^2_x}$. Thus by interpolation
\begin{equation}
\text{$\lerayy$\quad is $(L^1\cap L^2)(\R^2) \longrightarrow H^{-L}(\R^2)$-bounded for any $L>2$}.
\end{equation}

Fix any $T>0$ and consider $v:=u-\srad$. Then $v \otimes v \in \CC^0_t L^1_x\left([0,T] \times \R^2\right)$ and $\srad \otimes v + v\otimes \srad \in \CC^0_t L^2_x\left([0,T] \times \R^2\right)$, since $\srad \in L^\infty_x$ and 
\[
 \|v(t,\cdot)\otimes v(t,\cdot)-v(s,\cdot)\otimes v(s,\cdot)\|_{L^1_x}
 \le\bigl(\|v(t,\cdot)\|_2+\|v(s,\cdot)\|_{L^2_x}\bigr)\|v(t,\cdot)-v(s,\cdot)\|_{L^2_x}.
\]
Thus, the vector field
\begin{align*}
\varpi := -\lerayy \underbrace{\left(v \otimes v + \srad \otimes v + v\otimes \srad\right)}_{=: \mathcal{T}(t,x)}
\end{align*}
belongs to $\CC^0_t H^{-L}_x\left([0,T] \times \R^2\right)$. We show that it coincides with $\p_tv$.

For this purpose, as $\sigma$ is a stationary Euler field and $u$ is a distributional solution to the 2D Euler equations, for every
$\phi\in \CC_c^\infty(]0,T[\times\R^2;\R^2)$ with
$\operatorname{div}\phi=0$, it holds that
\begin{equation}\label{yy}
\int_0^T\!\int_{\R^2}
       \bigl(v\cdot\partial_t\phi+\mathcal{T}:\nabla\phi\bigr)\,\dd x\,\dd t=0.
\end{equation}
This can be established without assuming pointwise defined $\p_tv$. Take any $\eta\in \CC_c^\infty(]0,T[)$. Set
\begin{equation}\label{eq:time-averaged-w}
 w_\eta
    =-\int_0^T\eta'(t)v(t)\,\dd t
      -\int_0^T\eta(t)\varpi(t)\,\dd t
       \ \in H^{-L}(\R^2).
\end{equation}
Both $v$ and $\varpi$ are divergence-free, so
$\operatorname{div}w_\eta=0$. By testing against
$\phi(t,x)=\eta(t)\nabla^\perp\psi(x)$ with
$\psi\in \CC_c^\infty(\R^2)$, we also have 
$\operatorname{curl}w_\eta=0$.  Hence, $ \Delta w_\eta
   =\nabla\operatorname{div}w_\eta
       +\nabla^\perp\operatorname{curl}w_\eta=0$.  But there is no nonzero harmonic distribution in $H^{-L}(\R^2)$: indeed, $\widehat{w_\eta}$ vanishes on $\R^2_\xi \setminus \{0\}$, and $(1+|\xi|^2)^{-L/2}\widehat{w_\eta}(\xi)\in L^2_\xi$ implies that $\widehat{w}_\eta\in L^2_\loc$, hence cannot be supported on $\{0\}$. Thus $w_\eta=0$.

Then, since $\mathbb P\phi=\phi$, we have that
\[
 \langle \varpi,\phi\rangle=\int_0^T\!\int_{\R^2}\mathcal{T}:\nabla\phi\,\dd x\,\dd t.
\]
Thus $w=\partial_t v-\varpi$ annihilates every compact smooth solenoidal test vector field.

Now we have 
\begin{equation}\label{FTC}
 v(t)-v(s)=\int_s^t \varpi(\tau)\,\dd\tau\qquad\text{ for any }0\le s,t\le T, 
\end{equation}
understood as a Bochner integral taking values in $H^{-L}_x$; $L>2$. Note that $s=0$ can be attained. It follows that \[
 \lim_{t\searrow 0}
   \left\|\frac{v(t)-v(0)}t-\varpi(0)\right\|_{H^{-L}_x}=0,
\]
and the derivative is continuous up to $t=0^+$. Thus, by~\eqref{FTC} we have $v \in \CC^1_t H^{-L}_x \left([0,\infty)\times\R^2\right)$. Moreover, by construction $u(t) - u(s)=  v(t)-v(s)$, so it follows that $ u\in\operatorname{Lip}\left([0,T];H^{-L}_{\mathrm{loc}}(\R^2)\right)$ for any $0<T<\infty$. This completes the proof.

\smallskip
\noindent
{\bf Step~4}. For completeness, we show that a distributional pressure can be recovered.  Set $F:=\partial_tu+\operatorname{div}(u\otimes u)$
on $]0,\infty[\times\R^2$. By \eqref{eq:s4-weak-velocity}, $F$ annihilates any compactly supported, smooth, spatially solenoidal test vector fields.

Fix any $\theta\in \CC_c^\infty(\R)$ with $\int_\R \theta(x)\,\dd x=1$.
For any test function $g \in \CC^\infty_c(]0,\infty[\times \R^2)$ with zero spatial integral at each time,
let $G(t,x_2):=\int_{\R}g(t,s,x_2)\,\dd s$.

Define
\begin{align*}
\mathbf{Bog}[g](t,x):=\begin{bmatrix}
\int_{-\infty}^{x_1}
           \bigl(g(t,s,x_2)-\theta(s)G(t,x_2)\bigr)\,\dd s\\
        \theta(x_1)\int_{-\infty}^{x_2}G(t,s)\,\dd s.
\end{bmatrix}
\end{align*}
This Bogovskii map is continuous on the space of test functions, and
$\operatorname{div}  \mathbf{Bog}[g]=g$. Then, define $p \in \mathcal{D}'\left([0,\infty[\times \R^2\right)$ by duality: 
\begin{align*}
&\langle p,g\rangle=\langle F,\mathbf{Bog}[g_0]\rangle,\\
&g_0(t,x_1,x_2):=g(t,x_1,x_2) - \theta(x_1)\theta(x_2)\int_{\R^2}g(t,y)\,\dd y.
\end{align*}
Then, for an arbitrary test vector field $\Psi \in \CC^\infty_c([0,\infty[\times \R^2; \R^2)$, its divergence $\operatorname{div}\Psi$
has zero spatial integral and
$\Psi-\mathcal {\mathbf{Bog}}[\operatorname{div}\Psi]$ is solenoidal. Hence
$\langle p,\operatorname{div}\Psi\rangle=\langle F,\Psi\rangle$.
It follows that $\partial_tu+\operatorname{div}(u\otimes u)+\nabla p=0$ in the sense of distributions.  \end{proof}

\subsection{Attainability of the initial vorticity}

To conclude the proof of our Main Theorem~\ref{mainthm}, it remains to establish that
\begin{equation}\label{convergence of omega}
\omega(t,\cdot) \weakstar \omega_0\qquad\text{and}\qquad \|\omega(t,\cdot)\|_{L^1_x} \longrightarrow |\omega_0|(\R^2)\qquad\text{as } t \searrow 0.
\end{equation} 

For this purpose, let us first verify that the limiting vector field $u$ and the limiting scalar field $\omega$ given by~\eqref{def of u, omega as limits} satisfy that
\begin{equation}\label{u, omega, final}
u(t,\cdot) = \kbs\star\omega(t,\cdot)\qquad\text{for every } t>0.
\end{equation}

Indeed, fix any $t>0$. The smooth, finite-stage vorticity $\omega_q(t,\cdot)$ is compactly supported and has circulation $\Gamma \equiv \int_{\R^2}\,\dd\omega_0$. Also recall from~\eqref{b, def} that $b \in \CC^\infty_c(\R^2)$ is a radial scalar field satisfying $\int_{\R^2} b(x)\,\dd x = \G$. It follows that $\nu(x)=\omega_q(t,x)-b(x)$ has compact support, is smooth, and has zero
integral, and hence $\kbs\star\nu=\mathcal{O}(|x|^{-2})$ at spatial infinity. In particular, $\kbs\star\nu \in L^2(\R^2)$.

The difference \[
 U_q(t,\cdot)-\kbs\star\omega_q(t,\cdot)
   = [U_q(t,\cdot)-\srad]-\kbs\star(\omega_q(t,\cdot)-b)
\]
then belongs to $L^2$ and is harmonic (as it has zero curl and zero divergence), and hence it must be zero. Thus, $U_q \equiv \kbs \star \omega_q$ for each $q \in \mathbb{N}$.

Notice that for any bounded measurable set $A\subset \R^2$,
\[
 \sup_{y\in\R^2}\int_A |\kbs(x-y)|\,\dd x
 \leq \int_{|z|\leq 1}|\kbs(z)|\,\dd z + \frac{\leb^2(A)}{2\pi}<\infty.\] By Tonelli's theorem, $\|\kbs\star f\|_{L^1_x(A)}\leq C\|f\|_{L^1_x}$ for each $f\in L^1(\R^2)$, with $C$ depending on $\leb^2(A)$. Hence, $K*\omega_q(t,\cdot)\to K*\omega(t,\cdot)$ in $L^1_\loc$ and $U_q (t,\cdot) \to u(t,\cdot)$ in $L^2_\loc$ for each $t>0$. This allows us to send $q \to \infty$ in $U_q \equiv \kbs \star \omega_q$ to deduce~\eqref{u, omega, final}. Strong $L^1$-convergence of the curl also ensures that $\int_{\R^2}\omega(t,x)\,\dd x=\Gamma$ for each $t>0$.

Finally, for every spatial test function
$\varphi\in \CC_c^\infty(\R^2)$, it holds that
\begin{align*}
 \big|\langle\omega(t,\cdot)-\omega_0,\varphi\rangle\big|
 &=\left|\int_{\R^2}\big(u(t,x)-u_0(x)\big)\cdot\nabla^\perp\varphi(x)\,\dd x\right|\\
 & \leq \|u(t,\cdot)-u_0\|_{L^2_x}\,\|\nabla\varphi\|_{L^2_x}\longrightarrow 0.
 \end{align*}
Here we make crucial use of \eqref{u, omega, final} established above. Then, thanks to the uniform $L^1_x$-bound for $\omega(t,\cdot)$  and a density argument, we deduce that $\omega(t,\cdot) \weakstar \omega_0$ as measures. By lower semicontinuity of the weak convergence of Radon measures,  $|\omega_0|(\R^2) \leq \liminf_{t \searrow 0} \|\omega(t,\cdot)\|_{L^1_x}$. Furthermore, by~\eqref{eq:s4-near-zero-tv} we also have $\limsup_{t\searrow 0}\|\omega(t,\cdot)\|_{L^1_x}\leq |\omega_0|(\R^2)$. This proves~\eqref{convergence of omega}.

The proof of  Theorem~\ref{mainthm} is now complete.

\section{Concluding remarks}\label{sec: concluding remarks}

\subsection{The weak vorticity formulation \emph{\`{a} la} Delort and Schochet}
\label{subsec:delort-schochet}

The weak vorticity formulation (see~\eqref{schochet, kernel} in the introduction) goes back to Delort~\cite{Delort}. Here we follow the symmetric kernel formulation in Schochet~\cite{Schochet1995}; see Lopes Filho--Nussenzveig Lopes--Schochet~\cite[Definition~2.1
and equation~(2.9)]{LNS2007} for the formulation with initial data.
With our convention \(\kbs(z)=z^\perp/(2\pi|z|^2)\), for
\(\varphi\in\CC^\infty_c(\R^2)\) we define the symmetric kernel:
\begin{equation}\label{eq:ds-kernel}
 H_\varphi(x,y)=
 \begin{cases}
 \dfrac12\kbs(x-y)\cdot\bigl(\nabla\varphi(x)-\nabla\varphi(y)\bigr),
          & x\ne y,\\[3pt]
 0,       & x=y.
 \end{cases}
\end{equation}
The mean value
estimate for \(\nabla\varphi\) gives the pointwise bound
\begin{equation}\label{eq:ds-kernel-bound}
 \sup_{x,y\in\R^2}|H_\varphi(x,y)|
    \le \frac{1}{4\pi}\|D_x^2\varphi\|_{L^\infty_x}.
\end{equation}

The kernel is continuous off the diagonal. For a finite signed measure \(\nu\) without atoms, the diagonal is null for \(|\nu|\otimes|\nu|\), so its assigned value in \eqref{eq:ds-kernel}
is immaterial.

For a weak-\(*\) measurable family \(s\mapsto\nu_s\) of finite
signed vorticity measures, atomless for almost every \(s\) and
essentially bounded in total variation on \(]0,T[\), the
Delort--Schochet formulation is \eqref{eq:ds-weak} below, with
\(\dd\nu_s(x)\) in place of \(\omega(s,x)\,\dd x\) and
\(\dd\nu_s(x)\,\dd\nu_s(y)\) in place of
\(\omega(s,x)\omega(s,y)\,\dd x\,\dd y\). The initial datum is
paired with \(\Phi(0,\cdot)\). In our case the positive time
vorticity measures have densities \(\omega(s,\cdot)\).

\begin{proposition}
\label{propn:delort-schochet}
Let \(u\) and \(\omega=\curl u\) be the solution and vorticity
obtained in the proof of Theorem~\ref{mainthm}. For every \(T>0\) and every
scalar \(\Phi\in\CC^\infty_c([0,T[\times\R^2)\),
\begin{equation}\label{eq:ds-weak}
\begin{aligned}
 &\int_{\R^2}\Phi(0,x)\,\dd\omega_0(x)
   +\int_0^T\!\int_{\R^2}
           \omega(s,x)\,\p_s\Phi(s,x)\,\dd x\,\dd s\\
 &\quad+\int_0^T\!\iint_{\R^2\times\R^2}
    H_{\Phi(s,\cdot)}(x,y)\,
    \omega(s,x)\omega(s,y)\,\dd x\,\dd y\,\dd s
  =0.
\end{aligned}
\end{equation}
Moreover, for every \(t>0\) and \(\varphi\in\CC^\infty_c(\R^2)\),
\begin{equation}\label{eq:ds-fixed-time}
\begin{aligned}
 &\int_{\R^2}\omega(t,x)\varphi(x)\,\dd x
       -\int_{\R^2}\varphi(x)\,\dd\omega_0(x)\\
 &\qquad=\int_0^t\!\iint_{\R^2\times\R^2}
     H_\varphi(x,y)\,\omega(s,x)\omega(s,y)
                         \,\dd x\,\dd y\,\dd s .
\end{aligned}
\end{equation}
\end{proposition}

\begin{proof}[Proof of Proposition~\ref{propn:delort-schochet}]
Fix \(T>0\). By \eqref{def of u, omega as limits} and
\eqref{eq:s4-global-tv}, we have
\[
 u\in L^\infty(0,T;L^2_{\mathrm{loc}}(\R^2)),
 \qquad
 \omega\in L^\infty(0,T;\M(\R^2)).
\]
Since \(\omega=\curl u\), we also have
\(\omega\in L^\infty(0,T;H^{-1}_{\mathrm{loc}}(\R^2))\).
The Biot--Savart identity is \eqref{u, omega, final}, and the
initial datum satisfies \(u_0=\kbs\star\omega_0\).
Schochet's equivalence lemma \cite[Lemma~2.1]{Schochet1995}
therefore yields the equivalence between the weak velocity equation~\eqref{eq:s4-weak-velocity} and \eqref{eq:ds-weak}; see also
\cite[\S~2.1]{LNS2007} for the formulation for the Cauchy problem.

To obtain~\eqref{eq:ds-fixed-time}, let us fix
\(\varphi\in\CC^\infty_c(\R^2)\) and use the test function 
\(\Phi(s,x)=\vartheta(s)\varphi(x)\) in \eqref{eq:ds-weak}.
It follows that the distributional time derivative of
\(\int_{\R^2}\omega(s,x)\varphi(x)\,\dd x\) is
\[
 \iint_{\R^2\times\R^2}
 H_\varphi(x,y)\,\omega(s,x)\omega(s,y)\,\dd x\,\dd y .
\]
This is bounded by \eqref{eq:ds-kernel-bound} and
\eqref{eq:s4-global-tv}. The pairing is continuous for \(s>0\)
by \eqref{def of u, omega as limits}, and its limit at \(s=0\)
is \(\int_{\R^2}\varphi(x)\,\dd\omega_0(x)\) by
\eqref{convergence of omega}. Integrating the derivative shows~\eqref{eq:ds-fixed-time} at every \(t>0\).   \end{proof}

In fact, Proposition~\ref{propn:delort-schochet} can be proved without invoking \cite[Lemma~2.1]{Schochet1995}.

\begin{proof}[Alternative proof of Proposition~\ref{propn:delort-schochet} via direct computation]
Fix \(0<a<t\) and \(\Phi\in\CC^\infty_c([0,\infty[\times\R^2)\).

We first work on \([a,t]\) where the finite-stage triples are smooth;
see \S\ref{subsec: iteration scheme} and
Lemma~\ref{lem: induction}. Let \(\omega_q=\curl U_q\). By
incompressibility in \eqref{ER triple},
\(\operatorname{div}(U_q\otimes U_q)=(U_q\cdot\nabla_x)U_q\). The identity
\[
 (U_q\cdot\nabla_x)U_q
   =\nabla_x\!\left(\tfrac12|U_q|^2\right)+\omega_q U_q^\perp
\]
follows by expanding componentwise. As
\(\curl(\omega_q U_q^\perp)=\operatorname{div}_x(\omega_q U_q)\),
we take curl in \eqref{ER triple} to get
\begin{equation}\label{eq:ds-staged-curl}
 \p_s\omega_q+\operatorname{div}(\omega_q U_q)
   =\curl\operatorname{div} R_q .
\end{equation}

The proof of \eqref{u, omega, final} establishes the
identity \(U_q(s,\cdot)=\kbs\star\omega_q(s,\cdot)\) for every \(s>0\) and $q \in \mathbb{N}$.
At each such time, \(\omega_q(s,\cdot)\) is smooth and compactly
supported. In particular, splitting the integral against \(|\kbs|\)
into \(|x-y|\le1\) and \(|x-y|>1\) yields that 
\[
 \int_{\R^2}|\kbs(x-y)|\,|\omega_q(s,y)|\,\dd y
 \le \|\omega_q(s,\cdot)\|_{L^\infty_x}
          +\frac1{2\pi}\|\omega_q(s,\cdot)\|_{L^1_x}
\]
uniformly in $x$. Thus Fubini's theorem and the interchange of \(x\) and \(y\) are justified in the transport term. As \(\kbs\) is odd,
\begin{align*}
 &\int_{\R^2} \omega_q(s,x)U_q(s,x)\cdot\nabla_x\Phi(s,x)\,\dd x\\
 &\quad=\iint_{\R^2\times\R^2} \kbs(x-y)\cdot\nabla_x\Phi(s,x)\,
           \omega_q(s,x)\omega_q(s,y)\,\dd x\,\dd y\\
 &\quad=\iint_{\R^2\times\R^2} H_{\Phi(s,\cdot)}(x,y)\,
           \omega_q(s,x)\omega_q(s,y)\,\dd x\,\dd y .
\end{align*}

Multiplying \eqref{eq:ds-staged-curl} by \(\Phi\), using the Leibniz rule, and integrating on \([a,t]\), we obtain that
\begin{equation}\label{eq:ds-finite-stage}
\begin{aligned}
 &\int_{\R^2}\omega_q(t,x)\Phi(t,x)\,\dd x
       -\int_{\R^2}\omega_q(a,x)\Phi(a,x)\,\dd x\\
 &=\int_a^t\!\int_{\R^2}
          \omega_q(s,x)\,\p_s\Phi(s,x)\,\dd x\,\dd s\\
 &\quad+\int_a^t\!\iint_{\R^2\times\R^2}
       H_{\Phi(s,\cdot)}(x,y)\,
       \omega_q(s,x)\omega_q(s,y)\,\dd x\,\dd y\,\dd s\\
 &\quad+\int_a^t\!\int_{\R^2}
       R_q(s,x):\nabla_x(\nabla_x^\perp\Phi)(s,x)\,\dd x\,\dd s .
\end{aligned}
\end{equation}

We now let \(q\to\infty\) with \(a>0\) fixed. The endpoint terms
and the term involving \(\p_s\Phi\) converge by the uniform
\(L^1_x\)-convergence of \(\omega_q\) on \([a,t]\); see~\eqref{def of u, omega as limits}. By~\eqref{eq:s4-global-stress}, the last term of \eqref{eq:ds-finite-stage}
is controlled by
\[
 (t-a)\sup_{a\le s\le t}\|\nabla_x(\nabla_x^\perp\Phi)(s,\cdot)\|_{L^\infty_x}
                         \,\kappa^2 4^{-1-q},
\]
which tends to zero. For the quadratic term, add and subtract
\(\omega(s,x)\omega_q(s,y)\) in the difference of the two vorticity
products. Using \eqref{eq:ds-kernel-bound} and integrating each
product of absolute values in \(x,y\) yields that, at each \(s\in[a,t]\),
\begin{align*}
 &\left|
  \iint_{\R^2\times\R^2}H_{\Phi(s,\cdot)}(x,y)
     \bigl(\omega_q(s,x)\omega_q(s,y)-\omega(s,x)\omega(s,y)\bigr)
                                  \,\dd x\,\dd y\right|\\
 &\qquad\le
 \frac{\|D_x^2\Phi(s,\cdot)\|_{L^\infty_x}}{4\pi}\,
  \bigl(\|\omega_q(s,\cdot)\|_{L^1_x}
        +\|\omega(s,\cdot)\|_{L^1_x}\bigr)
          \|\omega_q(s,\cdot)-\omega(s,\cdot)\|_{L^1_x}.
\end{align*}
The total variation bound \eqref{eq:s4-global-tv} and the same
\(L^1_x\)-convergence show that this term tends to zero uniformly on
\([a,t]\). Therefore,  \eqref{eq:ds-finite-stage} holds with
\(\omega_q\) replaced by \(\omega\) and with the stress term omitted. The two limiting integrands are continuous at positive
times, thanks to 
\eqref{def of u, omega as limits} and \eqref{eq:ds-kernel-bound}.

It remains to let \(a\searrow0\). By \eqref{convergence of omega},
\(\omega(a,\cdot)\,\dd x\weakstar\omega_0\). We also have
\begin{align*}
 &\left|\int_{\R^2}\omega(a,x)\Phi(a,x)\,\dd x
                     -\int_{\R^2}\Phi(0,x)\,\dd\omega_0(x)\right|\\
 &\quad\le \|\omega(a,\cdot)\|_{L^1_x}\,
              \|\Phi(a,\cdot)-\Phi(0,\cdot)\|_{L^\infty_x}
     +\big|\big\langle\omega(a,\cdot)\,\dd x-\omega_0,
                                      \Phi(0,\cdot)\big\rangle\big|.
\end{align*}
The first term tends to zero by \eqref{eq:s4-global-tv} and
smoothness of \(\Phi\); the second tends to zero by the weak-star convergence against the fixed test function \(\Phi(0,\cdot)\). The contributions of \((0,a)\) to the time derivative and
quadratic terms are bounded, respectively, by
\[
 a\sup_{0\le s\le t}\|\p_s\Phi(s,\cdot)\|_{L^\infty_x}
       \sup_{\tau>0}\|\omega(\tau,\cdot)\|_{L^1_x}
\]
and
\[
 \frac{a}{4\pi}\sup_{0\le s\le t}\|D_x^2\Phi(s,\cdot)\|_{L^\infty_x}
       \left(\sup_{\tau>0}\|\omega(\tau,\cdot)\|_{L^1_x}\right)^2 .
\]
Both bounds tend to zero. This proves the limiting identity
between \(0\) and \(t\).

To obtain \eqref{eq:ds-weak}, extend an element of \(\CC^\infty_c([0,T[\times\R^2)\) by zero past \(T\) and take
\(t=T\). To obtain
\eqref{eq:ds-fixed-time}, take \(\Phi(s,x)=\vartheta(s)\varphi(x)\)
with \(\vartheta\in\CC^\infty_c([0,\infty[)\) equal to one on
\([0,t]\). Then \(\p_s\Phi=0\) on that interval and the terminal pairing is the one in \eqref{eq:ds-fixed-time}. Done.  \end{proof}

\subsection{Renormalisation property}

The global classical weak solution to the 2D Euler equation constructed in this work renormalises passive scalars in $L^1\cap L^\infty$. We refer to the recent work Crippa--Nobili--Seis--Spirito~\cite{CNSS} for related definitions and background.

\begin{proposition}\label{prop:passive-renormalization}
Let $u = \kbs \star \omega$ be the global classical weak solution to the 2D Euler equations constructed in Theorem~\ref{mainthm}. On every finite interval $[a,b]$ with
$0\leq a<b$, the continuity equation with prescribed velocity $u$,
\[
  \partial_t\rho+\operatorname{div}(u\rho)=0,\qquad
  \rho(a,\cdot)=\rho_a\in L^1(\R^2)\cap L^\infty(\R^2),
\]
has a unique distributional solution in
$L^\infty(]a,b[;L^1(\R^2)\cap L^\infty(\R^2))$, and this solution is
renormalized and Lagrangian in the sense of \cite{CNSS}.
\end{proposition}

\begin{proof}[Proof of Proposition~\ref{prop:passive-renormalization}]
We verify conditions (4.1)--(4.3) in
\cite[Theorem~4.1 and Corollary~4.2]{CNSS}.

Write $M_b=\operatorname*{ess\,sup}_{0<t<b}\|\omega\|_{L^1_x}$.
For a radius $r>0$, split the Biot--Savart kernel $\kbs=K_{\leq r}+K_{>r}$ by restricting  to $\{|x|\leq r\}$ and $\{|x|>r\}$. The two pieces satisfy
\[
 \|K_{\leq r}\|_{L^1}=r,\qquad
 \|K_{>r}\|_{L^\infty}\leq \frac{1}{2\pi r}.
\]
If $M_b=0$, then $u=0$ on $]a,b[$. Otherwise, for a level $\lambda>0$
we take $r=\frac{M_b}{\pi\lambda}$. The far-field convolution has
$L^\infty$-norm at most $\lambda/2$, and the near-field convolution
has $L^1$-norm at most $rM_b$. Since $u=\kbs\star \omega$, the Markov or Chebyshev inequality yields that
\begin{align*}
 \leb^2\Big(\big\{x:|u(t,x)|>\lambda\big\}\Big)\leq \frac{2M_b^2}{\pi\lambda^2}
 \quad\text{for a.e. }t\in ]a,b[.
\end{align*}
Integrating in time leads to
$u\in L^{2,\infty}(]a,b[\times\R^2)$. Hence,  (4.1) holds with $p=2$.

Next, note that $\partial_i u_j=(\partial_i \kbs^j) \star \omega$ in the sense of distributions. To see that these are admissible kernels in (4.2) of
\cite[Theorem~4.1 and Corollary~4.2]{CNSS}, let $\mathcal{J}$ denote
the counterclockwise $\pi/2$-rotation, so that $x^\perp= \mathcal{J} x$. For $x\ne0$,
\[
 \partial_i \kbs^j(x)
 =\frac{1}{2\pi}\left(\frac{\mathcal{J}_{ji}}{|x|^2}
           -2\frac{(\mathcal{J}x)^j x^i}{|x|^4}\right),
\]
which is smooth on $\R^2\setminus\{0\}$, is homogeneous of degree $(-2)$, and its derivative $\lesssim |x|^{-3}$. Its integral over every annulus about the origin is zero, because
$\int_{\mathbb S^1}(\mathcal{J}\vartheta)^j\vartheta^i\,\dd \vartheta=\pi \mathcal{J}_{ji}$.
Also, observe that $\partial_i \kbs^j$ equals the principal value of the above off-origin kernel plus $\tfrac12 \mathcal{J}_{ji}\delta_0$ in the sense of distributions.  Its Fourier multiplier is
$\xi^i(\mathcal{J}\xi)^j/|\xi|^2$ for $\xi\ne0$, which is bounded. These verify the conditions (K1)--(K5) of \cite{CNSS}; see also their Section~5. Moreover, as $ \int_a^b\|\omega(t,\cdot)\|_{L^1_x}\,\dd t\leq (b-a)M_b$, the input to the singular integral belongs to the required spacetime $L^1$-class. This verifies (4.2).

Finally, (4.3) holds since
$\operatorname{div}u=0$.

The above arguments use the vorticity only for
$a<t<b$; in particular, no $L^1$-density for the initial datum at $t=0$ is required. We thus conclude the renormalisation property from \cite[Theorem~4.1 and Corollary~4.2]{CNSS}.    \end{proof}

Now let us prove Corollary~\ref{cor: bourbaki}:  we gave a negative answer to the following question quoted verbatim from De Lellis' S\'{e}minaire Bourbaki report~\cite{bourbaki}:

\begin{quote}
A very interesting open question, naturally linked to Euler's equations, is whether
the renormalization property holds for divergence–free fields $b \in L^\infty\big(\R; L^2(\R^2;\R^2)\big)$ when
the vorticity of $b$ is a measure. 
\end{quote}

In the proof below we shall never compose \(\beta\) with a measure, so the argument goes through. It is compatible with the passive scalar renormalisation in Proposition~\ref{prop:passive-renormalization}.

\begin{proof}[Proof of Corollary~\ref{cor: bourbaki}]

Suppose for contradiction that $\omega$ were renormalised on $[a,b]$ for any $0<a<b$ such that $u$ agrees with a $\CC^\infty$-Euler flow for any $t>100b$. Consider any $\beta \in \CC^\infty_b(\R;\R)$ with $\beta\equiv 0$ on $[-M,M]$, $\beta >0$ on $\R \setminus [-M,M]$, and $\sup_{z_1 \in \R}|z_1\beta'(z_1)| + \sup_{z_2 \in \R}|\beta'(z_2)|<\infty$, where $M> \|\omega(b,\cdot)\|_{L^\infty_x}$. It follows from the estimates of $\omega$ in Theorem~\ref{mainthm}  that \(\beta(\omega) \in L^\infty_t(L^1_x \cap L^\infty_x)\) on $[a,b]\times\R^2$. However, by the $\CC^0_tL^1_x$-continuity of vorticity in positive time, the trace of $\beta(\omega)$ at $t=b$ is zero. By a time reversal argument, assuming that $\omega$ is renormalised (hence uniqueness holds for $\beta(\omega)$), we may apply Proposition~\ref{prop:passive-renormalization} to \(-u(b - \cdot)\) to conclude that $\beta(\omega) \equiv 0$. Thus, $\|\omega(t,\cdot)\|_{L^\infty_x} \leq M$ uniformly in $a$. Sending $a \searrow 0$ to conclude that the initial vorticity $\omega_0$ (in the sense of a weak-star trace) must have a bounded density.  Contradiction.   \end{proof}

\subsection{Open problems}

Many central problems concerning the 2D Euler equation with the general vortex sheet initial data yet remain open. For instance:

\begin{itemize}
\item
Does the energy inequality hold for any global classical weak solution that has zero circulation?
\item
Can our solution constructed in Theorem~\ref{mainthm} arise in the vanishing viscosity limit for the Navier--Stokes equations? 
\end{itemize}

\medskip
\noindent
{\bf Acknowledgement}.
The author thanks Professors~Gui-Qiang Chen and Yifu Zhou for insightful discussions on Delort's problem and on convex integration.

 SL is supported by NSFC Projects 12201399, 12331008, and 12411530065, Young Elite Scientists Sponsorship Program by CAST 2023QNRC001, National Key Research $\&$ Development Programs 2023YFA1010900 and 2024YFA1014900, Shanghai Rising-Star Program 24QA2703600, Qi-Guang Scholarship, and Shanghai Frontier Research Institute for Modern Analysis.

	\noindent
	{\bf Statement of competing interests}. 
	The author declares that there is no conflict of interest.

	\noindent
	{\bf Statement of data availability}.
	Our manuscript has no associated data.

	\noindent
	{\bf AI Statement}. This paper has grown out of a series of interactions of the author with ChatGPT 6.1 and 6.2 sol/ultra AI models. When initialising this project, the author had the idea that the recent developments on convex integration, especially the work~\cite{BCK} by Bru\`{e}--Colombo--Kumar, may lead to a (somewhat nonphysical) resolution of the general vortex sheet problem. The AI models helped design the convex integration scheme and the temporal glueing construction in their final form. The author has independently verified the mathematical content and written up the whole paper, and assumes full responsibility for it. 

\end{document}